\documentclass[3p,times]{elsarticle}

\usepackage{mathtools}
\usepackage{amssymb}
\usepackage{amsfonts}

\usepackage{dsfont}
\usepackage{amsthm}
\usepackage{cases}
\usepackage{bm}
\usepackage{geometry}
\usepackage{indentfirst}
\usepackage{tabularx}
\usepackage{booktabs}
\usepackage{graphicx}
\usepackage{graphicx}
\usepackage{ulem}
\usepackage{xspace}
\usepackage[dvipsnames]{xcolor}
\usepackage{aliascnt}
\usepackage{multirow}
\usepackage{xspace}
\usepackage{dutchcal}
\usepackage{url}
\usepackage{accents}
\usepackage{bbm}
\usepackage{enumitem}
\usepackage{caption}
\usepackage{transparent}

\newcommand{\es}[1]{E{#1}}

\newcommand{\x}{\mathbf{x}}
\newcommand{\ct}{\mathbf{c}}
\newcommand{\Q}{\mathbf{Q}} 

\newcommand{\F}{\mathbf{F}}

\newcommand{\f}{\mathbf{f}}

\newcommand{\N}{\mathbb{N}}
\newcommand{\R}{\mathbb{R}}
\newcommand{\CFL}{\textrm{CFL}}

\newcommand{\mbf}[1]{\textbf{#1}}
\newcommand{\n}{\mathbf{n}}

\renewcommand{\u}{\mathbf{u}}
\newcommand{\q}{\mathbf{q}}

\newcommand{\de}[2]{\frac {\partial #1}{\partial#2}}

\usepackage[tracking=true]{microtype}

\newcommand{\be}{\begin{equation} \begin{aligned} }
		\newcommand{\ee}{\end{aligned} \end{equation}}

\usepackage{cancel}
\usepackage{soul}
\setstcolor{red}

\newcommand{\RIIcolor}[1]{{\leavevmode\color{black} #1}}
\newcommand{\RIcolor}[1]{{\leavevmode\color{black} #1}}
\newcommand{\Rall}[1]{{\color{black} #1}}

\journal{Computer Methods in Applied Mechanics and Engineering}

\begin{document}
	
	\begin{frontmatter}
		\title{On the treatment of topology changes on 3D polyhedral moving meshes\\ via 4D~space-time \textit{hole-like} elements in direct ALE ADER-DG methods}		 
		\author[univr1]{Elena Gaburro$^*$}
		\ead{elena.gaburro@univr.it}
		\cortext[cor1]{Corresponding author}
		
		\author[univr1]{Matej Klima}
		\ead{matej.klima@univr.it}
		
		\author[univr1]{Mauro Bonafini}
		\ead{mauro.bonafini@univr.it}
		
		\author[univr2]{Maurizio Tavelli}
		\ead{maurizio.tavelli@univr.it}

		\address[univr1]{Department of Computer Science, University of Verona, Strada le Grazie 15, Verona, 37134, Italy}		
		\address[univr2]{Department of Engineering for Innovation Medicine, University of Verona, Strada le Grazie 15, Verona, 37134, Italy}				
		
		\begin{abstract}  
			
			This work investigates a novel approach for the \textit{high order} evolution of hyperbolic PDEs using ADER discontinuous Galerkin schemes within a direct Arbitrary-Lagrangian-Eulerian (ALE) framework on 3D \textit{moving} polyhedral meshes with \textit{topology changes}. 
			\RIIcolor{ 
				Our \textit{direct} ALE method is based on the PDE integration over 4D (3D+time) space-time control volumes connecting the elements of two subsequent tessellations, so to \textit{simultaneously evolve} the solution both in time and between the two different meshes in an effective and high order manner. In this way, we also avoid any complex and expensive projection-reconstruction techniques and any mesh intersection operation typical of \textit{indirect} ALE schemes. 
				
				The crucial step consists in the strategy for building space-time control volumes that also connect elements with different shapes and neighborhoods due to a change in topology. In fact, simply linking existing elements by collapsing or expanding their edges would leave a "hole" in the space-time domain. To fill it, we introduce additional degenerate elements that we call \textit{hole-like} elements. These are 4D objects with zero 3D volume at both the beginning and end of the timestep, but which possess a strictly non-zero 4D space-time volume. 			
				Given the uniqueness of this space-time approach in 3D+time and the necessity of characterizing the geometry of such elements, the main objective of this paper is the formal geometrical and numerical description of the method. Specifically, we show that despite being non-trivial to visualize, these elements are \textit{well-defined} and ultimately even easy to manage.
			}
			In particular, we provide here a detailed characterization of the \textit{hole-like} elements needed to connect two polyhedral tessellations, corresponding to 2-3, 3-2, and 4-4 flips on the underlying Delaunay tetrahedralization. 	
			We describe how to identify and construct them, how to adapt the numerical method to their degenerate geometry by introducing for them a \textit{locally implicit} treatment, and we provide new and intuitive \textit{visualization} strategies. 	
			
			Finally, a set of benchmarks, simple in their wave structure but challenging for moving meshes, is used to prove that the method is fully conservative, satisfies the \RIIcolor{Geometric} Conservation Law condition, and maintains the correct order of convergence even in the presence of frequent topology changes. 
		\end{abstract}


		\begin{keyword}				
			ADER discontinuous Galerkin (DG) schemes\sep
			direct Arbitrary-Lagrangian-Eulerian (ALE) \sep			
			locally~implicit globally explicit methods \sep
			hyperbolic PDEs \sep
			3D moving polyhedral meshes \sep 
			topology changes \sep 
			space-time control volumes \sep 
			degenerate space-time geometry 
		\end{keyword}
	\end{frontmatter}

	\section{Introduction} \label{sec_intro}	
	
	In this work, we propose a novel algorithm designed to discretize any first-order 
	hyperbolic partial differential equation that can be cast into the following form
	\begin{equation} \label{eq.generalform}
		\partial_t{\Q}(\x, t) + \nabla \cdot \F(\Q(\x, t)) = \mathbf{0},
	\end{equation}
	by means of a high order accurate ADER discontinuous Galerkin method on 
	moving 3D meshes. 
	In the above formulation, $\x = (x_1,x_2,x_3) \in \R^3$ denotes the spatial 
	position vector, $t\ge0$ represents the time, $\Q \colon \Omega \times [0,+\infty) \to \R^\nu$ is the vector valued function representing the $\nu \in \N$ conserved variables. 
	Moreover, for given directional fluxes $\f_{i} \colon \R^\nu \to \R^\nu$, $i = 1,2,3$, we denote the flux function as $\F \colon \R^\nu \to [\R^{\nu}]^3, \Q \mapsto (\f_{1}(\Q), \f_{2}(\Q), \f_{3}(\Q))$ and interpret its Jacobian to be a vector of matrices in the form $\partial\F/\partial \Q = (\partial \f_{1}/\partial \Q, \partial \f_{2}/\partial \Q, \partial \f_{3}/\partial \Q)$. In particular, the spatial operator $\nabla = (\partial_{x_1}, \partial_{x_2}, \partial_{x_3})$ computes the divergence as $\nabla \cdot \F(\Q) = \partial_{x_1}\f_{1}(\Q) + \partial_{x_2}\f_{2}(\Q) + \partial_{x_3}\f_{3}(\Q)$.

	\paragraph{Arbitrary-Lagrangian-Eulerian methods}
	
	The use of Lagrangian methods, where the computational domain moves according to the local fluid velocity, dates back to the seminal works of von~Neumann, Richtmyer and Wilkins~\cite{Neumann1950, Wilkins1964}. These methods are particularly attractive for the simulation of non-linear hyperbolic conservation laws because, thanks to the flow-driven motion of the mesh, they allow to significantly reduce the numerical errors due to convection terms and to sharply capture moving interfaces and contact discontinuities~\cite{munz94, CaramanaShashkov1998, Maire2010, chengshu2, scovazzi2, despres2017numerical, morgan2021origins, cremonesi2010lagrangian, cremonesi2017explicit, Dobrev3, plessier2023implicit, del2023triangular, despres2024lagrangian,kincl2025numerical,KLIMA2024LAG,gaburro2026innovations}. Moreover, Lagrangian schemes exhibit negligible numerical dissipation at contact waves and material interfaces, and they can be made automatically entropy stable, rotationally invariant, and discretely Galilean invariant.
	
	However, Lagrangian methods in their pure form, where the mesh is forced to move exactly with the fluid velocity under a fixed connectivity, are commonly affected by mesh distortion or mesh tangling. These issues may slow down the computation or cause it to halt entirely, in particular if strong shear flows or vortical flows are simulated for long times. To overcome these limitations, the Arbitrary-Lagrangian-Eulerian (ALE) technique was developed, starting from the works~\cite{wilkins1964methods, trulio1966air, hirt1974arbitrary} in the context of finite difference discretizations, and later adopted in the finite element community~\cite{donea1977lagrangian, hughes1981lagrangian, belytschko1978computer, donea1982arbitrary}. An extensive review of the original ALE literature can be found in~\cite{benson1992computational}.
	
	ALE techniques can generally be distinguished into two categories. The first is that of \textit{indirect} ALE schemes, characterized by a rezoning procedure and a subsequent remapping phase, where the numerical solution is transferred from the old Lagrangian mesh onto the new optimized grid~\cite{LoubereShashkov2004,ReALE2010, ReALE2011, ReALE2015, ShashkovCellCentered, ShashkovRemap1, MaireMM2, indALE-AWE2016, wu2021cell, kenamond2021positivity,lei2023high,KLIMA2020ALE,BARLOW2018ALE}. The second category is that of \textit{direct} ALE schemes, where the discretization method responsible for the PDE evolution directly provides the updated solution on the new mesh, without the need of a projection-reconstruction procedure, but just via the PDE integration over space-time control volumes~\cite{ALELTS1D,Lagrange2D,Lagrange3D,ALELTS2D, springel2010pur,gaburro2017direct}. 
	
	Note that, a fundamental tool of modern ALE methods, especially required when complex flows lead to extreme mesh deformations,
	consists in allowing \textit{topology} and neighborhood changes during the mesh motion, as seen in the indirect ReALE methods~\cite{ReALE2010, ReALE2015,KUCHARIK2008CONNECTIVITY}, and in the context of the AREPO code~\cite{springel2010pur}.
	Here, because of the complexity of developing high order accurate remapping techniques, the direct ALE framework represents a convenient option when interested in high order accurate numerical methods, such as those based on Finite Volume (FV) or Discontinuous~Galerkin (DG)~schemes.
	
	In particular, FV and DG schemes combined with the ADER approach~\cite{mill, toro1, toro3, dumbser2008unified} are particularly well-suited for the direct ALE framework thanks to their intrinsic space-time formulation. 
	ADER FV and DG methods have been extensively adopted in recent literature on Eulerian meshes for a wide variety of applications, 
	starting from the classical Euler and shallow water equations~\cite{rannabauer2018ader,fernandez2022arbitrary,boscheriAFE2022,micalizzi2023efficient,rio2024high,ciallella2024very},
	multiphase models \cite{kemm2020simple,rio2024high,gaburro2024discontinuous}, and magnetohydrodynamics \cite{balsara2012self,fambri2020discontinuous}, 
	including also the treatment of dispersive, turbulent and reactive flows \cite{yuan2022hybrid,busto2022new,popov2023space}, 
	up to the study of much more complex systems, such as unified models for continuum mechanics \cite{dumbser2016high,dumbser2017high,tavelli2020space,busto2020high,chiocchetti2021high,lakiss2024ader}
	and various first-order reformulations of the Einstein field equations for general relativity \cite{dumbser2019glm,dumbser2023WBGR,zanotti2025new,muzzolon2026high},
	thus demonstrating their broad interest and capabilities.
	
	Notably, ADER FV and DG schemes have been successfully extended also to the direct ALE framework on moving Voronoi meshes with topology changes in~\cite{gaburro2020high, gaburro2021unified, gaburro2021bookchapter,gaburro2025high}, where the key ingredient is the PDE integration over space-time control volumes, possibly degenerate, connecting even different tessellations.	
	These methods are conservative, satisfy the Geometric Conservation Law (GCL), and achieve high order of accuracy. Moreover, their stability, also in the presence of degenerate space-time geometries, has been recently proved in~\cite{bonafini2026stability}.
	For the sake of completeness, we would also like to acknowledge the pioneering work of Alauzet \textit{et al.} \cite{olivier2011new} who, already in 2011, albeit in a different context, proposed the intuition of using a degenerate control volume to connect 2D polygonal meshes undergoing topology changes.
	
	Finally, we remark that many of the works cited above were validated just in the two-dimensional framework. Specific remarkable developments in the 3D ALE context have been recently achieved in~\cite{anderson2018high,di20243d,guisset2024cell,vargas2025multi,breil20253d,colaitis2026cell-1,colaitis2026cell-2,eder20263d}, all featuring the indirect approach.
	
	\paragraph{Mesh generation and optimization}
	
	Fundamental pillars of direct and indirect ALE methods concern i)~generating a simulation-ready mesh, and then applying optimization techniques to ii)~correct element movement and iii)~update element connectivity.
	In our case, as explained in Section~\ref{sec.domain}, we use a polyhedral mesh for the computations, 
	but we generate it as the dual of a tetrahedralization to take advantage of the more extensive literature available for the latter. Here, we recall the works that have inspired our strategies.
	
	In~\cite{freitagGooch}, an empirical study was performed showing that the combination of mesh smoothing and transformations, specifically edge/face removal operations, is a viable technique for tetrahedral mesh optimization. That study was subsequently extended in~\cite{klingerSchewchuk} with a more detailed look at possible transformations and the inclusion of additional features such as vertex addition and removal. Nowadays, software tools based on these methods, such as the TetGen library~\cite{hang2015tetgen}, are readily available and can be used for efficient tetrahedral mesh generation and optimization.
	
	These techniques are generally global and lead to a complete transformation of the initial mesh. 
	This is not suitable for the context of our work, which requires a continuous and incremental optimization, controlled and localized, performed at each timestep. To our knowledge, this approach has only been explored in 2D in~\cite{olivier2011new, gaburro2020high}. In 3D, this methodology, which is part of the novelty of our work, remains largely under-explored, with the exception of an incremental tetrahedral mesh optimization proposed in~\cite{alauzet2014changing}. 
	
	Therefore, we have opted to implement a step-by-step version of the more traditional edge removal procedure. Edge removal operations in tetrahedral meshes are well-defined and were already described several decades ago in~\cite{delisleGeorge}; they are currently employed in state-of-the-art mesh optimizers such as~\cite{hang2015tetgen,DASSI20182}. In particular, in our present work, we build upon the algorithm described by Dassi \textit{et al.} in~\cite{DASSI20182}.

	\paragraph{Challenges and novelties of the present work and structure of the paper}
	
	The objective of this work is the development of a 3D \textit{direct} ALE ADER-DG scheme on polyhedral meshes, where the mesh velocity is designed to follow the fluid flow, as closely as possible, but can be modified to allow for mesh optimization techniques driven by the quality of the underlying Delaunay tetrahedralization (as explained in Section~\ref{sec.domain}).
	To grant the optimization process full flexibility, the mesh motion is characterized by frequent topology changes. 
	In this context, our goal is to propose a method that remains high order accurate, fully conservative, and satisfies the GCL even in the vicinity of topology changes, without the use of traditional but expensive projection-reconstruction techniques. 
	Instead, we connect subsequent meshes through space-time control volumes, yet necessitating also the inclusion of degenerate elements (as described in Section~\ref{sec.space-time}).
	
	While the conceptual foundations of using \textit{degenerate} space-time elements (formerly denoted as \textit{slivers}) to handle topology changes were established in previous two-dimensional works of the authors~\cite{gaburro2020high, gaburro2025high, bonafini2026stability}, the extension to 3D+time entails major challenges.
	First, it requires managing a complex, dynamic data structure that must track the non-trivial impact of each Delaunay flip on the dual polyhedral tessellation at each timestep (see Section~\ref{s-sec.optimization-flips}). Second, we need to formally characterize and build these new 3D+time=4D \textit{hole-like} elements so to completely fill the space-time domain, while also proposing strategies to visualize them, see in particular Section~\ref{s-sec-topology-flips}.
	
	Then, we must address the redesign of the numerical method (described in Section~\ref{sec.method-description}), specifically in the neighborhood of these degenerate elements. 
	In particular, here, we introduce for the first time a novel \textit{locally implicit} treatment of our \textit{hole-like} elements, see Section~\ref{s-sec.predictorcrazy}.
	This strategy represents a significant improvement of our paradigm: by coupling the \textit{hole-like} evolution with the so-called \textit{predictors} approximating the solution in the neighboring elements, we automatically achieve full conservation on topology changes without any \textit{a posteriori} correction, while strictly preserving the globally explicit character of the overall method.
	
	Finally, a set of benchmarks is presented in Section~\ref{sec_num-results} to demonstrate that the method is fully conservative, satisfies the GCL, and maintains the correct convergence order even during frequent topology changes. 
	The paper ends with concluding remarks and an outlook to future work in Section~\ref{sec_concl}.
	

	
	\section{Spatial discretization of the domain and mesh movement}
	\label{sec.domain}
	
	\subsection{An optimized tetrahedralization and a computational polyhedral tessellation} 
	\label{s-sec.spatialdomain}
	
\begin{figure}[!b]
		\centering
		\includegraphics[width=0.97\linewidth]{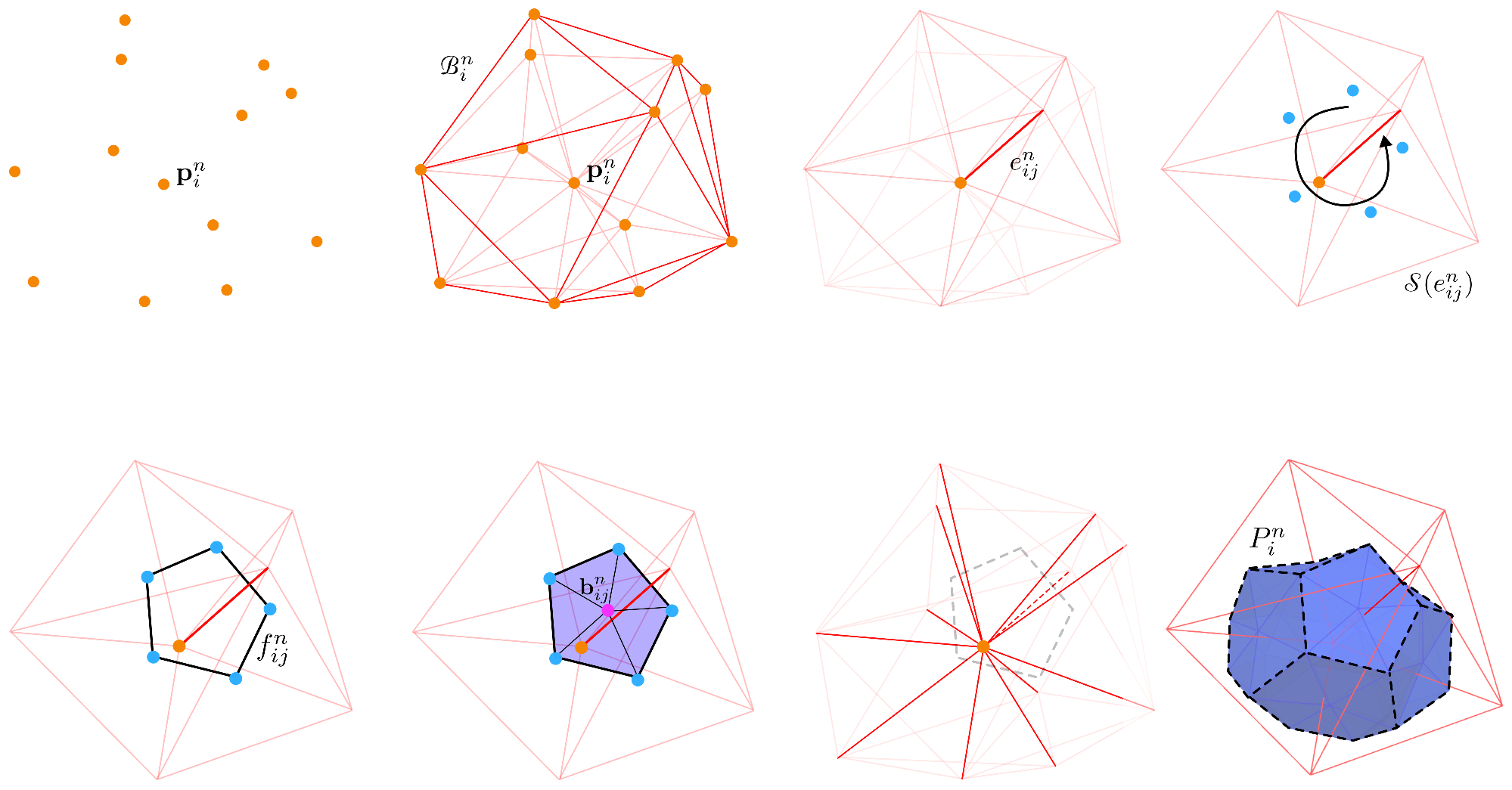}%
		\caption{ 
			From left to right, top to bottom.
			We cover our domain with a set of generator points (orange) and we connect them via a tetrahedralization.
			We then fix a generator $\textbf{p}_i^n$ along with the set of its surrounding tetrahedra $\mathcal{B}_{i}^n$. Then, we select one of the edges $e_{ij}^n$ that are connected to $\textbf{p}_i^n$ and construct the star region $\mathcal{S}(e_{ij}^n)$ consisting of the tetrahedra sharing this edge.
			We consider the (generally not co-planar) centroids of the tetrahedra in $\mathcal{S}(e_{ij}^n)$ and we connect them orderly. 
			We compute the barycenter $\textbf{b}_{ij}^n$ (pink) of these points and we use it to build a 2D surface $f_{ij}^n$ made of triangular planar facets connecting each boundary segment with the barycenter. 
			By repeating the process over each outgoing edge, we obtain our final polyhedral computational element $P_i^n$ around the selected generator $\textbf{p}_i^n$, with $|\mathcal{V}_i^n|$ neighbors $\{P_j^n\}_{j\in \mathcal{V}_i^n}$ and lateral surfaces $\mathcal{F}_i^n = \{f_{ij}^n\}_{j \in \mathcal{V}_i^n}$.
		}
		\label{fig.mesh}
	\end{figure}

	{
		\linespread{1.1}\selectfont
		At each time level $t^n$ we cover our 3D computational domain $\Omega^n$ with two different meshes, 
		a tetrahedralization~$\mathcal{T}^n$ and a polyhedral tessellation~$\mathcal{P}^n$, 
		each of them with a different purpose. 
		We detail their construction here below and we also refer the reader to Figure~\ref{fig.mesh} 
		and its caption for a visual interpretation and a concise description. \par
	}
	
	{
		\linespread{1.1}\selectfont
		We start by populating $\Omega^0$ at time $t=0$ with $N_P$ points $\textbf{p}_i^n$, $i=1,\dots, N_P$,
		that we call \textit{generator} points (the orange points in Figure~\ref{fig.mesh}). 
		The number of generators $N_P$ is fixed all along a given simulation, while their position may vary in time, as described later in Section~\ref{s-sec.mesh_movement}.
		At the beginning of each timestep $t^n$, we first connect the generators via a tetrahedralization $\mathcal{T}^n$
		made of $N_T^n$ tetrahedra ${T}_i^n$, $i = 1, \dots, N_T^n$.	
		The set of indices of all tetrahedra connected to a generator point $\mathbf{p}_i^n$, which form a sort of ball around the point, is denoted as $\mathcal{B}_i^n$.
		The tetrahedralization is generated to satisfy as much as possible a set of desired quality criteria as described in Section~\ref{s-sec.optimization-criteria}.
		In particular, at time $t=0$ we build it by directly employing the software TetGen~\cite{hang2015tetgen}, 
		whereas at later timesteps we deduce $\mathcal{T}^{n+1}$ from $\mathcal{T}^n$ by using a set of in-house routines, inspired both to those of TetGen and to the strategy forwarded in~\cite{DASSI20182}, as explained in~\ref{s-sec.optimization-flips}.
		In this way, we can precisely control the connectivity changes (flips) 
		that occur between one timestep and the next one.	
		The \textit{tetrahedralizations} $\mathcal{T}^n$ are employed to ensure the \textit{quality} of the discretization and to oversee the connectivity transformations eventually performed.\par
	}
	
	Then, we deduce from $\mathcal{T}^n$ a (dual) \textit{polyhedral tessellation} $\mathcal{P}^n$,
	which is the mesh that will be used for \textit{computational} purposes. 
	Indeed, in contrast to tetrahedra, polyhedral cells offer several numerical advantages~\cite{boscheriAFE2022} in the context of cell-centered discretizations: 
	they allow for larger timestep size, exhibit significantly less mesh imprinting on the solution 
	and, thanks to the higher number of neighbors, provide a superior capability to capture multidimensional phenomena even when using multidimensional Riemann solvers~\cite{gallice2022entropy,Gaburro2026}. 
	Moreover, within ALE methods with topology changes, they are usually the preferred choice~\cite{ReALE2010,ReALE2015,springel2010pur} because they are more robust and allow for a more physically-driven information flow between corresponding (but different) cells.
	
	{
		\linespread{1.1}\selectfont		
		Each polyhedron $P_i^n \in \mathcal{P}^n$, $i=1,\dots, N_P$, is built around the corresponding generator $\textbf{p}_i^n$ and
		it has as many neighbors as the edges outgoing from $\textbf{p}_i^n$.	
		In particular, with $\mathcal{V}_i^n$ we denote the set of indices of the neighbors $P_j^n$ of~$P_i^n$. 
		To build $P_i^n$, we loop over the tetrahedralization edges $\{e_{ij}^n\}_{j \in \mathcal{V}_i^n}$ outgoing from $\textbf{p}_i^n$ and, 
		for each fixed outgoing edge~$e_{ij}^n$, we consider the set 
		$\mathcal{S}(e_{ij}^n)$ of tetrahedra sharing it, also termed the \textit{star region} of the edge $e_{ij}^n$.
		We then compute the \textit{centroids} of the tetrahedra in the star region and connect them orderly to build the countour of the 2D surface $f_{ij}^n$ of~$P_i^n$. 
		Note that we consider each of these surfaces $f_{ij}^n$ as being made of triangular planar facets connecting each two consecutive points of its contour with the barycenter of the face $\mathbf{b}_{ij}^n$. 
		By repeating the process over each outgoing edge, we obtain our polyhedral element $P_i^n$ as the region enclosed within the set of faces $\mathcal{F}_i^n = \{f_{ij}^n\}_{j \in \mathcal{V}_i^n}$.
		We call $\textbf{c}_i^n$ the center of mass of each element $P_i^n$ (not necessarily coinciding with $\textbf{p}_i^n$). 
		Moreover, we define $\mathcal{D}_i^n = \{ \textbf{d}_{ij}^n \}_j$ as the set of all the vertices of $P_i^n$, 
		that is the set of all the \textit{centroids} of $\mathcal{B}_i^n$. 
		Finally, for the sole purpose of selecting quadrature nodes for volume integral calculations with optimal accuracy and high flexibility, we often consider $P_i^n$ as a union of tetrahedra, each obtained by connecting a triangular facet of $\mathcal{F}_i^n$ with $\mathbf{c}_i^n$.
		For a graphical representation of $\mathbf{c}_i^n$ and the sub-tetrahedra, see also Figure~\ref{fig.controlvolumes}.	 
		\par 
	}

	Note that classically named Voronoi cells are built in this same way starting from a given 
	tetrahedralization, but taking the \textit{circumcenters} of the tetrahedra. 
	The advantage of Voronoi cells is that they have planar surfaces that would not need to be further sub-triangulated; 
	however, they offer less control over the cell shape.
	Thus, in order to have rounder cells, whose quality is more directly linked to the 
	tetrahedra quality, we have made the choice of building a so-called \textit{centroid-based}
	Voronoi-like tessellation. This choice actually makes our code more general: 
	the classical Voronoi case is simply a particular instance of our framework and could be handled automatically.

	\subsection{Generator dynamics and updated tetrahedralization connectivity}
	\label{s-sec.mesh_movement}

	{
		\linespread{1.05}\selectfont			
		The position of the generators evolves at each timestep from $ \mathbf{p}_i^n$ to a novel value $\mathbf{p}_i^{n+1}$ 
		determined in order to balance \textit{two} often contrasting \textit{requirements}: 
		maintaining a Lagrangian approach while ensuring high mesh quality. 
		Specifically, our motion is obtained by:
		i) a \textit{high order} integration of the generator \textit{trajectories} to closely follow the fluid flow (see Section~\ref{ssec.HighOrderTraj}), 
		followed by ii) an \textit{optimization} of their position so that the mesh, to be built at the next timestep, 
		$\mathcal{T}^{n+1}$ features high-quality elements, minimizing potential numerical errors due to mesh distortion (see Section~\ref{s-sec.smoothing}).
		
		Once the new positions of the generators $\mathbf{p}_i^{n+1}$ have been fixed, we can construct $\mathcal{T}^{n+1}$. This is achieved by performing at most one elementary flip per generator with respect to $\mathcal{T}^{n}$, as described in Section~\ref{s-sec.optimization-flips}.
		The purpose of this operation is to determine a connectivity $\mathcal{T}^{n+1}$ for the points $\mathbf{p}_i^{n+1}$ that offers better mesh quality than simply retaining the topology of $\mathcal{T}^{n}$.
		By limiting the number of flips per region (a constraint that could be relaxed in future developments), 
		we maintain a manageable link between the two successive configurations. 
		This allows for a completely new way of connecting the two tessellations in space-time (see in particular Sections~\ref{s-sec-topology-flips}~and~\ref{s-sec.predictorcrazy}) which, while requiring a certain conceptual effort, remains manageable compared to the complexity of standard projection-reconstruction procedures.\par
	}

	\subsubsection{High order integration of the generator trajectories}
	\label{ssec.HighOrderTraj}
	
	Due to the ALE framework of the present work, the mesh can in principle be moved with an \textit{arbitrary} velocity $\mathbf{v}_i^n$
	\begin{equation*}
		\mathbf{p}_i^{n+1} = \mathbf{p}_i^n + \Delta t \, \mathbf{v}_i^n, \ \text{ with }\ \Delta t \text{ given by}\ \eqref{eq:timestep}, 
		\label{eqn.xpnew} 
	\end{equation*}
	and there is not a specific necessity of moving the grid in a fully Lagrangian fashion.
	Nevertheless, it is suitable to compute accurately this velocity to better exploit the Lagrangian benefits of the scheme; to do this we employ a high order accurate approximation of the generator trajectories obtained by a high order Taylor method~\cite{boscheri2013semi, tavellihigh, gaburro2020high}. 	
	The use of this technique improves the overall Lagrangian behavior of the algorithm and in vortical flows also the mesh quality.
	
	The Taylor expansion of the new (temporary) generator position $\hat{\mathbf{p}}_i^{n+1}$ 
	with respect to its position at time $t^n$ can be written as 
	\begin{equation}
		\hat{\mathbf{p}}_i^{n+1} = 	\mathbf{p}_i^{n} + \Delta t \frac{d \x}{dt}  + \frac{\Delta t^2}{2} \frac{d^2\x}{dt^2} 
		+ \frac{\Delta t^3}{6} \frac{d^3\x}{dt^3}  + \frac{\Delta t^4}{24} \frac{d^4\x}{dt^4} + \mathcal{O}(5),
		\label{eqn.xcnew_high} 
	\end{equation}
	which achieves \textit{fourth} order of accuracy in time.
	Now, the high order time derivatives in~\eqref{eqn.xcnew_high} are replaced by high order spatial derivatives, via the \textit{Cauchy-Kovalevskaya} procedure, using repeatedly the trajectory equation 
	\begin{equation*}
		\frac{d \x}{dt} = \textbf{v} (\textbf{x}(t)),
	\end{equation*}
	and assuming a stationary velocity field (i.e. $\partial_t \textbf{v} = 0$). Hence, by denoting a vector $\mathbf x = (x_1, x_2, x_3)$ as $[x_i]_i$, we compute
	\begin{equation*}
		\dfrac{d\mbf{x}}{dt}     = \mbf{v} = [\varv_i]_i,\qquad
		\dfrac{d^2\mbf{x}}{dt^2} = \frac{d}{dt}\left(\dfrac{d\mbf{x}}{dt}\right) = \de{\mbf{v}}{\mbf{x}}\,\de{\mbf{x}}{t} 
		= \left[ \sum_j \frac{\partial \varv_i}{\partial x_j}\,{\varv_j}\right]_i.
	\end{equation*}
	By applying the chain rule again, the third derivative of the position reads
	\begin{equation*}
		\dfrac{d^3\mbf{x}}{dt^3} = \frac{d}{dt}\left(\dfrac{d^2\mbf{x}}{dt^2}\right) = 
		\left[ \sum_{j,k} \left( \frac{\partial^2 \varv_i}{\partial{x_j}\,\partial{x_k}}\,\varv_j\,\varv_k + \frac{\partial{\varv_i}}{\partial{x_j}}\,\frac{\partial{\varv_j}}{\partial{x_k}}\,\varv_k \right) \right]_i,
	\end{equation*}
	and similarly, the fourth derivative reads
	\begin{equation*}
		\dfrac{d^4\mbf{x}}{dt^4} =
		\left[\sum_{j,k,l} \left( \frac{\partial^3 \varv_i}{\partial x_j\,\partial x_k\, \partial x_\ell}\,\varv_j\,\varv_k\,\varv_\ell + 
		\frac{\partial^2 \varv_i}{\partial x_j\,\partial x_k}\,\frac{\partial \varv_k}{\partial x_\ell}\,\varv_\ell\,\varv_j + 
		2\,\frac{\partial^2 \varv_i}{\partial x_j\,\partial x_k}\,\varv_k\,\frac{\partial{\varv_j}}{\partial x_\ell}\,\varv_\ell + 
		\frac{\partial \varv_i}{\partial x_j}\,\frac{\partial^2 \varv_j}{\partial x_k\,\partial x_\ell}\,\varv_\ell\,\varv_k + 
		\frac{\partial \varv_i}{\partial x_j}\,\frac{\partial \varv_j}{\partial x_k}\,\frac{\partial \varv_k}{\partial x_\ell}\,\varv_\ell \right)\right]_i.
	\end{equation*}
	Finally, the partial derivatives of $\textbf{v}$ evaluated at each $\textbf{p}_i^n$ are recovered from the high order polynomial representation~$\u_i^n$~\eqref{eq.un} of the conserved variables $\mathbf{Q}$ inside each cell $P_i^n$ and thus from the local fluid velocity.

	\subsubsection{Mesh smoothing}
	\label{s-sec.smoothing}
	
	Since in a ALE framework the velocity is not constrained to follow the local fluid velocity \textit{exactly}, we have the possibility to apply some mesh optimization techniques to improve the quality of the moving mesh.   
	In this work, the mesh regularization procedures are implemented in such a way to slightly modify, 
	at each timestep, the Lagrangian coordinates of the generator points, that is, the vertices of the tetrahedralization, so to optimize the quality of the novel tetrahedralization $\mathcal{T}^{n+1}$ to be built.

	Thus, we start by computing $\hat{\textbf{p}}_i^{n+1}$ for all $i\in \{0, \dots, N_P\}$ according to the previous Section~\ref{ssec.HighOrderTraj}.
	Then, we apply to the new Lagrangian coordinates a smoothing technique to recover, for each generator, an alternative candidate position $\accentset{\star}{\mathbf{p}}_i^{n+1}$ (see a few paragraphs below).	
	We say that $\accentset{\star}{\mathbf{p}}_i^{n+1}$ is a location for the generator that is optimal 
	in the sense of mesh quality, as opposed to optimal
	in following the flow of the fluid, which is the role of $\hat{\textbf{p}}_i^{n+1}$.
	The final novel position of the generator is then recovered as 
	\begin{equation} 
		\label{eq.blending-use}
		{\textbf{p}}_i^{n+1} = (1 - \mu)\,\hat{\textbf{p}}_i^{n+1} + \mu\,\accentset{\star}{\mathbf{p}}_i^{n+1},
	\end{equation} 
	with $\mu$ being a blending factor that yields the balance between 
	the amount of mesh motion due to fluid flow with the one due to smoothing. 
	
	To determine this \textit{blending factor} we take into account the fluid velocity, the current grid configuration, as well as the specific explicit timestep restriction in use, following the approach described in detail in~\cite{gaburro2020high}. 
	Specifically, we compute 	
	\begin{equation} \label{eq:relaxationparameter}
		\mu = \min \left( 1, \sqrt{ \frac{ U_\ast \, \Delta t }{ \displaystyle \min_{i=1, \dots, N_P} h_i^n } \, \kappa } \, \right),
	\end{equation}	
	with ${U}_\ast$ being a rough scaling estimate for the mesh velocity, obtained at each timestep as the maximum \Rall{non uniform} velocity encountered for all generator points, $\Delta t$ the timestep size and $h_i^n$ the characteristic length of a mesh cell (defined in~\eqref{eq.char.length}). The underlying idea is that we want to balance, during each timestep, the spatial scaling of fluid flow, with a characteristic length representative of the mesh motion due to pure smoothing in the smallest cells of the domain, which we implicitly assume to be the most delicate.
	In this way, we have a blending factor $\mu$ which is a non-dimensional \textit{parameter} and that can be modulated through the factor~$\kappa$ which we set to a small value, e.g. in this work \Rall{we mostly use~$\kappa = 1/100$}.
	In general, the factor~$\kappa$ may depend on the necessity of the test-case, namely on the desired amount of Lagrangian flavor we wish for the final results, but it is independent on the mesh size and the order of the method, being these features already automatically handled by the formula~\eqref{eq:relaxationparameter}.

	Concerning the computation of $\accentset{\star}{\mathbf{p}}_i^{n+1}$, we mainly follow the vertex smoothing technique presented by Alauzet in~\cite{alauzet2014changing}, that we summarize here for the sake of completeness.
	The underlying idea consists in optimizing the quality of the current tetrahedralization $\mathcal{T}^n$ used to connect, instead of the old coordinates $\mathbf{p}_i^n$, the newly calculated Lagrangian coordinates~$\hat{\mathbf{p}}_i^{n+1}$.
	In order to find $\accentset{\star}{\mathbf{p}}_i^{n+1}$, 
	the first step requires to compute, for each tetrahedron $T_j^n \in \mathcal{B}_i^n$, 
	$\accentset{\star}{\mathbf{p}}_{i,j}^{n+1}$ which represents the ideal coordinate of the point in order to maximize just the quality of the tetrahedron $T_{j}^n$.
	Then, $\accentset{\star}{\mathbf{p}}_i^{n+1}$ is obtained as a weighted average of the these optimal point locations 
	\begin{equation} 
		\accentset{\star}{\mathbf{p}}_i^{n+1} = \frac{\sum_{j \in \mathcal{B}_i^n} \max{\left(Q\left(\hat{{T}}_j^{n+1}\right), Q_{\max}\right)} \; \accentset{\star}{\mathbf{p}}_{i,j}^{n+1}}{\sum_{j \in \mathcal{B}_i^n} \max{\left(Q\left(\hat{{T}}_j^{n+1}\right), Q_{\max}\right)}}, 
		\label{eqn.smooth_param}
	\end{equation} 
	where $Q_{\max} = 100$ is a cutoff parameter and 
	$Q(\hat{{T}}_j^{n+1})$ is a quality coefficient computed for the tetrahedron $\hat{{T}}_j^{n+1}$, that is a tetrahedron retaining the same connectivity of ${T}_j^n$ but whose vertex coordinates are $\hat{\mathbf{p}}_i^{n+1}$, i.e., those moved with Lagrangian velocity.
	The quality coefficient is defined as
	\begin{equation*} 
		Q(\hat{{T}}_j^{n+1}) = \frac{\sqrt{3}}{216} \frac{\sum\limits_{k,l \in \hat{{T}}_j^{n+1}, \, k \neq l } \left\| \hat{\mathbf{p}}_k^{n+1} - \hat{\mathbf{p}}_l^{n+1} \right\|^2}{\left|\hat{{T}}_j^{n+1}\right|},
	\end{equation*} 
	where $|\hat{{T}}_j^{n+1}|$ is the volume of $\hat{{T}}_j^{n+1}$. This coefficients approaches $1$ for a perfectly regular tetrahedron and grows for irregular elements.

	\subsubsection{Incremental tetrahedralization optimization via elementary flips performed over multiple timesteps} 
	\label{s-sec.optimization-flips}

	At this point, the new coordinates $\mathbf{p}_i^{n+1}$ have been fixed through~\eqref{eq.blending-use}, and now our goal is to find a new way to connect them through a tetrahedralization $\mathcal{T}^{n+1}$ that is deduced from $\mathcal{T}^n$ and improves its quality.
	
	Unlike classical mesh optimization methods, which operate globally on the entire mesh, 
	repeatedly rearranging the connections between all generators
	until the desired quality is achieved, 
	often resulting in a complete transformation where the connection between initial and final elements is completely lost, 
	we focus here on a fundamentally different approach.
	In particular, we aim for a \textit{continuous} and \textit{incremental} mesh improvement performed via a \textit{controlled number} of connectivity changes per region at each timestep. If this limited amount of transformation is insufficient to reach the \textit{target quality} within a single step, we rely on the fact that the optimization will persist across subsequent timesteps, reaching the target quality incrementally over time \RIIcolor{(this feature will be shown through numerical evidence in our benchmarks)}.
	Furthermore, since we start from a high-quality initial mesh and apply our optimization strategy at every stage, the mesh quality remains consistently under control; consequently, the minimal interventions to which we limit ourselves in this approach appear to be sufficient for our current purposes.
	
	\paragraph{Elementary flips} 
	First of all, we impose that the total number of generators $N_P$ do not vary all along the simulation:
	the Delaunay transformations that preserve this number can be 
	all expressed as a sequence of elementary operations called~2-3, 3-2 flips and potentially 4-4- flips.
	Thus, without loss of generality, since any complex mesh transformation can be achieved through a sequence of such flips, we limit ourselves to \textit{performing only elementary flips}, whose full sequence required to optimize a specific region will be completed across multiple successive timesteps.

	\begin{figure}[bt]
		\centering
		\includegraphics[width=1.0\linewidth]{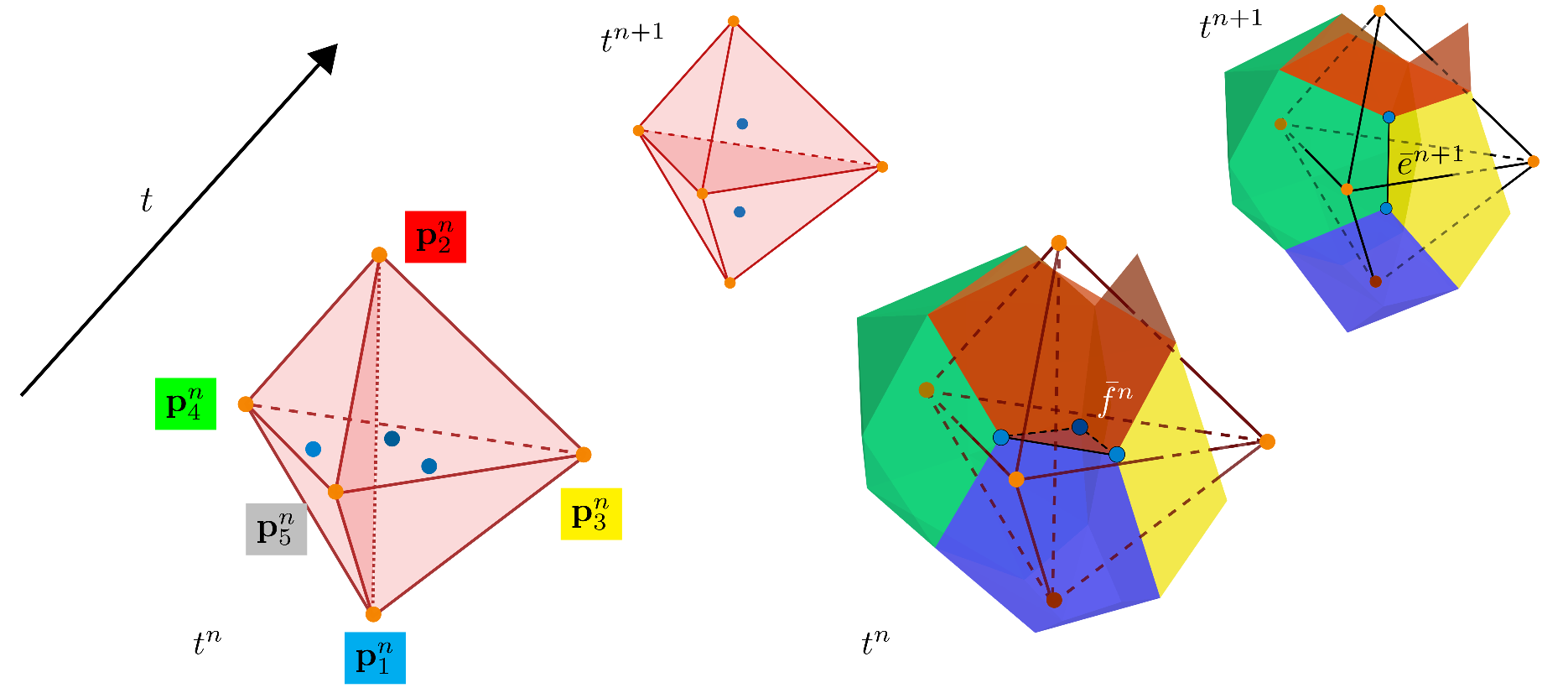}		
		\caption{	
			In this figure, we depict the portion of tetrahedralization (left) where a \textit{3-2 flip} occurs between time $t^n$ and $t^{n+1}$ and the corresponding effect on the polyhedral tessellation (right). In particular, we show 5 generators at time $t^n$ (bottom-left) and at time $t^{n+1}$ (top-left), and 4 of the 5 corresponding polyhedral elements at time $t^n$ (bottom-right) and time $t^{n+1}$ (top-right). 
			At time $t^n$, the five generators are connected via 3 tetrahedra sharing the central edge $\left\{\mathbf{p}_1^n, \mathbf{p}_2^n\right\}$, while at time $t^{n+1}$, they are connected by only 2 tetrahedra sharing the central face, from which the name 3-2 flip. 
			This flip also modifies the topology of the corresponding polyhedra, specifically: 
			i) the blue and the red polyhedron, at time $t^n$, are neighbors sharing a triangular face $\bar{f}^n$; conversely, at time $t^{n+1}$, this face disappears and they no longer touch; moreover, a new edge is created $\bar{e}^{n+1}$ which can be used as link between the blue and the red element; ii) the yellow, green and gray elements share the edge $\bar{e}^{n+1}$ at time $t^{n+1}$, while at time $t^n$ they are not in direct contact; 
			their connection at $t^n$ is represented by the fact that they all have as edge one boundary of the face $\bar{f}^n$. 
			This is the situation where a \textit{hole-like} element should be constructed in order to establish a space-time connection between $\bar{f}^n$ and $\bar{e}^{n+1}$ and avoid leaving a gap in space-time.
		}
		\label{fig.3-2flip-origin}
	\end{figure}	
	The configuration of a single \textit{3-2 flip} on the tetrahedral mesh is depicted on the left of Figure~\ref{fig.3-2flip-origin} 
	and it consists of the transformation from a set of \textit{three} tetrahedra, 
	defined by the generator points
	\begin{equation*}
		\left\{\left\{\mathbf{p}_1^n, \mathbf{p}_2^n, \mathbf{p}_3^n, \mathbf{p}_4^n \right\}, \left\{\mathbf{p}_1^n, \mathbf{p}_2^n, \mathbf{p}_4^n, \mathbf{p}_5^n \right\}, \left\{\mathbf{p}_1^n, \mathbf{p}_2^n, \mathbf{p}_5^n, \mathbf{p}_3^n \right\} \right\},
	\end{equation*}
	to a set of just \textit{two} tetrahedra 
	\begin{equation*}
		\left\{\left\{\mathbf{p}_1^{n+1}, \mathbf{p}_3^{n+1}, \mathbf{p}_4^{n+1}, \mathbf{p}_5^{n+1} \right\}, \left\{\mathbf{p}_2^{n+1}, \mathbf{p}_5^{n+1}, \mathbf{p}_4^{n+1}, \mathbf{p}_3^{n+1} \right\}\right\},
	\end{equation*}
	which are obtained by removing the edge $\left\{\mathbf{p}_1^n, \mathbf{p}_2^n\right\}$ from the configuration. 
	The resulting effect on the polyhedral tessellation, depicted on the right of Figure~\ref{fig.3-2flip-origin}, 
	will be described in details in Section~\ref{s-sec.elementaryflip3-2}.
	
	In addition, we note that the so-called \textit{2-3 flip} corresponds to the same operation carried out in the reverse order. 
	We also remark that the 3-2 flip is a special case of \textit{edge removal} which, however, is an operation not applicable to all connections between the generator points because an edge can have more neighboring tetrahedra than just three, preventing its removal by a single elementary operation.
	
	\begin{figure}[bt]
		\centering
		\includegraphics[width=1.0\linewidth]{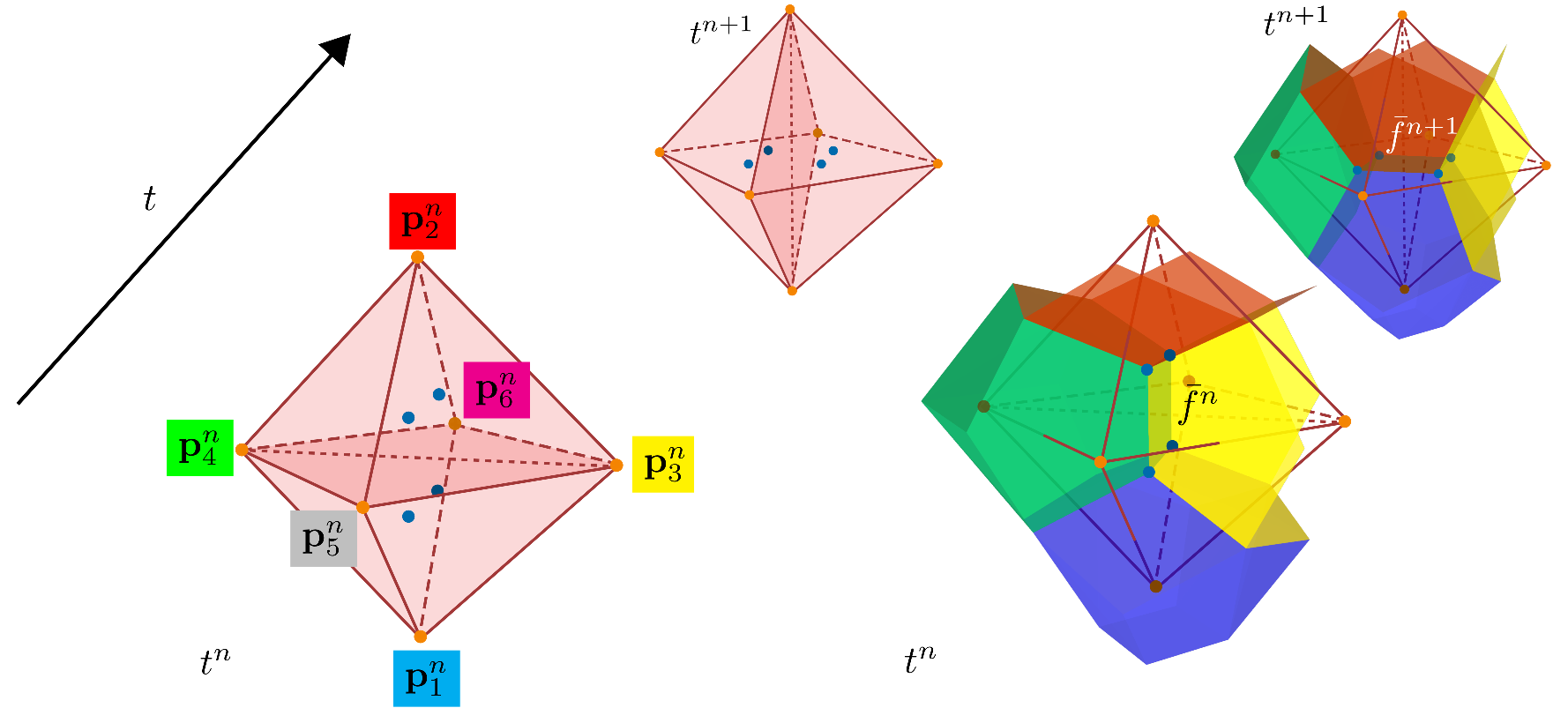}		
		\caption{	
			In this figure, we depict the portion of tetrahedralization (left) where a \textit{4-4 flip} occurs between time $t^n$ and $t^{n+1}$ and the corresponding effect on the polyhedral tessellation (right). In particular, we show 6 generators at time $t^n$ (bottom-left) and at time $t^{n+1}$ (top-left), and 4 of the 6 corresponding polyhedral elements at time $t^n$ (bottom-right) and time $t^{n+1}$ (top-right). 
			At time $t^n$, the six generators are connected via 4 tetrahedra sharing the central edge $\left\{\mathbf{p}_3^n, \mathbf{p}_4^n\right\}$, while at time $t^{n+1}$, they are connected by 4 tetrahedra in a different configuration (the edge connects different generators). 
			This flip also modifies the topology of the corresponding polyhedra, specifically: the green and the yellow polyhedron, at time $t^n$, are neighbors sharing a quadrangular face $\bar{f}^n$; conversely, at time $t^{n+1}$ we find a new face $\bar{f}^{n+1}$ connecting now the red and blue polyhedra.
		}
		\label{fig.4-4flip-origin}
	\end{figure}	
	Moreover, by combining a sequence of a 2-3 flip followed by a 3-2 flip, we obtain another operation that can be considered elementary, which is called \textit{4-4 flip}. As shown in Figure~\ref{fig.4-4flip-origin}, this flip involves 6 generators which are connected via 4 tetrahedra whose connectivity is, however, different before and after the flip: in particular, an edge that centrally connects 2 generators, in the figure the edge $\{ \mathbf{p}_3^n, \mathbf{p}_4^n\}$, will instead connect 2 different generators at the subsequent time. 
	
	Notably, we are able to construct the space-time volumes that succeed in connecting $\mathcal{P}^n$ and $\mathcal{P}^{n+1}$ in the presence of any of these single flips, as we will describe in Section~\ref{s-sec-topology-flips}, and therefore we use them as building blocks of a more complex mesh optimization algorithm, described here below and operating over multiple timesteps.

	\paragraph{Quality criteria for the tetrahedralization}
	\label{s-sec.optimization-criteria}
	We now compute a quality criterion $\alpha_i^n$ for any $T_i^n \in \mathcal{T}^n$ used to connect, instead of the old coordinates $\mathbf{p}^n$, the new Lagrangian and optimized coordinates $\mathbf{p}^{n+1}$; we call $\mathbbm{c}_i^n$ the center of its circumsphere with respect to the new coordinates and $\mathbbm{r}_i^n$ its radius.
	
	Our quality measure penalizes first the tetrahedra that violate the Delaunay condition. 
	In addition, there exist also flat tetrahedra (slivers), in the sense that their maximum dihedral angle $\beta^{\text{dih}}_{i,\text{max}} \rightarrow 180^{\circ}$, which, despite satisfying the Delaunay condition, should be avoided. 
	Therefore our quality parameter $\alpha_i^n$, different from the smoothing parameter~\eqref{eqn.smooth_param},  is computed as
	\begin{equation}
		\label{eq_mesh_quality}
		\alpha_i^n = \min \left(\  1, \ 
		\min \limits_{\,\substack{k \, \in \, \bigcup \limits_{\ell} \mathcal{B}_\ell^n}} \frac{\left\| \, {\mathbf{p}}_k^{n+1} - \mathbbm{c}_i^n \, \right\|}{\mathbbm{r}_i^n}, \ 
		\frac{1 + \cos \beta^{\text{dih}}_{i,\text{max}} }{\,1 + \cos \beta^{\min}\,} \ \right),	
	\end{equation}
	where the index $\ell$ of $\mathcal{B}_\ell^n$ runs over all the vertices of $T_i^n$ and $\beta^{\min} = \cos^{-1} (-0.7) \sim 135^{\circ}$.
	
	According to~\eqref{eq_mesh_quality}, we establish a priority quality queue from the worst tetrahedra, corresponding to the lowest values of $\alpha_i^n$, to the best ones, corresponding to the maximum value of $\alpha_i^n$.

	\paragraph{Priority-based flip sequence search}
	
	Once the priority queue of the tetrahedra that need to be optimized is determined, we begin our search for the elementary transformations that would be necessary to indeed optimize them; we start from the element with the worst quality. 
	
	Our goal is now to find a sequence of elementary operations that maximizes the quality criterion~\eqref{eq_mesh_quality} of this tetrahedron and any others that will be created as a result of the subsequent flips.
	However, we should note that the sequence of the quality values corresponding to the operations leading to a final improvement of the mesh  is not guaranteed to be monotonic, because it is sometimes necessary to temporarily introduce worse quality tetrahedra. 
	This fact prevents us from using a simple algorithm that considers only elementary flips that immediately improve the mesh quality. Thus, we need to predict the transformations that will follow in future steps and evaluate the resulting mesh quality of the possible sequence. Even though we cannot perform the sequence directly, having decided to limit ourselves to one elementary operation per generator and timestep, finding the current optimal one allows us to select the appropriate flip for the current timestep.
	
\begin{figure}[t]\centering  
		\includegraphics[width=1.0\linewidth]{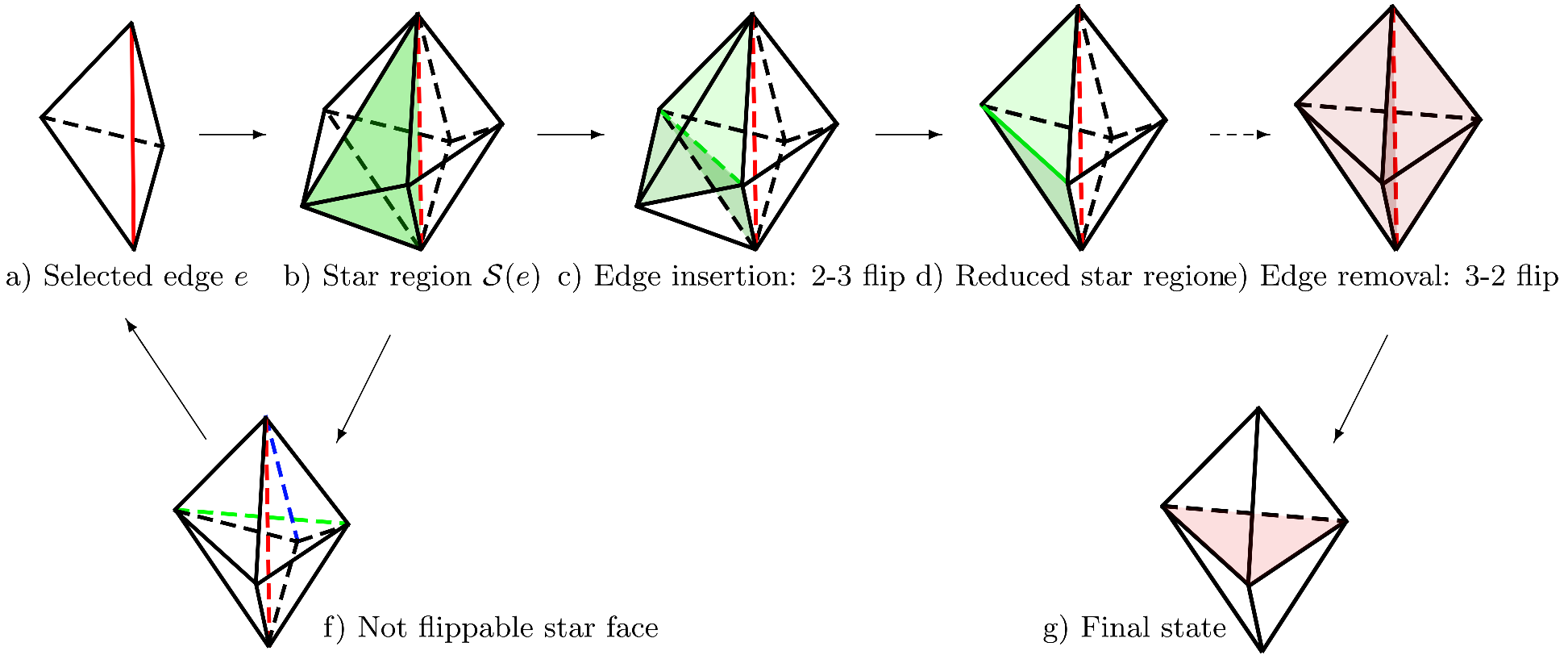}		
		\caption{
			In this image we present an example of the sequence of elementary operations needed to perform a \textit{removal} of a selected edge $e$, shown in red in a), belonging to a bad quality tetrahedron with the scope of improving the quality of the region, as later achieved in f).
			We first construct its star region $\mathcal{S}(e)$ shown in b) defined as the set of all the tetrahedra of the mesh that contain $e$. 
			To remove an edge via a 3-2 flip, the star region must consist of only 3 tetrahedra; if this is not the case, we need to reduce its size by progressively inserting new edges (shown in green in c and d) via 2-3 flips. 
			After an edge is inserted, the outer tetrahedra created by the flip are no more part of the star region, as they no longer contain $e$.           
			However, even if they were not present in the original tetrahedralization $\mathcal{T}^n$, they must be considered when calculating the mesh quality for determining the optimal sequence of flips. 
			By iterating this process the size of the star region is reduced to 3 tetrahedra and $e$ can now be removed by a 3-2 flip, see e).
			Note that, sometimes the neighboring tetrahedra are \textit{not flippable}, i.e., as shown in f) the edge to be inserted (green) would lie outside the star region so the flip cannot be performed right away. In this case we select the closest edge (shown in blue) that does not share any points with the green edge and we proceed with a recursive call of the edge removal algorithm for the blue edge instead.
		}
		\label{fig_edge_remove}
	\end{figure}

	To find the most appropriate elementary operation, we follow the logic of the edge removal search as presented in~\cite{DASSI20182}, and for each tetrahedron, we investigate the consequences of removing each of its edges so to find the optimal sequence. 
	An edge removal sequence for an edge $e$, as also shown in Figure~\ref{fig_edge_remove}, consists of $|\mathcal{S}(e)| - 2$ flips, with possible recursive calls in the case of a non-convex star region. We explore all possible combinations of successive 2-3 flips to reduce the number of tetrahedra in the star region $\mathcal{S}(e)$ to 3, to then perform a final 3-2 flip that removes the edge.
	
	The search conditions are narrowed by a specific requirement to our incremental algorithm that is not present in common global mesh optimization frameworks. Indeed, it is necessary to guarantee that no elementary flip is performed on a tetrahedron that has already been involved in a previous transformation in the current timestep. 
	%
	If the algorithm encounters a configuration that does not satisfy this condition, that sequence is not explored further.
	
	Here, we summarize the mesh optimization algorithm performed during a timestep of our numerical method. 
	\begin{enumerate}[itemsep=2pt, topsep=4pt]
		\item Continue edge removal from previous steps: if there are any edges that have been already marked for removal, continue with the edge removal algorithm for these edges first, performing only the next flip in their respective sequences (see step 3 below).
		\item Create a priority queue of bad tetrahedra, ordered according to the computed quality parameter~\eqref{eq_mesh_quality}.
		\item Iterate through the tetrahedra in the queue, starting from the worst one, and search for the edge removal sequence that results in the best local mesh configuration (see also Figure~\ref{fig_edge_remove}): 
		\begin{enumerate}[label=\alph*), leftmargin=2em, itemsep=2pt, topsep=1pt]
			\item examine all edges $e$ of the tetrahedron;
			\item determine their star region $\mathcal{S}(e)$;
			\item explore all possible 2-3 flips within the star region;
			\item obtain a configuration with a reduced number of tetrahedra;
			\item repeat steps c) and d) until a configuration with 3 tetrahedra is reached.
			\item If it is \textit{not} possible to perform a 2-3 flip on a face, recursively call the edge removal algorithm for the edge that is closest to the non-convex part of the star region.
			\item Identify the final 3-2 flip to remove the edge.
		\end{enumerate}
		\item Once the optimal sequence is found, perform only the first flip of the sequence and mark the edge to continue with its removal in the following timesteps.
	\end{enumerate}

	\section{Space-time connectivity via control volumes filling the 4D space-time slices}
	\label{sec.space-time} 
	
	{
		\linespread{1.1}\selectfont	
		At this point we have $\mathcal{T}^n$ and $\mathcal{P}^n$ discretizing $\Omega^n$ at $t^n$, 
		and $\mathcal{T}^{n+1}$ and $\mathcal{P}^{n+1}$ discretizing $\Omega^{n+1}$ at $t^{n+1}$.
		Moreover, we precisely know that $\mathcal{P}^{n}$ and $\mathcal{P}^{n+1}$ have the same number of polyhedral elements and 
		that any topological difference between the two is due to at most a single 2-3, 3-2, or 4-4 flip affecting each generator.
		
		As will be explained in the next Section~\ref{sec.method-description}, a fundamental aspect for 
		designing a numerical scheme that explicitly evolves the PDE from $t^n$ to $t^{n+1}$, 
		moving the solution directly from each element $P_i^n$ to the corresponding $P_i^{n+1}$, 
		is the capability of \textit{integrating} the PDE \textit{over space-time control volumes} that connect the spatial elements. 
		While we defer the full methodological description to Section~\ref{sec.method-description}, 
		we focus here on the precise identification, construction, and characterization 
		of all the required control volumes. This includes a control volume 
		corresponding to each $P_i^n$, as well as additional special space-time 
		control volumes referred to as \textit{hole-like} elements.
		\par }
	
	\subsection{Classical 4D space-time control volumes and their 3D space-time lateral surfaces} 
	\label{s-sec.standardcontrolvolumes}
	
	{
		\linespread{1.05}\selectfont	
		We first consider the case of elements $P_i^n$ and $P_i^{n+1}$ with the \textit{same topology}, 
		i.e., the same shape and the same connectivity at the beginning and at the end 
		of the timestep. In this case, we can naturally establish a correspondence 
		between their faces: in particular, the face $f_{ij}^n$ of $P_i^n$ separating 
		it from its neighbor $P_j^n$ corresponds to the face $f_{ij}^{n+1}$ of 
		$P_i^{n+1}$ separating it from $P_j^{n+1}$, for every $j \in \mathcal{V}_i^n$.
		This situation is the classical one and it is verified when: 
		\begin{enumerate}
			\item[i)] $P_i^n$ and $P_i^{n+1}$ have the same set of neighbors at $t^n$ and $t^{n+1}$, 
			i.e., $\mathcal{V}_i^n = \mathcal{V}_i^{n+1}$, 
			\item[ii)] each corresponding face has the same number of facets at $t^n$ and $t^{n+1}$, namely the star region around $e_{ij}^{n}$ and $e_{ij}^{n+1}$ has the same cardinality, i.e. $|\mathcal{S}(e_{ij}^n)| = |\mathcal{S}(e_{ij}^{n+1})|$ for every $j \in \mathcal{V}_i^n$.
		\end{enumerate}	
		Of course, instead, the coordinates of the nodes in $\mathcal{D}_i^n$ and $\mathcal{D}_i^{n+1}$ may be all different.
		\par }
	
	\begin{figure}[tb]
		\centering
		\includegraphics[width=1.0\linewidth]{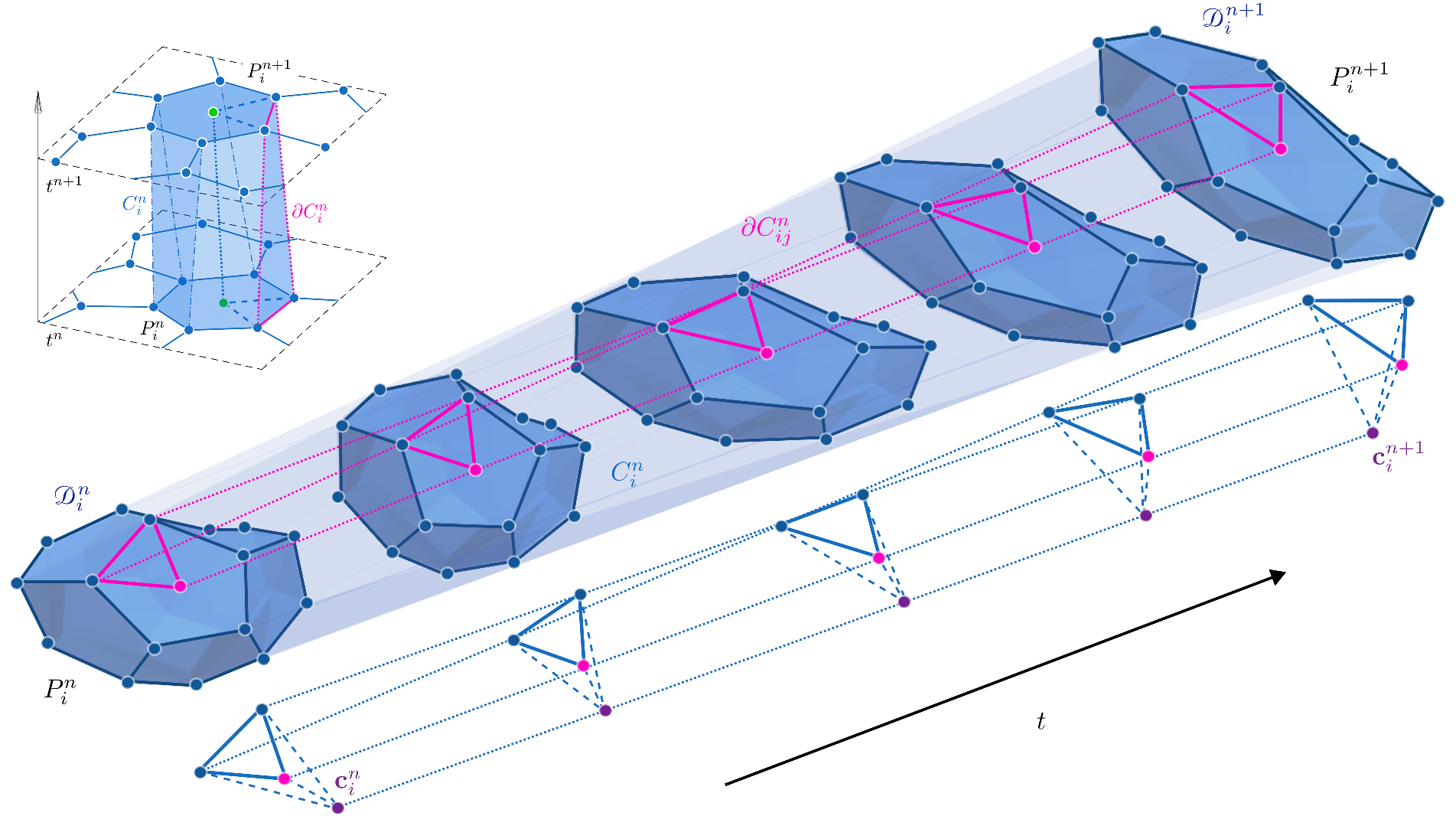}		
		\caption{In this figure, we depict in blue the concept of the 4D space-time control volume $C_i^n$ 
			connecting in space and time the 3D elements $P_i^n$ (bottom-left) and $P_i^{n+1}$ (top-right). 
			To better convey the idea we also report, in the top-left part of the image,
			the sketch that illustrates a 3D 
			space-time control volume $C_i^n$ connecting 2D polygons, as presented in 
			previous works by the authors~\cite{gaburro2020high,gaburro2025high}. 
			The 4D control volumes, which are the subject of this work, cannot be fully 
			visualized on a two-dimensional sheet of paper; however, to provide an 
			intuition, we adopt the strategy of displaying five snapshots of it at different 
			time instants $t^{n+\tau}$ in between $t^n$ and $t^{n+1}$. 
			We also depict in pink its so-called space-time lateral surface $\partial C_{ij}^n$ which is a sort of prism connecting in space and time 
			a triangular facet of $P_i^n$ with the corresponding facet in $P_i^{n+1}$. 
			Finally, we highlight one sub-tetrahedron in which $P_i^{n/n+1}$ could be subdivided and we show its space-time evolution by providing five of its snapshots. 					
		}		
		\label{fig.controlvolumes}
	\end{figure}	

	\subsubsection{4D space-time control volumes}
	{
		\linespread{1.05}\selectfont	
		In this case, we can build the control volume $C_i^n$ by simply connecting $P_i^n$ 
		and $P_i^{n+1}$ in space and time. This straightforward construction results in a 
		4D volume, see also Figure~\ref{fig.controlvolumes}, 
		that at $t=t^n$ coincides with $P_i^n$, at $t=t^{n+1}$ coincides with $P_i^{n+1}$, 
		and in between is obtained as their convex combination.
		For example, at any fixed time $t= t^{n+\tau} = t^n + \tau\Delta t$ with $\tau \in [0,1]$, 
		the 3D time-slice of $C_i^n$ is given by the polyhedral element bounded by the vertices 
		\begin{equation*} 
			\textbf{d}_{ij}^{n+\tau} = (1-\tau) \, \textbf{d}_{ij}^{n} + \tau \, \textbf{d}_{ij}^{n+1}, \quad \text{for } \textbf{d}_{ij}^n \in \mathcal{D}_i^n, \textbf{d}_{ij}^{n+1} \in \mathcal{D}_i^{n+1}, j \in \left\{1,\dots, |\mathcal{D}_i^n|\right\} .
		\end{equation*}
		\par }

	\subsubsection{3D space-time lateral surfaces}\label{s-sec.3dlateralsurface}
	The total lateral surface of $C_i^n$ is denoted with $\partial C_i^n$ and is made as follow
	\be 
	\partial C_i^n =  \bigcup_{j \, \in \, \mathcal{W}_i^{n}} \partial C_{ij}^n \cup P_i^n \cup P_i^{n+1},
	\label{eqn.dC}
	\ee		
	where $\mathcal{W}_i^{n}$ indicates the set of indices of the space-time control volumes neighboring $C_i^n$ which here coincides with the control volumes built in correspondence of each index of $\mathcal{V}_i^n$, thus $\mathcal{W}_i^n = \mathcal{V}_i^n = \mathcal{V}_i^{n+1}$.	
	
\begin{figure}[tb]
		\centering
		\begin{minipage}[c]{0.67\textwidth}
			\vspace{0pt} 
			\centering
			\includegraphics[width=1.0\linewidth]{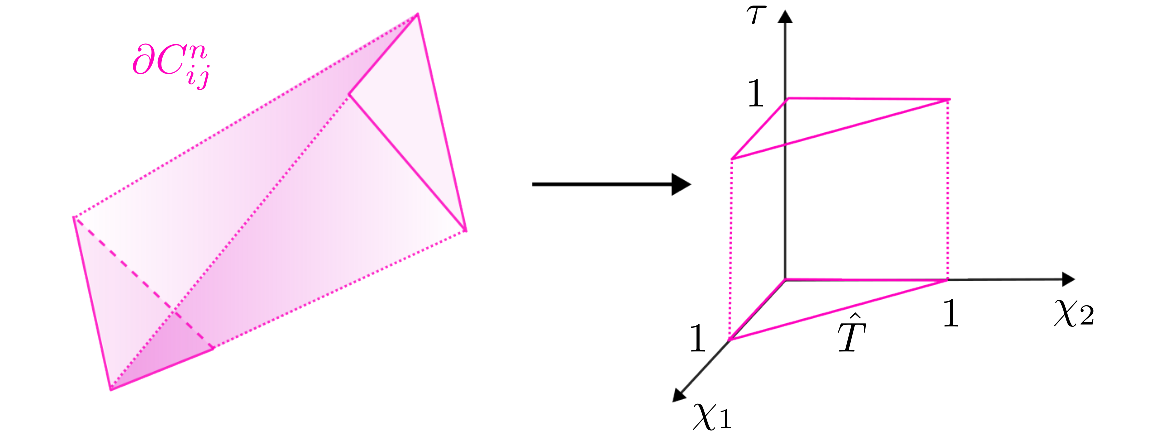}
		\end{minipage}
		\hfill 
		\begin{minipage}[c]{0.3\textwidth}
			\vspace{0pt} 
			\caption{Mapping of a space-time lateral surface~$\partial C_{ij}^n$ (left) to a reference element composed of a triangle $\hat{T}$ extruded in time over $[0,1]$ (right).						
			}	
			\label{fig.lateralsurf_standard}
		\end{minipage}
	\end{figure}

	{
		\linespread{1.1}\selectfont		
		The space-time lateral surfaces $\partial C_{ij}^n$ of these 4D control volumes are 3D curvilinear objects, see both Figures~\ref{fig.controlvolumes} and~\ref{fig.lateralsurf_standard}, characterized by a non-trivial curvature. They can be conceptualized as a union of curvilinear prisms whose "bases" are the triangular facets of each face of $P_i^{n/n+1}$. Specifically, the "bottom base" of the prism coincides with a triangular facet of $P_i^n$ at time $t^n$, while the "top base" coincides with the corresponding facet at time $t^{n+1}$ and the "height" is given by $\Delta t$.
		\par }
	
	The lateral faces of these prisms are generally not planar; they possess a curvature that requires a more sophisticated description. To accurately represent these surfaces, we employ a set of bi-linear basis functions \cite{Lagrange3D}, known as the $\beta$-basis functions. 
	In this 3D space-time context, any point $\hat \x$ on $\partial C_{ij}^n$ is defined through a mapping from a reference element (a reference triangle $\widehat{T}$ extruded in time over $[0,1]$). 
	Using the space-time coordinates of the 6~vertices, namely $\hat{\x}_1, \hat{\x}_2$ and $\hat{\x}_3$ vertices of a triangular facet at time $t^n$, each connected (respectively) to $\hat{\x}_4, \hat{\x}_5$ and $\hat{\x}_6$ vertices of the same facet at time $t^{n+1}$, the mapping is expressed as:
	\begin{equation}
		\label{eq.beta}
		\hat \x(\zeta_1, \zeta_2, \tau) = \sum_{j=1}^{6} \beta_j(\zeta_1, \zeta_2, \tau) \hat\x_j,
	\end{equation}
	where $(\zeta_1, \zeta_2)$ are the spatial coordinates on the reference triangle and $\tau \in [0,1]$ is the normalized time. The $\beta$-basis functions takes the form
	\begin{equation*}
		\begin{cases} 
			\beta_1(\zeta_1, \zeta_2, \tau) = (1-\tau) \zeta_1 \\ 
			\beta_2(\zeta_1, \zeta_2, \tau) = (1-\tau) \zeta_2 \\ 
			\beta_3(\zeta_1, \zeta_2, \tau) = (1-\tau) (1-\zeta_1-\zeta_2)
		\end{cases} \quad \text{and} \quad 
		\begin{cases} 
			\beta_4(\zeta_1, \zeta_2, \tau) = \tau \zeta_1 \\ 
			\beta_5(\zeta_1, \zeta_2, \tau) = \tau \zeta_2 \\ 
			\beta_6(\zeta_1, \zeta_2, \tau) = \tau (1-\zeta_1-\zeta_2)
		\end{cases}.
	\end{equation*}
	The resulting surface is bi-linear: it remains linear in space for any fixed instant $\tau$ and linear in time for any fixed spatial point $(\zeta_1, \zeta_2)$. This mathematical characterization is essential for the consistent calculation of the space-time normal vectors and to ensure the Geometric Conservation Law (GCL) is satisfied during the mesh movement.

	\paragraph{Geometrical information on the space-time lateral surfaces}
	The use of the $\beta$-basis functions is particularly advantageous for the computation of the space-time normal vectors ${\mathbf{n}}_{ij}^n$ and the measures of the lateral surfaces $\partial C_{ij}^n$. 
	Indeed, at each point $\hat\x(\zeta_1, \zeta_2, \tau) = (x_1(\zeta_1, \zeta_2, \tau), x_2(\zeta_1, \zeta_2, \tau), x_3(\zeta_1, \zeta_2, \tau), t(\zeta_1, \zeta_2, \tau))$ on $\partial C_{ij}^n$, we can compute the normal space-time vector $\n(\zeta_1, \zeta_2, \tau)$ through the generalized cross product of its partial derivatives: i.e.,
	\begin{equation}
		(\partial_{\zeta_1} \hat\x \times \partial_{\zeta_2} \hat\x \times \partial_{\tau} \hat\x)(\zeta_1, \zeta_2, \tau) = \text{det} 
		\begin{pmatrix} 
			\mathbf{i} & \mathbf{j} & \mathbf{k} & \mathbf{l}  \\
			\partial_{\zeta_1} x_1(\zeta_1, \zeta_2, \tau) & \partial_{\zeta_1} x_2(\zeta_1, \zeta_2, \tau) & \partial_{\zeta_1} x_3(\zeta_1, \zeta_2, \tau) & \partial_{\zeta_1} t(\zeta_1, \zeta_2, \tau)  \\
			\partial_{\zeta_2} x_1(\zeta_1, \zeta_2, \tau) & \partial_{\zeta_2} x_2(\zeta_1, \zeta_2, \tau) & \partial_{\zeta_2} x_3(\zeta_1, \zeta_2, \tau) & \partial_{\zeta_2} t(\zeta_1, \zeta_2, \tau)  \\
			\partial_{\tau} x_1(\zeta_1, \zeta_2, \tau)    & \partial_{\tau} x_2(\zeta_1, \zeta_2, \tau)    & \partial_{\tau} x_3(\zeta_1, \zeta_2, \tau) & \partial_{\tau}    t(\zeta_1, \zeta_2, \tau)
		\end{pmatrix},
		\label{eq.surfaceintegral-term}
	\end{equation}
	where $\mathbf{i, j, k, l}$ are the unit vectors of the 4D space-time coordinate system. The space-time unit normal vector ${\n}(\zeta_1, \zeta_2, \tau)$ is then obtained by normalizing the previous result. In particular, the generalized cross product $\partial_{\zeta_1} \hat\x \times \partial_{\zeta_2} \hat\x \times \partial_{\tau} \hat\x$ represents the local Jacobian of the transformation mapping defining $\partial C_{ij}^n$, hence
	\begin{equation}
		|\partial C_{ij}^n| = \int_0^1 \int_{\widehat{T}} |(\partial_{\zeta_1} \hat\x \times \partial_{\zeta_2} \hat\x \times \partial_{\tau} \hat\x)(\zeta_1, \zeta_2, \tau)| \, d\zeta_1 d\zeta_2 d\tau.
		\label{eq.surfaceintegral}
	\end{equation}
	The same Jacobian determinant appears when dealing with integration of a given function defined on $\partial C_{ij}^n$.

	\subsection{Topology changes, \textit{hole-like} elements and their neighborhood}
	\label{s-sec-topology-flips}
	
	Let us now consider the situations in which, due to the mesh optimization driven by the motion of the generators, 
	some \textit{flips} occur. 
	These are rearrangements of the connectivity of the tetrahedralizations such that $\mathcal{T}^n$ differs from $\mathcal{T}^{n+1}$, leading to topology changes that affect $\mathcal{P}^n$ and $\mathcal{P}^{n+1}$, making them different. 
	Hence, in this section, we deal with the neighborhood of an element $P_i^{n/n+1}$ that do not satisfies the conditions listed in Section~\ref{s-sec.standardcontrolvolumes}.
	
	The issue in these situations is that, even if we find a way to connect $P_i^n$ and $P_i^{n+1}$ with different topology (different shape and connectivity) through volumes $C_i^n$, 
	as indeed done in Section~\ref{s-sec.controlvolumes-neigh-hole}, 
	their union would not be sufficient to cover the entire space-time slice between $t^n$ and $t^{n+1}$, 
	thus leaving some gaps. 
	Our primary goal is, therefore, to define specific space-time control volumes capable of exactly filling these 4D gaps. 
	We refer to these elements as \textit{hole-like} elements, and we denote them with $H_i^n$, where indices $i > N_P$ span the set of all necessary \textit{hole-like} elements that compose the collection $\mathcal{H}^n = \{ H_{i}^n \}_i$.
	
	These \textit{hole-like} elements are characterized by the fact that at $t^n$ and $t^{n+1}$ they do not coincide with any element $P_i^{n/n+1}$, nor with anything having a non-zero 3D volume; 
	for this reason, we refer to them as \textit{spatially degenerate}. 
	As we will see, at time $t^n$ or $t^{n+1}$, they will coincide either with faces or with edges of the tessellations $\mathcal{P}^{n/n+1}$. 
	However, these \textit{hole-like} elements $H_i^n$ possess a \textit{space-time volume} measured in the 4D space that, while small, is \textit{strictly non-zero}. 
	Neglecting them would lead to a lack of conservativity in the numerical method. 
	Furthermore, their lateral surfaces have a significant dimension, comparable to the size of the faces of standard elements; ignoring them would therefore mean omitting fundamental flux exchanges between the elements involved in the 
	topology change.
	
	We remark that, according to the mesh optimization strategy described in Section~\ref{s-sec.optimization-flips}, only one elementary flip of type 2-3, 3-2, or 4-4 occurs in each neighborhood at a time. 
	Consequently, the entire space-time slice can be completely covered by considering just three types of \textit{hole-like} elements (as detailed in Sections~\ref{s-sec.elementaryflip3-2} and~\ref{s-sec.elementaryflip4-4}) which, while potentially multiple at each timestep, cannot be mutually adjacent. This non-adjacency is indeed enforced by our optimization, which is limited to one elementary flip per generator. The neighbors of each \textit{hole-like} element will be the space-time control volumes $C_i^n$, which connect polyhedra in a manner essentially identical to that described in Section~\ref{s-sec.standardcontrolvolumes}, requiring only slight modifications to the 3D lateral faces shared with the \textit{hole-like} elements~(see~Section~\ref{s-sec.controlvolumes-neigh-hole}).

	\subsubsection{Hole-like elements corresponding to 3-2 and 2-3 flips}
	\label{s-sec.elementaryflip3-2}	
	
	\begin{figure}[tb]
		\centering		
		\includegraphics[width=1.0\linewidth]{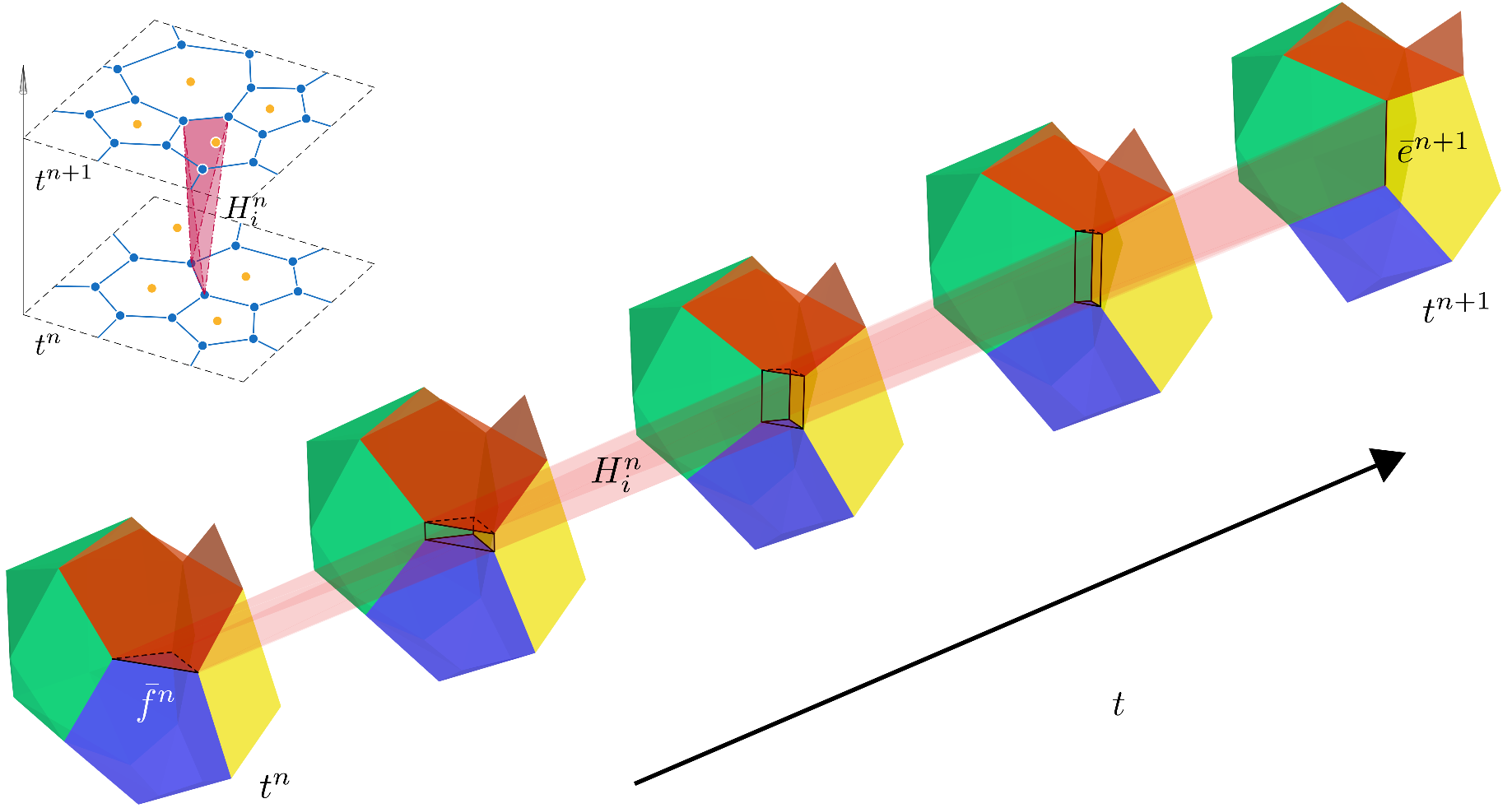}		
		\caption{
			In this figure, we depict the concept of the 4D \textit{hole-like} element $H_i^n$ through the space-time slab shown in pink and 5 of its snapshots between $t^n$ and $t^{n+1}$, together with the geometry of 4 of its 5 neighbor elements. 
			To provide a clearer intuition of our novel 4D element, we also report the equivalent 3D \textit{hole-like} element in the top-left of the figure, introduced in previous works of the authors~\cite{gaburro2020high,gaburro2025high}, which served to fill the space-time gap between two 2D tessellations in correspondence with a topology change.
			Here, our novel \textit{hole-like} element is a 4D object that, at each fixed instant in time, has the shape of a triangular prism, traced with solid black lines; this element is thus delimited at each instant by 6 vertices with space-time coordinates.
			In particular, as $t \rightarrow t^n$, it degenerates into a prism with zero height, corresponding to the planar surface $\bar{f}^n$; therefore, at each corner of the triangle, we have two overlapping vertices of the {flat prism}. As $t \rightarrow t^{n+1}$, it instead degenerates into a prism with zero-area triangular faces, corresponding to the segment $\bar{e}^{n+1}$; thus, at each extremity of the edge, we have three overlapping vertices of the {slender prism}.		
		}
		\label{fig.3-2hole}
	\end{figure}
	\begin{figure}[tb]\centering
		\includegraphics[width=1.0\linewidth]{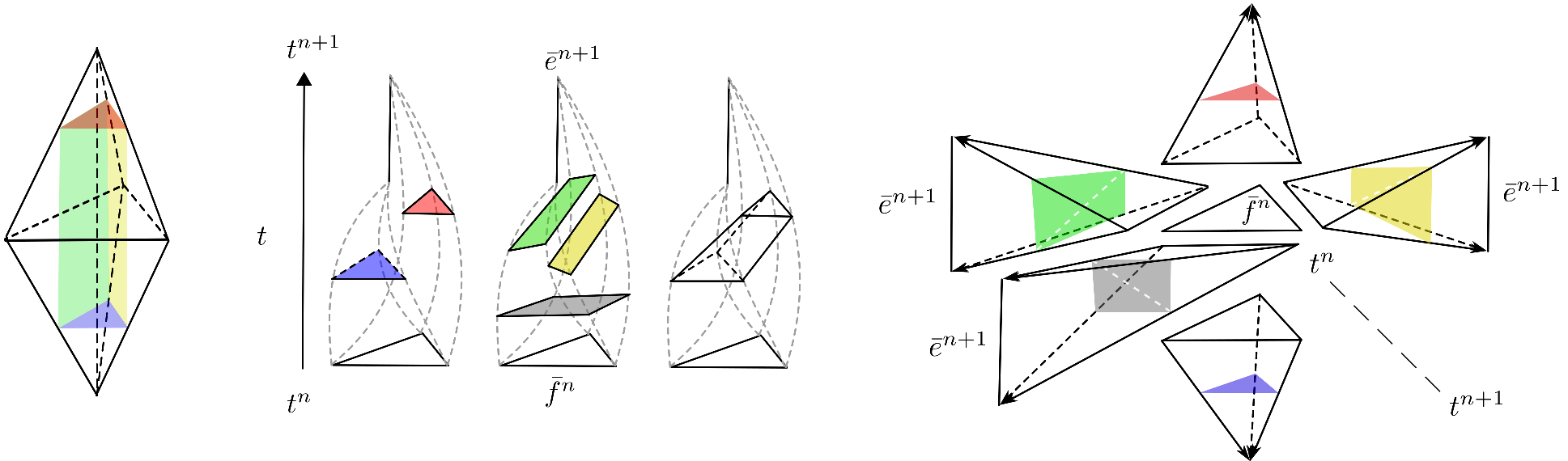}		
		\caption{ 
			In this figure, we provide three additional different visualization strategies for the \textit{hole-like} element $H_i^n$ and its space-time lateral surfaces $\partial H_{ij}^n$.
			In particular, on the left, we propose a spatial projection of the \textit{hole-like} element in the form of a 4-simplex (or \textit{pentatope}). 
			It is showing the spatial 3D volume swept by the \textit{hole-like} element along its time trajectory.
			In the middle, where time evolves along the vertical axis, we show the trajectories spanned by the surfaces $\partial H_{ij}^n$ (with the same color-code of the previous images) along with the trajectory of the \textit{hole-like} volume.
			Finally, on the right, we provide an exploded representation where the triangular face $\bar{f}^n$ is visible at the center, while at the left and right extremities we find the edge $\bar{e}^{n+1}$, and above and below, the endpoints of said edge; here, the time direction evolves from the center outwards. In this third representation, the 3D space-time lateral surfaces $\partial H_{ij}^n$ are clearly displayed, providing the necessary details for their construction. In fact, the red and blue surfaces connect $\bar{f}^n$ with the endpoints of $\bar{e}^{n+1}$ and consist of a single 3D curvilinear element; conversely, the yellow, green, and gray ones each connect a side of $\bar{f}^n$ with $\bar{e}^{n+1}$ and are each composed of two parts.
		}
		\label{fig.3-2holesurfaces}
	\end{figure}
	
	Type 2-3 and type 3-2 flips are reciprocal one another: in this section, we focus on the 3-2 flip and the opposite case can be straightforwardly derived by inverting $t^n$ and $t^{n+1}$.	
	We also refer to Figure~\ref{fig.3-2flip-origin} for a visual interpretation of this process. 
	
	The situation involves five generators $\mathbf{p}_j^{n/n+1}$, $j=1 \dots 5$,  and, accordingly, the corresponding polyhedra $P_j^{n/n+1}$. 
	In the figure (right), we show portions of 4 of these 5 polyhedra, each marked with a different color: 
	blue ($P_1^n$), red ($P_2^n$), yellow ($P_3^n$), and green~($P_4^n$). 
	We omit the fifth one ($P_5^n$), which we refer to as gray, since it is located in front of the others and would obstruct the view. 
	On the left of the figure, we illustrate the portion of the tetrahedralization where the flip occurs.
	
	Our primary interest lies in understanding the effect this (standard) elementary transformation of the tetrahedralization has on $\mathcal{P}^{n/n+1}$: 
	indeed, the five elements corresponding to these generators are all affected by a change in their topology. 
	Specifically:
	\begin{itemize}
		\item {$P_1^n$ (blue) and $P_2^n$ (red)}, at time $t^n$, are neighbors sharing a triangular face $\bar{f}^n$; conversely, at time $t^{n+1}$, this face disappears and they no longer touch. A connection between them can only be found in the fact that each of them has as vertex one of the extremities of the new edge $\bar{e}^{n+1}$.
		
		\item {$P_3^{n+1}, P_4^{n+1}$, and $P_5^{n+1}$} share the edge $\bar{e}^{n+1}$ at time $t^{n+1}$, while at time $t^n$ they are not in direct contact. Their connection at~$t^n$ is represented by the fact that they all have as edge one boundary of the face $\bar{f}^n$.
	\end{itemize}	
	In this scenario, if we were to construct only the standard control volumes $C_i^n$ connecting points for which a correspondence had eventually been found, we would leave a gap in the space-time domain. We therefore aim to characterize this gap and fill it with a \textit{hole-like} element.
	
	Our novel \textit{hole-like} element $H_i^n$ is shown in Figure~\ref{fig.3-2hole}, adopting again the strategy of illustrating its configuration at several intermediate instants \textit{in between} $t^n$ and $t^{n+1}$. It is a 4D object that, at each fixed instant in time $t^{n+\tau}$, has \textit{the shape of a triangular prism} (thus having 6 vertices with space-time coordinates).
	In particular, as $t \rightarrow t^n$, the element degenerates into a prism with zero height, corresponding to the planar surface $\bar{f}^n$; therefore, at each corner of the triangle, we have two overlapping vertices of the {flat prism}. 
	As $t \rightarrow t^{n+1}$, it instead degenerates into a prism with zero-area triangular faces, corresponding to the segment $\bar{e}^{n+1}$; thus, at each extremity of the edge, we have three overlapping vertices of the {slender prism}.
	In between, the prism is defined by the convex combination, via the parameter~$\tau$, of these 6 moving points. Specifically, it connects the 6 points at $t^n$ (which overlap in pairs) with the 6 points at $t^{n+1}$ (which overlap in triplets).
	
	This kind of \textit{hole-like} element always has 5 neighbors, so $|\mathcal{W}_i^n|=5$, specifically all the 5 control volumes $C_j^n$ originating from the polyhedra $P_j^n$ for $j=1 \dots 5$, and thus 5 space-time lateral faces $\partial H_{ij}^n$. 
	In particular, each $\partial H_{ij}^n$ represents the space-time interface between $H_i^n$ and $C_j^n$.
	Interestingly, these lateral surfaces do not require a dedicated formulation; they are ultimately constructed and parameterized exactly like standard ones, with the only difference being their potential degeneracy at time $t^n$ and/or $t^{n+1}$.
	Indeed, in this case as well, they consist of a collection of 3D triangular prisms (with curved surfaces described by the $\beta$-basis functions as in \eqref{eq.beta}). 
	For a visual interpretation, we refer to Figure~\ref{fig.3-2holesurfaces}. 
	The "height" of these curvilinear prisms is always given by $\Delta t$. 
	Regarding the "bottom" and "top" "bases" of the prisms, we distinguish between two cases:	
	\begin{itemize} 
		\item $\partial H_{i1}^n$ and $\partial H_{i2}^n$, at time $t^n$, correspond to the face $\bar{f}^n$ (which serves as base of the prism) and degenerate~for $t \rightarrow t^{n+1}$ into the two extremities of $\bar{e}^{n+1}$, respectively. Although these extremities are points, they can still be~viewed as degenerate triangles; therefore, each of these surfaces can be fully described by a single curvilinear~prism.
		
		\item $\partial H_{i3}^n, \partial H_{i4}^n$, and $\partial H_{i5}^n$ are each composed of two curvilinear prisms with both bottom and top "triangular" faces completely degenerated into segments. Specifically, the top faces all coincide with the segment $\bar{e}^{n+1}$ at time $t^{n+1}$, while the bottom faces correspond to one of the edges of the face $\bar{f}^n$ at time $t^n$. However, these {slender elements} at the extrema of the time slab expand into two {triangular prisms} in space-time.
	\end{itemize}

	\subsubsection{Hole-like elements corresponding to 4-4 flips}
	\label{s-sec.elementaryflip4-4}
	
	\begin{figure}[!tb] 
		\centering	
		\begin{minipage}{\linewidth}
			\centering
			\includegraphics[width=1.0\linewidth]{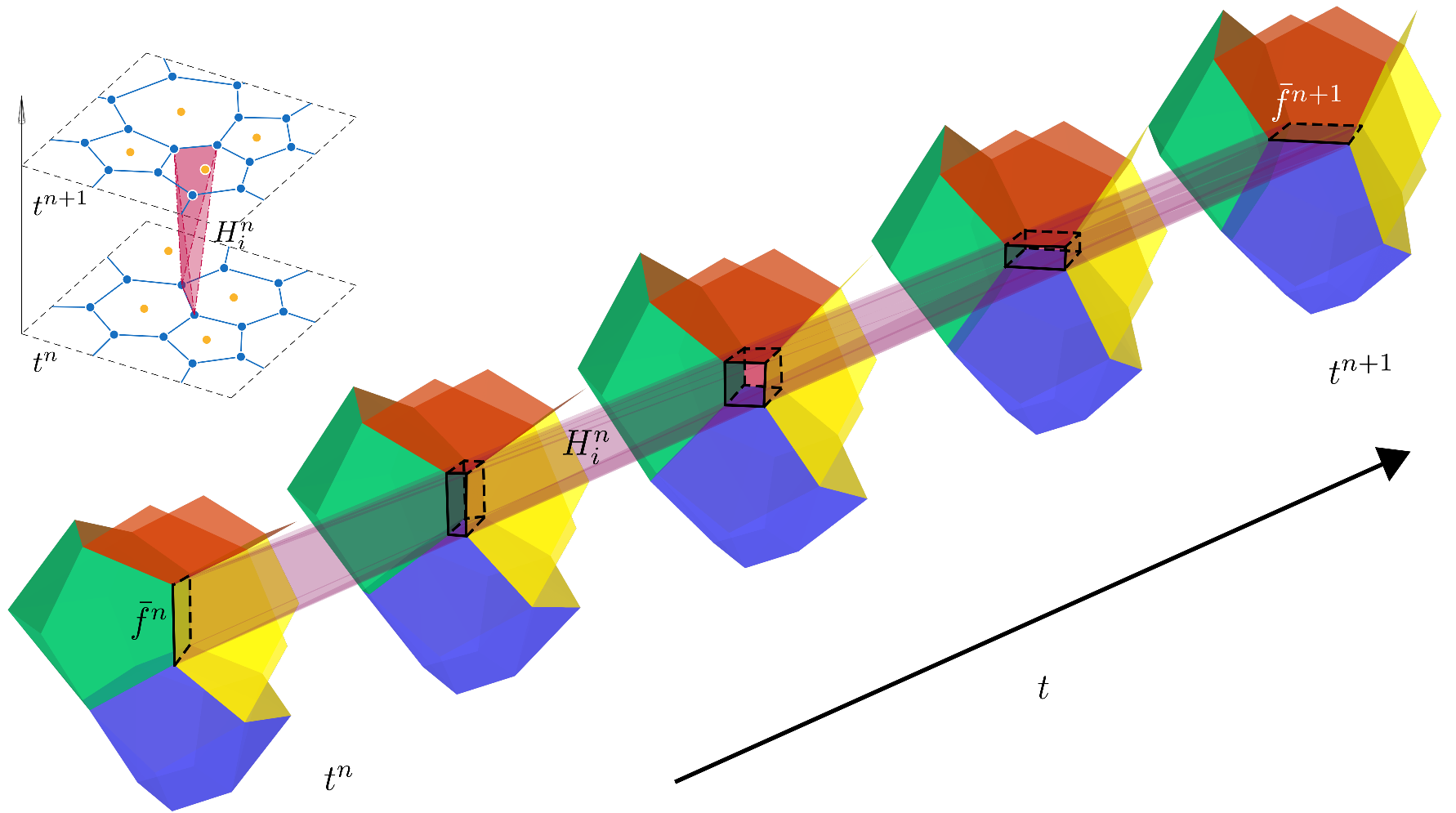}		
			\caption{%
				In this figure, the 4D \textit{hole-like} element $H_i^n$ corresponding to a 4-4 flip is depicted as a pink space-time slab, along with five time snapshots of the element between $t^n$ and $t^{n+1}$. At each fixed instant, the \textit{hole-like} element has the shape of a parallelepiped, displayed in solid black lines; this element is thus delimited at each instant by 8 vertices with space-time coordinates. Additionally 4 of its 6 neighboring elements are displayed (excluding the front and rear neighbors). 
				In particular, as $t \rightarrow t^n$ or $t \rightarrow t^{n+1}$, the element degenerates into a flat shape, corresponding respectively to the quadrangular surfaces $\bar{f}^n, \bar{f}^{n+1}$.}
			\label{fig.4-4hole}
		\end{minipage}	
		\vspace{2em} 
		\begin{minipage}{\linewidth}
			\centering	
			\includegraphics[width=0.65\linewidth]{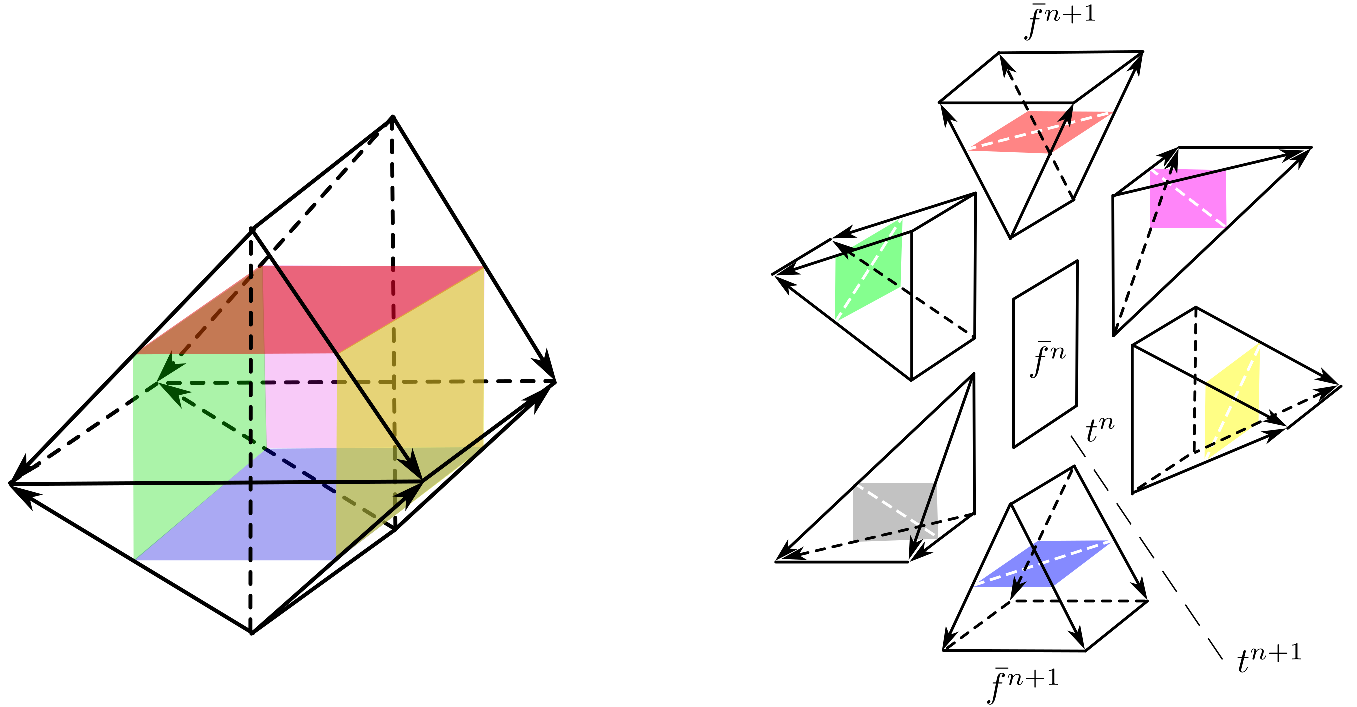}		
			\caption{
				In this figure, we provide additional visualization strategies for the \textit{hole-like} element corresponding to a 4-4 flip and its space-time lateral surfaces.
				On the left, a spatial projection of the 4-dimensional \textit{hole-like} element is shown, representing the spatial 3D volume swept by the element during its evolution in time.
				On the right we provide an exploded representation where the time direction evolves from the center outwards. 
				The initial face $\bar{f}^n$ is visible at the center, while at the top and bottom we see the final face $\bar{f}^{n+1}$. On left, right, front and back,  we can see the 4 edges of face $\bar{f}^{n+1}$. 
				Together here we display the 3D space-time lateral surfaces $\partial H_{ij}^n$ as the 3D volumes connecting $\bar{f}^n$ (or its edges) with $\bar{f}^{n+1}$ (or its edges).
			}        
			\label{fig.4-4holesurface}
		\end{minipage}
	\end{figure}

	The 4-4 flips can be interpreted as a succession of a 2-3 flip followed by a 3-2 flip, eventually performed at the subsequent timestep.
	However, these two flips can also be treated in a single step because they result in a easy to handle and unique \textit{hole-like} element.
	
	Here, the process involves 6 generators $\textbf{p}_j^{n/n+1}$, $j=1 \dots 6$, 
	and the corresponding 6 polyhedra $P_j^{n/n+1}$, $j=1 \dots 6$, 
	as shown in Figure~\ref{fig.4-4flip-origin}.
	In the tetrahedralization, the points are always connected by 4 tetrahedra but arranged in different configurations before and after the flip. 
	This transition produces topological changes in all the 6 involved polyhedra.
	
	Following a reasoning similar to the previous case, the required \textit{hole-like} element must connect 
	the planar quadrangular face $\bar{f}^n$ shared by polyhedra $P_3^n$ and $P_4^n$ (yellow and green) at time $t^n$, with the planar quadrangular face $\bar{f}^{n+1}$ shared by $P_1^{n+1}$ and $P_2^{n+1}$ (red and blue) at time $t^{n+1}$. 
	The resulting \textit{hole-like} element is a 4D object (defined by 8 vertices with space-time coordinates) that, 
	at any fixed instant $t \in (t^n, t^{n+1})$, takes the shape of a \textit{rectangular prism}.
	
	For $t \rightarrow t^n$, the \textit{hole-like} element degenerates into the face $\bar{f}^n$, which has zero 3D volume (with two nodes of the prism coinciding at each of its four vertices). Similarly, for $t \rightarrow t^{n+1}$, it degenerates into the face $\bar{f}^{n+1}$ (again, with two nodes coinciding at each vertex). However, in the space-time domain, the \textit{hole-like} element has a non-zero volume bounded by the trajectories connecting the 8 points at time $t^n$ (which overlap in pairs) to the 8 points at time $t^{n+1}$ (which also overlap in different pairs) via a convex combination with parameter $\tau \in [0, \Delta t]$.
	
	This kind of \textit{hole-like} elements always has 6 neighbors so $|\mathcal{W}_i^n|=6$, specifically all the 6 control volumes $C_j^n$ originating from the polyhedra $P_j^n$, $j=1 \dots 6$, and thus 6 space-time lateral faces $\partial H_{ij}^n$. 
	For their visual interpretation, we refer to Figure~\ref{fig.4-4holesurface}. In this case, each $\partial H_{ij}^n$ is composed of two 3D curvilinear  triangular prisms, whose "height" is defined by $\Delta t$. Their bases coincide with one half of the quadrangular faces $\bar{f}^n$ and $\bar{f}^{n+1}$. This construction requires a coherent ordering of the delimiting nodes and careful attention to their association between the $t^n$ and $t^{n+1}$ time levels.

	\subsubsection{Neighbors of \textit{hole-like} elements}
	\label{s-sec.controlvolumes-neigh-hole}
	
	{
		\linespread{1.1}\selectfont
		A final point concerns the construction of control volumes $C_j^n$ corresponding to the polyhedra $P_j^n$ affected by a topology change when a flip occurs.
		In these cases, the associated control volumes $C_j^n$ possess a new face neighboring the \textit{hole-like} element $H_i^n$. 
		While $H_i^n$ does not exist spatially and is not included in the spatial neighbors $\mathcal{V}_j^{n/n+1}$, it exists in the space-time domain and must be taken into account into the list of space-time neighbors of $C_j^n$, which is given by the common neighbors at time $t^{n/n+1}$ plus the \textit{hole-like} element, i.e., $\mathcal{W}_j^n = (\mathcal{V}_j^n \cap \mathcal{V}_j^{n+1}) \cup \{i\}$. This introduces a space-time interface (potentially split into two facets) between $C_j^n$ and the \textit{hole-like} element, consisting of 3D curvilinear prisms with possibly degenerate faces at $t^n$ or $t^{n+1}$.
		
		Furthermore, in correspondence of these degenerate triangular facets we have to build additional sub-tetrahedra of $P_j^{n/n+1}$ 
		which have zero 3D volume at time $t^{n/n+1}$ but span a non-zero 4D volume during their space-time evolution,
		thus contributing to the total 4D volume of each $C_j^n$. 
		
		To sum up, the concepts of standard space-time control volumes generalize automatically to these cases. The only required adjustments are considering the correct set of neighbors $\mathcal{W}_j^n$ and ensuring that the node coordinates for the face neighboring the \textit{hole-like} element are correctly stored and mapped.
		\par
	}

	\section{Direct Arbitrary-Lagrangian-Eulerian ADER discontinuous Galerkin method}
	\label{sec.method-description}
	
	This section is devoted to the description of our ALE ADER-DG numerical method, 
	which allows us to solve any PDE system cast in the form \eqref{eq.generalform} 
	via the integration of its \textit{space-time divergence form}
	\begin{equation}	\label{eq.general-divergenceform}
		\tilde\nabla \cdot \tilde\F(\Q(\x,t))  = \mathbf{0}, \ \text{ with } \ \tilde\nabla = (\nabla, \partial_t) = (\partial_{x_1}, \partial_{x_1}, \partial_{x_1}, \partial_{t}) \ \text{ and } \  \tilde\F(\Q) = (\F(\Q), \Q),
	\end{equation}
	over space-time control volumes connecting 
	the different subsequent tessellations originated from our controlled incremental moving mesh optimization.
	
	To this end, in Section~\ref{s-sec.basis} we begin by introducing the families of basis functions necessary to describe 
	the i)~approximate solution in space and in ii) space-time, 
	as well as the iii) \textit{moving} test functions employed for the weak formulation of the DG scheme,
	which represent a peculiar key ingredient of our ALE scheme on moving meshes. Within an ADER-DG framework, the method is characterized by two main steps: a \textit{predictor} step and a \textit{corrector} step. In Section~\ref{sec.predictor} we describe the local predictor step over each classical control volume $C_i^n$, in Section~\ref{s-sec.predictorcrazy} we introduce a novel locally implicit procedure for the \textit{hole-like} elements and in Section~\ref{sec.corrector} we describe the final explicit corrector step of the method. We close this section by discussing how we deal with the numerical integration.
	
	
	{
		\linespread{1.1}\selectfont
		We also remark that the presented method evolves the solution in an \textit{explicit} way, from $t^n$ to $t^{n+1} = t^n + \Delta t$, 
		by selecting
		\begin{equation}
			\Delta t < \frac{\textnormal{CFL}^N_{\textup{max}}}{d} 
			\frac{ h_i^n }{ \lambda_{\max,i} } 
			\qquad \text{for every } P_i^n \in \mathcal{P}^n,
			\label{eq:timestep}
		\end{equation}
		where $h_i^n$ is a characteristic length of the element $P_i^n$ defined in \eqref{eq.char.length} and
		\[
		\lambda_{\max,i} = 
		\max_{\n \in \R^3, \,||\n ||=1} \left( 
		\max \left(\rho\left( \frac{\partial \F\left(\mathbf{c}_i^n\right)}{\partial \Q} \cdot \n\right),
		\max_{j} \rho\left( \frac{\partial \F\left(\mathbf{d}_{i,j}^n\right)}{\partial \Q} \cdot \n\right),
		\max_{j} \rho\left( \frac{\partial \F\left(\mathbf{b}_{i,j}^n\right)}{\partial \Q} \cdot \n\right) \right) \right),
		\]
		with $\rho(A)$ the maximal spectral radius of a given square matrix $A$.
		Here, following the $1$D ALE ADER-DG stability analysis in \cite{bonafini2026stability}, which also accounts for \textit{hole-like} elements, we set $\textnormal{CFL}^0_{\textup{max}} = 1.0$, $\textnormal{CFL}^1_{\textup{max}}=0.333$, $\textnormal{CFL}^2_{\textup{max}} = 0.170$ and $\textnormal{CFL}^3_{\textup{max}} = 0.104$ (which are just slightly more restrictive than the Runge--Kutta $\CFL$ of $1/(2N+1)$). Moreover, in higher dimension $d > 1$, the maximal allowed $1$D $\CFL$ has to be scaled by $1/d$, see~\cite{dumbser2008unified}. 
		\par
	}
	
	\subsection{Families of basis functions and representation of the approximate solution}
	\label{s-sec.basis}
	
	We represent the conserved variables $\mathbf{Q}$
	via polynomials of degree at most $N$ in an ambient space $\mathbb{R}^d$, 
	which can be either the 3-dimensional space with $d=3$ or the space-time with $d=3+1=4$.
	Hence, for any $N, d \in \N$, we define the set of indices
	\[
	I_N^d = \left\{ \bm{\ell} = (\ell_1, \dots, \ell_d) \in [\{0,\dots,N\}]^d \text{ s.t. } \sum_{k=1}^d \ell_k \leq N \right\} 
	\  \text{ of cardinality } \
	\mathcal{L}_N^d = \frac{1}{d!}\prod_{k=1}^{d}(N+k).
	\]
	We also fix a suitable bijective index mapping
	$
	\ell_N^d \colon  I_N^d \to \{0, \dots, \mathcal{L}_N^d-1\}
	$
	induced by the co-lexicographical ordering, so that
	\[
	\ell_N^{d+1}\left( \left(\ell_1,\dots,\ell_d,0\right) \right) = \ell_N^{d}\left( \left(\ell_1,\dots,\ell_d\right) \right) \quad \text{ for any } d > 1 \text{ and any } (\ell_1, \dots, \ell_d) \in I_N^d.
	\]
	The map $\ell_N^d$ allows us to go back and forth from the set of possible exponents of the polynomial representation and a linear set of indices, 
	which we use to define three different families of basis functions. 

	\paragraph{Spatial basis functions} 
	First, we introduce, over each polyhedron $P_i^n$, the set of spatial basis functions $\{\phi_{i,\ell}^n\}_\ell$, defined~as
	\begin{equation*}
		\phi_{i,\ell}^n \colon P_i^n \to \R, \,\x \mapsto \prod_{k=1}^3 \left(\frac{x_k - (\textbf{c}_i^n)_k}{h_{i}^n} \right)^{\ell_k} \quad \text{for } \ell = \ell_N^3(\bm{\ell}) \in \{0, \dots, \mathcal{L}_N^3-1\}, \bm{\ell} = (\ell_1, \ell_2, \ell_3) \in I_N^3.
	\end{equation*}
	Here, $\textbf{c}_i^n$ is the center of mass of $P_i^n$ and $h_i^n$ is the polyhedral characteristic length, 
	which is computed as the minimum distance between the cell center $\textbf{c}_i^n$ and the face barycenters $\textbf{b}_{ij}^n$ :
	\begin{equation}\label{eq.char.length}
		h_i^n = 2 \cdot \min_{j \, \in \, \mathcal{V}_i^n} \| \ct_i^n - \mathbf{b}_{i,j}^n \|_2 .
	\end{equation}
	
	
	
	\paragraph{Space-time basis functions} 
	Next, we define the set of space-time basis functions $\{\theta_{i,\ell}^n\}_\ell$.
	Over each classical control volume $C_i^n$, $i=1, \dots, N_P$, we have 
	\begin{equation*}
		\begin{aligned}
			&\theta_{i,\ell}^n \colon C_i^n \to \R, \, (\x,t) \mapsto \left[ \prod_{k=1}^3 \left(\frac{x_k - (\textbf{c}_i^n)_k}{h_{i}^n} \right)^{\ell_k} \right] \left( \frac{t-t^n}{\Delta t}\right)^{\ell_4}
			\quad \text{for } \ell = \ell_N^4(\bm{\ell}) \in \{0, \dots, \mathcal{L}_N^4-1\}, \bm{\ell} = (\ell_1, \ell_2, \ell_3, \ell_4) \in I_N^4.
			\\
		\end{aligned}
	\end{equation*}
	Similarly, over each \textit{hole-like} element $H_i^n$, $i > N_P$, we define
	\begin{equation*}
		\begin{aligned}
			&\theta_{i,\ell}^n \colon H_i^n \to \R, \, (\x,t) \mapsto \left[ \prod_{k=1}^3 \left(\frac{x_k - (\textbf{c}_i^n)_k}{h_{i}^n} \right)^{\ell_k} \right] \left( \frac{t-t^n}{\Delta t}\right)^{\ell_4}
			\quad \text{for } \ell = \ell_N^4(\bm{\ell}) \in \{0, \dots, \mathcal{L}_N^4-1\}, \bm{\ell} = (\ell_1, \ell_2, \ell_3, \ell_4) \in I_N^4,
		\end{aligned}
	\end{equation*}
	where $\textbf{c}_i^n$ is computed as the arithmetic average of the coordinates of its generating points, namely: 
	for a 3-2 flip, the vertices of $\bar{f}^{n}$ and $\bar{e}^{n+1}$, 
	for a 2-3 flip, the vertices of $\bar{e}^{n}$ and $\bar{f}^{n+1}$, 
	and for a 4-4 flip, the vertices of $\bar{f}^{n}$ and $\bar{f}^{n+1}$;
	$h_i^n$ is the minimum between the length of $\bar{e}^{n/n+1}$ and/or the incircle of $\bar{f}^{n/n+1}$ (see Figures \ref{fig.3-2hole} and \ref{fig.4-4hole}).

	\paragraph{Representation of the approximate solution}
	
	The conserved variables $\Q$ are represented over each element $P_i^n$ via a piecewise polynomial function $\u_i^n$ in the form
	\begin{equation}\label{eq.un}
		{\u}_i^n(\x)  = \sum_{\ell = 0}^{\mathcal{L}_N^3-1} \phi_{i,\ell}^n(\x) \, \hat{\u}^{n}_{i,\ell} \quad  \text{for } \x \in P_i^n, \,i \in \{1,\dots, N_P\},
	\end{equation}
	where $\{ \hat{\u}^{n}_{i,\ell} \}_{n,i,\ell} \subset \R^\nu$ are vectorial degrees of freedom 
	also denoted as $\hat{\u}_i^n \in [\R^\nu]^{\mathcal{L}_N^3}$.
	
	Moreover, we also make use of approximations in the space-time 4-dimensional domain, 
	thus defined on each control volume $C_i^n$ and on each \textit{hole-like} element $H_i^n$,
	that take the form of space-time polynomials $\q_i^n$ written as
	\begin{equation}\label{eq.qn}
		\q_i^n(\x,t) = \sum_{\ell = 0}^{\mathcal{L}_N^4-1} \theta_{i,\ell}^n (\x, t) \hat{\q}_{i,\ell}^n \quad \text{for }
		\left\{
		\begin{aligned}
			&(\x,t) \in C_i^n &\quad &\text{if } i \in \{1,\dots, N_P\}, \\
			&(\x,t) \in H_i^n &\quad &\text{if } i > N_P,
		\end{aligned}
		\right.
	\end{equation}
	where again the vector of degrees of freedom are given by $\{ \hat{\q}^{n}_{i,\ell} \}_{n,i,\ell} \subset \R^\nu$ also denoted with $\hat{\q}_i^n\in [\R^\nu]^{\mathcal{L}_N^4}$.

	\paragraph{Moving basis functions} 
	Finally, we define an additional set of basis functions, which we refer to as \textit{moving} basis functions, $\{\psi_{i,\ell}^n\}_\ell$, such that
	\begin{equation*}
		\begin{aligned}
			& \psi_{i,\ell}^n \colon C_i^n \to \R, \, , \,(\x, t) \mapsto \prod_{k=1}^3 \left(\frac{x_k - \left(\tilde{\ct}_i^n(t)\right)_k}{h_{i}^n} \right)^{\ell_k} \quad \text{for } \ell = \ell_N^3(\bm{\ell}) \in \{0, \dots, \mathcal{L}_N^3-1\}, \bm{\ell} = (\ell_1, \ell_2, \ell_3) \in I_N^3
			\\[0.6ex]
			& \text{with } \ \tilde \ct_i^n(t) = \left(1-\frac{t-t^n}{\Delta t}\right) \ct_i^{n} + \frac{t-t^n}{\Delta t} \ct_i^{n+1}.
		\end{aligned}
	\end{equation*}
	These functions will be used as \textit{test functions} in the weak formulation of the DG scheme~\eqref{eqn.corrector} (corrector phase) and are the essential ingredient that allows the spatial functions $\{\phi_{i,\ell}^n\}_\ell$ to be always automatically referred to the center $\textbf{c}_i^n$ of~$P_i^n$ at each $n$, without the need to perform an $L_2$ projection to pass from a representation (e.g. $\u_i^n$) referred to $\textbf{c}_i^n$ to the subsequent one (e.g. $\u_i^{n+1}$) referred to~$\textbf{c}_i^{n+1}$.

	\subsection{Direct ALE ADER-DG method: predictor and corrector phases}

	The ADER approach, originally proposed by \cite{toro3,toro2005ader} and used here in the form introduced by Dumbser et al.~\cite{dumbser2008unified}, is an explicit and intrinsically space-time approach to evolve \eqref{eq.generalform} and construct a fully discrete and arbitrarily high order scheme not only in space, but also in time. This method is ideal in a moving mesh context because, when accounting also for the space-time divergence form of the PDE~\eqref{eq.general-divergenceform}, it allows to \textit{directly} pass from time $t^n$ with the solution $\u_i^n$ defined on $P_i^n$, to time $t^{n+1}$ with the newly found solution $\u_i^{n+1}$ already defined on the new cell $P_i^{n+1}$.
	
	The method is classically divided into two phases: the first is called \textit{predictor} (see Section~\ref{sec.predictor}) and is responsible for constructing a high order space-time polynomial $\q_i^n$ as in \eqref{eq.qn} defined on each space-time control volume which, since it is computed without flux exchanges, is valid only within it; this predictor will then be used in the subsequent phase, called \textit{corrector} (see Section~\ref{sec.corrector}), in which the actual update takes place via a \textit{one-step} procedure precisely thanks to the availability of the predictor.
	This two phase procedure is applied over the control volumes $C_i^n$ with~$i~=~1, \dots, N_P$.
	
	The new feature of the scheme proposed here lies in inserting, between the predictor and the corrector phase, 
	a \textit{novel} procedure applied independently to each \textit{hole-like} element (see Section~\ref{s-sec.predictorcrazy}). 
	These elements will be treated via an \textit{implicit space-time DG} method; 
	thus, their evolution will be coupled with that of their neighbors, thereby providing a valid and complete solution of the PDE over them. 
	Nevertheless, our strategy relies on our continuous incremental optimization, 
	which is designed to ensure that \textit{hole-like} elements always have only standard elements as neighbors, 
	for which the predictor can be computed in the classical way and subsequently used as boundary information for the implicit update of each \textit{hole-like} element. 
	This allows the implicit step on each \textit{hole-like} element to involve just a small system of $\nu\mathcal{L}_N^4$ equations for $\nu\mathcal{L}_N^4$ unknowns, that can be solved independently from the others.
	Consequently, while the treatment is locally implicit for the \textit{hole-like} elements, the overall method remains globally explicit.
	
	Finally, we emphasize that the use of this \textit{locally implicit} technology, applied here for the first time to the treatment of \textit{hole-like} elements \textit{in combination with} an \textit{explicit} scheme, makes the method \textit{automatically conservative} even in the vicinity of these degenerate elements, eliminating the need for an \textit{a posteriori} conservative correction which was instead required in previous works by the authors.
	
	\subsubsection{Predictor step on classical space-time control volumes}\label{sec.predictor}
	
	Exploiting the local causality typical of the ADER predictor, 
	for each classical space-time control volume $C_i^n$ with~$i~=~1,\dots, N_P$, 
	we seek for a space-time polynomial $\mathbf{q}_i^n$ of the form \eqref{eq.qn} that satisfies the element-local Cauchy problem
	\begin{equation} \label{eq.Cauchy}
		\left\{ \begin{aligned} 
			& \partial_t \q_i^n(\x,t) + \nabla \cdot \F(\q_i^n(\x,t)) = 0 &\quad& \text{over } C_i^n,	\\
			& \q_i^n(\x,t^n) = \u_i^n(\x) &\quad& \text{on } \ P_i^n.
		\end{aligned} \right.
	\end{equation}
	Here we have imposed that $\mathbf{q}_i^n$ should satisfy the PDE inside $C_i^n$ and, as the sole external condition influencing it, 
	an information coming from the past, namely the known solution on $P_i^n$ at time $t^n$.
	
	To solve~\eqref{eq.Cauchy}, we multiply the PDE by the space-time test functions $\theta_{i,k}^n$ and integrate over the space-time control volume $C_i^n$, obtaining 
	\begin{equation}\label{eqn.pred-step1}
		\int_{C_i^n} \theta_{i,k}^n(\mathbf{x},t) \left[ \partial_t \mathbf{q}_i^n(\mathbf{x},t) + \nabla \cdot \mathbf{F}(\mathbf{q}_i^n(\mathbf{x},t)) \right]\,d\mathbf{x} dt = 0, 
		\quad \text{for } k = 0,\dots, \mathcal{L}_N^4-1.
	\end{equation}
	Then, we rewrite the first term in \eqref{eqn.pred-step1} by taking into account a potential jump of $\q_i^n$ at the boundary $P_i^n$ of $C_i^n$
	via a simplified path-conservative approach~\cite{Pares2006,Castro2006,Castro2008},
	yielding 
	\begin{equation}\label{eqn.pred-step2}
		\begin{aligned} 
			& \int_{C_i^n} \theta_{i,k}^n(\x,t) \partial_t \q_i^n(\x,t) \, d\x dt + 
			\int_{P_i^n} \theta_{i,k}^n(\x,t^n) \, \left( \q_i^{n}(\x,t^n) - \u_i^n(\x) \right) \, d\x  \\ 
			&+ \int_{C_i^n} \theta_{i,k}^n(\x,t) \nabla \cdot \F(\q_i^n(\x,t)) \, d\x dt = \, 0,
			\quad \text{for } k = 0,\dots, \mathcal{L}_N^4-1.
		\end{aligned} 
	\end{equation}
	Now, given $\mathbf{q}_i^n$, we introduce the vector of coefficients $\hat{\mathbf{F}}_{i}^n = (\hat{\mathbf{F}}_{i,0}^n, \dots, \hat{\mathbf{F}}_{i,\mathcal{L}_N^4-1}^n)$, with each $\hat{\mathbf{F}}_{i,\ell}^n \in [\mathbb{R}^\nu]^3$, so that the polynomial representation 
	\begin{equation}\label{eq.fhat}
		\mathbf{F}_i^n(\mathbf{x},t) = \sum_{\ell = 0}^{\mathcal{L}_N^4-1} \theta_{i,\ell}^n (\mathbf{x}, t) \hat{\mathbf{F}}_{i,\ell}^n
	\end{equation}
	is the $L^2$-projection of $\mathbf{F}(\mathbf{q}_i^n)$ onto the finite-dimensional space spanned by the basis functions $\{\theta_{i,0}^n,\dots,\theta_{i,\mathcal{L}_N^4-1}^n\}$.
	By replacing the term $\F(\q_i^n)$ with $\F_i^n$ in~\eqref{eqn.pred-step2}, we get
	\begin{equation}\label{eqn.predictor}
		\begin{aligned} 
			& \int_{C_i^n} \theta_{i,k}^n(\x,t) \partial_t \q_{i}^n(\x,t) \, d\x dt  + 
			\int_{P_i^n} \theta_{i,k}^n(\x,t^n) \q_{i}^n(\x,t^n) \,d\x \\		
			&= \int_{P_i^n} \theta_{i,k}^n(\x,t^n) \u_{i}^n(\x) \, d\x
			- \int_{C_i^n} \theta_{i,k}^n(\x,t) \nabla \cdot \F_{i}^n(\x,t)\, d\x dt,  \quad \text{for } k = 0,\dots, \mathcal{L}_N^4-1.  
		\end{aligned} 
	\end{equation}
	By expanding $\q_i^n$ and $\u_i^n$ in \eqref{eqn.predictor} using respectively~\eqref{eq.qn} and \eqref{eq.un}, we obtain
	\begin{equation*}
		\begin{aligned} 
			& \sum_{\ell = 0}^{\mathcal{L}_N^4-1 } \left[ \int_{C_i^n} \theta_{i,k}^n(\x,t) \partial_t \theta_{i,\ell}^n (\x, t)  \, d\x dt \right] \hat{\q}_{i,\ell}^n  + 
			\sum_{\ell = 0}^{\mathcal{L}_N^4-1 } \left[ \int_{P_i^n} \theta_{i,k}^n(\x,t^n) \theta_{i,\ell}^n (\x, t^n) \,d\x \right] \hat{\q}_{i,\ell}^n \\		
			&= \sum_{\ell = 0}^{\mathcal{L}_N^3-1 } \left[ \int_{P_i^n} \theta_{i,k}^n(\x,t^n)  \phi_{i,\ell}^n(\x) \, d\x \right] \hat{\u}^{n}_{i,\ell}
			- \sum_{j=1}^3 \sum_{\ell = 0}^{\mathcal{L}_N^4-1 }  \left[ \int_{C_i^n} \theta_{i,k}^n(\x,t) \partial_{x_j}\theta_{i,\ell}^n(\x,t) \, d\x dt \right] \left[\hat{\F}_{i,\ell}^n\right]_j, \quad \text{for } k = 0,\dots, \mathcal{L}_N^4-1, 
		\end{aligned} 
	\end{equation*}
	where the coefficients $\hat{\F}_{i,\ell}^n$ are a generally nonlinear function of $\hat{\q}_{i}^n$, i.e. $\hat{\F}_{i,\ell}^n = \hat{\F}_{i,\ell}^n\left(\hat{\q}_{i}^n\right)$,
	thus yielding an implicit system of equations for the unknown coefficient vector $\hat \q_i^n$.
	Hence, to approximate $\hat{\q}_{i}^n$ we employ a fixed point Picard iteration, as detailed in~\cite{dumbser2008unified,hidalgo2011ader,busto2020high}: for each $r>0$, one solves for $\hat \q_i^{n,(r+1)}$ the linear system
	\begin{equation*}
		\begin{aligned} 
			& \sum_{\ell = 0}^{\mathcal{L}_N^4-1 } \left[ \int_{C_i^n} \theta_{i,k}^n(\x,t) \partial_t \theta_{i,\ell}^n (\x, t)  \, d\x dt \right] \hat{\q}_{i,\ell}^{n,(r+1)}  + 
			\sum_{\ell = 0}^{\mathcal{L}_N^4-1 } \left[ \int_{P_i^n} \theta_{i,k}^n(\x,t^n) \theta_{i,\ell}^n (\x, t^n) \,d\x \right] \hat{\q}_{i,\ell}^{n,(r+1)} \\		
			&= \sum_{\ell = 0}^{\mathcal{L}_N^3-1 } \left[ \int_{P_i^n} \theta_{i,k}^n(\x,t^n)  \phi_{i,\ell}^n(\x) \, d\x \right] \hat{\u}^{n}_{i,\ell}
			- \sum_{j=1}^3 \sum_{\ell = 0}^{\mathcal{L}_N^4-1 }  \left[ \int_{C_i^n} \theta_{i,k}^n(\x,t) \partial_{x_j}\theta_{i,\ell}^n(\x,t) \, d\x dt \right] \left[\hat{\F}_{i,\ell}^n\left(\hat{\q}_{i}^{n,(r)}\right)\right]_j, \quad \text{for } k = 0,\dots, \mathcal{L}_N^4-1.  
		\end{aligned} 
	\end{equation*}
	Here, we can consider the vector $\hat \q_i^{n,(0)}:=(\hat \u_{i,0}^n, \dots, \hat \u_{i,\mathcal{L}_N^3-1}^n, 0, \dots, 0)$ as starting point for the iteration; moreover, we remark that such fixed point procedure has already been proved to be convergent and to yield the desired order of 
	accuracy, see~\cite{jackson2017eigenvalues,busto2020high,han2021dec} for more details.
	
	Upon convergence, the predictor step provides us with a set of locally defined high order polynomials $\q_i^n$ over each space-time control volume $C_i^n$. These polynomials will serve as approximant of the solution in the interior of $C_i^n$ and will be used in the next phases for the computation of the numerical fluxes at the interfaces.
	
	\subsubsection{Locally implicit space-time DG method on \textit{hole-like} elements}
	\label{s-sec.predictorcrazy}
	
	Once the predictors for all classical space-time control volumes have been computed as described above, 
	we consider the \textit{hole-like} elements $H_{i}^n$ for $i > N_P$. 
	We recall, in particular, that due to how we organized our continuous incremental mesh optimization, 
	we produce just one \textit{hole-like} element per neighborhood, 
	and consequently they have as neighbors exclusively classical elements, 
	on which, at this stage, we have already computed the predictors. 
	
	Since the \textit{hole-like} elements are spatially degenerate, 
	coinciding just with a face or an edge of the tessellation instead of a real element, 
	we do not have an inflow information to impose as initial value for the Cauchy problem~\eqref{eq.Cauchy} at time $t^n$, 
	thus we cannot apply the previous strategy over them. 
	Furthermore, as noted in previous works by the authors~\cite{gaburro2020high,gaburro2025high}, the corrector step would also not be well-defined for these elements since they are spatially degenerate even at time $t^{n+1}$; 
	moreover, such corrector phase is unnecessary for the \textit{hole-like} elements which, living only in between two time steps, 
	do not need to produce any approximation of the solution in the form~\eqref{eq.un} at $t^{n+1}$.
	
	Thus, for \textit{hole-like} elements, only a space-time representation $\mathbf{q}_i^n$ of the form~\eqref{eq.qn} is required within $H_i^n$ for $i > N_P$, so that it can be used by neighboring cells during their corrector step. 
	Crucially, this information must be high order accurate and ensure the conservativity of the method. 
	A key novelty of the approach presented here is that such properties are achieved intrinsically, 
	without requiring any \textit{a posteriori} correction. 
	To obtain this result, the procedure for $\mathbf{q}_i^n$ for $i > N_P$, unlike the predictor typically employed for computing these space-time quantities, must already account for the space-time information available in the neighbors, 
	exactly as the corrector will later do for the standard~elements. 
	
	Hence, we derive our scheme for the \textit{hole-like} elements as an \textit{implicit} space-time DG scheme, 
	which is however implicit only \textit{element-wise}, locally over each \textit{hole-like} element.
	This is because all the space-time information from the neighboring classical control volumes are already known and available, 
	thus decoupling each region needing an implicit treatment from the others and making the update of each \textit{hole-like} element independent. 
	
	In detail, for each $i>N_P$, we multiply the space-time divergence form of the governing PDE~\eqref{eq.general-divergenceform} by a spacetime basis function $\theta_{i,k}^n$, $k = 0,\dots,\mathcal{L}_N^4-1$, and integrate over the \textit{hole-like} element $H_i^n$ 
	\begin{equation*}
		\int_{H_i^n} \theta_{i,k}^n(\x,t) \; \tilde\nabla \cdot \tilde\F(\Q(\x,t)) \,d\x dt = 0 \quad \text{for } k = 0,\dots, N.
	\end{equation*}
	Next, integrating by parts, we get
	\begin{equation*}
		\int_{\partial H_i^n} \theta_{i,k}^n \tilde\F(\Q) \cdot \n_i^n \,ds - \int_{H_i^n} \tilde\nabla \theta_{i,k}^n(\x,t) \cdot \tilde\F(\Q(\x,t)) \,d\x dt = 0\quad \text{for } k = 0,\dots, N,
	\end{equation*}
	where $\mathbf{n}_i^n$ denotes the outward pointing unit normal vector on the spacetime surfaces composing the boundary $\partial H_i^n$.
	Now, we consider that the PDE solution is represented, both on $H_i^n$ and 
	on its neighbors $C_j^n$, $j \in \mathcal{W}_i^n$, with polynomials $\q_{i}^n$ and $\q_{j}^n$ in the form \eqref{eq.qn}, 
	with all the $\q_j^n$ already available and the coefficients $\hat{\q}_i^n$ as the only unknowns; 
	also, similarly to~\eqref{eq.fhat}, we introduce the polynomial representation  $\tilde\F_i^n = [\F_{i}^n, \q_{i}^n]$.
	In this manner,  we obtain
	\begin{equation}\label{eqn.qnsliver}
		\begin{aligned}
			& \sum_{j \, \in \,  \mathcal{W}_i^n} \int_{\partial H_{ij}^n}
			\theta_{i,k}^n \mathcal{F}(\q_i^n,\q_j^n, \mathbf{n}_{i}^n)\, ds
			= \int_{H_{i}^n} 
			\tilde\nabla \theta_{i,k}^n(x,t) \cdot \tilde \F_{i}^n(\x,t) \,d\x dt 
			\quad \text{for }k = 0, \dots, \mathcal{L}_N^4-1,
		\end{aligned}
	\end{equation}
	where $\mathcal{F} \colon \R\times \R \times S^3 \to \R$ is a numerical flux function that we 
	compute via an ALE Riemann solver applied to the inner and outer boundary-extrapolated data $\q_L, \q_R \in \R^\nu$ at the interfaces.
	Here, the simplest choice consists in adopting a Rusanov-type~\cite{Rusanov:1961a} ALE flux: for $\q_L,\q_R \in \R^\nu$ and $\mathbf{n} = (\mathbf n_\x, \mathbf n_t) = (\n_{x_1}, \n_{x_2}, \n_{x_3}, \n_{t})$ with $||\mathbf n||=1$, we set
	\begin{equation}
		\label{eq.rusanov} 
		\mathcal{F}(\q_L,\q_R,\mathbf{n})  =  
		\frac{1}{2} \left( \F(\q_L) + \F(\q_R)  \right) \cdot \mathbf n_\x
		+ \frac{1}{2} \left( \q_L + \q_R  \right) \mathbf n_t
		- \frac{1}{2} s_{\max} \left( \q_R - \q_L \right),  
	\end{equation} 
	where $s_{\max}$ is the maximum of the spectral radii of the ALE Jacobian matrices $\mathbf A(\q_L, \mathbf{n})$ and $\mathbf A(\q_R, \mathbf{n})$, with
	\[
	\mathbf A(\q, \mathbf{n}) = \left[ \frac{\partial\F(\q)}{\partial\Q} \cdot \mathbf{n}_\x - (\mathbf{v} \cdot \mathbf{n}_\x) \mathbf{I} \right],
	\]
	where $\mathbf{I}$ is the $\nu\times \nu$ identity matrix and $\mathbf{v}$ is the local grid velocity, obtained as $\mathbf{v} = -\mathbf{n}_t/\|\mathbf{n}_\x\|$.
	
	By expanding in~\eqref{eqn.qnsliver} the unknown $\q_i^n$ according to~\eqref{eq.qn}, we get
	\begin{equation}
		\label{eqn.qnsliver-expanded}
		\begin{aligned}
			& \sum_{j \, \in \,  \mathcal{W}_i^n} \int_{\partial H_{ij}^n}
			\theta_{i,k}^n \mathcal{F}\left(\sum_{\ell = 0}^{\mathcal{L}_N^4-1} \theta_{i,\ell}^n \hat{\q}_{i,\ell}^n,\q_j^n, \mathbf{n}_{i}^n\right)\, ds
			= \sum_{\ell = 0}^{\mathcal{L}_N^4-1}
			\left[
			\int_{H_{i}^n} \partial_t \theta_{i,k}^n(\x,t)  \theta_{i,\ell}^n (\x, t)  \,d\x dt 
			\right]
			\hat{\q}_{i,\ell}^n \\
			&+ \sum_{j=1}^3 \sum_{\ell = 0}^{\mathcal{L}_N^4-1}
			\left[
			\int_{H_{i}^n} \partial_{x_j} \theta_{i,k}^n(\x,t) \theta_{i,\ell}^n (\mathbf{x}, t) \,d\x dt 
			\right]
			\left[\hat{\F}_{i,\ell}^n\left(\hat{\q}_{i}^n\right) \right]_j
			\quad \text{for }k = 0, \dots, \mathcal{L}_N^4-1.
		\end{aligned}
	\end{equation}%
	{%
		\linespread{1.1}\selectfont%
		This is a non-linear system of $\nu\,\mathcal{L}_N^4$ equations that \textit{implicitly} define the unknown coefficients vector $\hat\q_i^n \in \R^{\nu\,\mathcal{L}_N^4}$. Here, we recall once again that each \textit{hole-like} element is surrounded exclusively by classical space-time control volumes, hence, all the neighbours predictors $\{\q_j^n\}_{j\in\mathcal{W}_i^n}$ have already been computed in the previous step of the algorithm, leaving only $\hat\q_i^n$ to be determined. To do so, we solve \eqref{eqn.qnsliver-expanded} for the unknown $\hat\q_i^n$ by means of a Newton method, where the initial guess is computed as a weighted average of the predictors of the neighbouring cells $\{\hat\q_j^n\}_{j\in\mathcal{W}_i^n}$. Each step of the algorithm only requires the inversion of a \textit{local} $\left(\nu\,\mathcal{L}_N^4 \times \nu\,\mathcal{L}_N^4\right)$ linear system involving the Jacobian matrix of~\eqref{eqn.qnsliver-expanded} with respect to $\hat\q_i^n$. To build the Jacobian matrix we differentiate explicitly both the flux function $\mathcal{F}$ (with respect to its first argument) and the projection operator $\hat\q_i^n \mapsto \hat\F_i^n$ \eqref{eq.fhat}. Note that the only non-smooth operation in~\eqref{eqn.qnsliver-expanded} is introduced when computing $s_{\max}$ in~\eqref{eq.rusanov}, for which we can replace the maximum function by a soft-max.
		\par
	}
	
	To sum up, solving~\eqref{eqn.qnsliver-expanded} yields a high order accurate approximation of the solution $\mathbf{q}_i^n$ over the \textit{hole-like} element~$H_i^n$, which can be used in the subsequent phases. In particular, the first $\nu$ equations in~\eqref{eqn.qnsliver-expanded}, those corresponding to $k = 0$, i.e. $\theta_{i,0}^n(\x,t) \equiv 1$, ensure that the overall flux balance on $H_i^n$ is actually zero. This, combined with the corrector step here below, guarantees that the method is \textit{conservative} by construction.
	\RIcolor{
		\paragraph{Heuristic on the practical behavior of our local implicit solver on \textit{hole-like} elements}
		As observed during the numerical simulations presented hereafter, the local Newton iteration consistently exhibits quadratic convergence, requiring on average only three iterations to achieve the desired tolerance (in general, 1.0E-10). 
		We found that our neighbor predictor initial guess leads to convergence in one iteration less than a non-tailored random initial guess, however the method appears to be robust for any initial guess within the phase space. Finally, the local Jacobian matrix of the system in equation \eqref{eqn.qnsliver-expanded} remains well-conditioned throughout all our simulations, indicating that the implicit step is stable and reliable in practice.
	}
	
	\subsubsection{Final solution update via the corrector step}
	\label{sec.corrector}
	
	Now that a high order space-time approximation is available for each control volume, we can proceed with the final update of the solution,
	which is an operation required only for those elements that indeed exist at time $t^{n+1}$, namely those with $i \in \{1,\dots,N_P\}$.
	
	For each $i \in \{1,\dots,N_P\}$, we multiply the space-time divergence form of the governing PDE~\eqref{eq.general-divergenceform} by the \textit{moving} basis function $\psi_{i,k}^n$, $k = 0,\dots,\mathcal{L}_N^3-1$, and integrate over the control volume $C_i^n$ to obtain
	\begin{equation*}
		\int_{C_i^n} \psi_{i,k}^n(\x,t) \; \tilde\nabla \cdot \tilde\F(\Q(\x,t)) \,d\x dt = 0 \quad \text{for } k = 0,\dots, \mathcal{L}_N^3-1.
	\end{equation*}
	Next, integrating by parts, we get
	\begin{equation*}
		\int_{\partial C_i^n} \psi_{i,k}^n \tilde\F(\Q) \cdot \mathbf n_{i}^n \,ds - \int_{C_i^n} \tilde\nabla\psi_{i,k}^n(\x,t) \cdot \tilde{\F}(\Q(\x,t)) \,d\x dt = 0\quad \text{for } k = 0,\dots, \mathcal{L}_N^3-1,
	\end{equation*}
	where $\mathbf{n}_i^n$ denotes the outward pointing unit normal vector on the spacetime surfaces composing the boundary $\partial C_{i}^n$. We now recall i) the decomposition~\eqref{eqn.dC} of $\partial C_{i}^n$ into $P_i^n$, $P_i^{n+1}$ and the lateral faces $\partial C_{ij}^n$, ii) that $\u_i^n$ (respectively $\u_i^{n+1}$) approximates $\Q$ on $P_i^n$ (respectively $P_i^{n+1}$) and iii) that the space-time polynomials $\q_i^n$ approximates $\Q$ inside the control volumes $C_i^n$ and $H_i^n$. 
	Hence, using again the numerical flux~\eqref{eq.rusanov} we get
	\begin{equation}
		\label{eqn.corrector}
		\begin{aligned} 
			&\int_{P_i^{n+1}} \psi_{i,k}^n(\x,t^{n+1}) \u_{i}^{n+1}(\x) \, d\x = 
			\int_{P_{i}^n} \psi_{i,k}^n(\x,t^n) \u_{i}^{n}(\x) \, d\x
			+ \int_{C_i^n} \tilde\nabla\psi_{i,k}^n(\x,t) \cdot \tilde{\F}(\q_i^n(\x,t)) \,d\x dt \\
			& - \sum_{j \in \mathcal{W}_i^n} \int_{\partial C_{ij}^n}
			\psi_{i,k}^n \mathcal{F}(\q_i^n,\q_j^n, \mathbf{n}_{i}^n)\, ds  \quad  \text{for }k = 0, \dots, \mathcal{L}_N^3-1,
		\end{aligned} 
	\end{equation}
	or equivalently with the expansion~\eqref{eq.un}
	\begin{equation*}
		\begin{aligned} 
			&\sum_{\ell = 0}^{\mathcal{L}_N^3-1} \left[ \int_{P_i^{n+1}} \psi_{i,k}^n(\x,t^{n+1}) \phi_{i,\ell}^{n+1}(\x) \, d\x \right]  \hat{\u}^{n+1}_{i,\ell}  = 
			\sum_{\ell = 0}^{\mathcal{L}_N^3-1} \left[ \int_{P_{i}^n} \psi_{i,k}^n(\x,t^n)  \phi_{i,\ell}^n(\x) \, d\x \right] \hat{\u}^{n}_{i,\ell}
			+ \int_{C_i^n} \partial_t \psi_{i,k}^n(\x,t) \q_i^n(\x,t) \,d\x dt \\
			& +\sum_{j=1}^3 \int_{C_i^n} \partial_{x_j} \psi_{i,k}^n(\x,t) \q_i^n(\x,t) \,d\x dt - \sum_{j \in \mathcal{W}_i^n} \int_{\partial C_{ij}^n}
			\psi_{i,k}^n \mathcal{F}(\q_i^n,\q_j^n, \mathbf{n}_{i}^n)\, ds  \quad  \text{for }k = 0, \dots, \mathcal{L}_N^3-1,
		\end{aligned} 
	\end{equation*}
	where the coefficient vector $\hat \u_i^{n+1}$ can be computed {explicitly} upon knowing $\hat \u_i^n$ and by integrating the remaining terms that depend only on the already available predictors~$\q_i^n$.
	
	\subsection{Numerical integration and quadrature formulas}
	\label{s-sec.quadrature}
	
	All integrals involved in our ADER-DG scheme are approximated via Gaussian quadrature rules \cite{stroud}. The number of Gaussian points is in general determined to obtain an exact result whenever integrating polynomials of degree at least up to~$2N$ on each domain of integration. In particular, each of our domains can be obtained as (possibly a union or family of) tetrahedra or triangular prisms. Hence, we start by fixing the Gaussian quadrature on these two reference elements:
	\begin{itemize}
		\item[i)] over the reference tetrahedron with vertices $\{(0,0,0), (1,0,0), (0,1,0), (0,0,1)\}$ we consider $\max(8, (N+1)^3)$ Gaussian quadrature points;
		\item[ii)] over the reference triangular prism with vertices $\{(0,0,0), (1,0,0), (0,1,0), (0,0,1), (1,0,1), (0,1,1)\}$ we consider $\max(2, N+1) \cdot \max(4, (N+1)^2)$ Gaussian quadrature points, obtained as tensor product of the $\max(4, (N+1)^2)$ Gaussian quadrature points on each triangular face and the $\max(2, N+1)$ Gaussian quadrature points along the height dimension.
	\end{itemize}
	Then, each element in the physical domain is derived as the image of one of these two reference elements under a suitable transformation map. We now describe how to handle each physical element explicitly.

	\smallskip
	\noindent
	\textit{Polyhedral domains $P_i^n$.}
	We look at $P_i^n$ as a union of sub-tetrahedra, each obtained by connecting a triangular facet of $\mathcal{F}_i^n$ with $\mathbf{c}_i^n$, and integrate over each sub-tetrahedron separately.
	
	\smallskip
	\noindent
	\textit{Classical space-time control volumes $C_i^n$.}
	As in Figure~\ref{fig.controlvolumes}, we look at each 3D time-slice at time $t^{n+\tau}$ of the 4D volume~$C_i^n$ as a union of sub-tetrahedra, obtained as the convex combination of the corresponding sub-tetrahedra in $P_i^n$ and $P_i^{n+1}$. We now fix $N_\tau = \max(2,N+1)$ Gaussian quadrature points $t^n \leq t^{n+\tau_1} < \ldots < t^{n+\tau_{N_\tau}} \leq t^{n+1}$ in time, and for each of them we compute the 3D volume integral of the corresponding 3D time-slice. Then, we sum up these contributions weighted according to the Gaussian weights in time.
	
	\smallskip
	\noindent
	\textit{Classical space-time lateral surfaces $\partial C_{ij}^n$.} As described in Section \ref{s-sec.3dlateralsurface}, each classical space-time lateral surface $\partial C_{ij}^n$ can be decomposed as a union of 3D curvilinear  surfaces obtained as images of the reference triangular prism, over which we integrate separately.
	
	\smallskip
	\noindent
	\textit{Hole-like elements $H_i^n$.}
	In the case of a 3-2 flip, as shown in Figure~\ref{fig.3-2hole}, each 3D time-slice at time $t^{n+\tau}$ of the 4D volume~$H_i^n$ is described as a triangular prism. Hence, we fix $N_\tau = \max(2,N+1)$ Gaussian quadrature points $t^n \leq t^{n+\tau_1} < \ldots < t^{n+\tau_{N_\tau}} \leq t^{n+1}$ in time, and for each of them we compute the 3D volume integral of the corresponding 3D triangular prism. Then, we sum up these contributions weighted according to the Gaussian weights in time. For a 2-3 flip the same holds true, while for a 4-4 flip, according to Figure~\ref{fig.4-4hole}, we can see each 3D time-slice of $H_i^n$ as the union of two triangular prisms, and the same procedure can be implemented.
	
	\smallskip
	\noindent
	\textit{Hole-like space-time lateral surfaces $\partial H_{ij}^n$.}
	As described in the right panel of Figures~\ref{fig.3-2holesurfaces} and \ref{fig.4-4holesurface}, each part of the lateral surface $\partial H_{ij}^n$ can be viewed as the image of a 3D triangular prism under a suitable transformation map, possibly shrinking the lower or upper triangular basis to a point or a line. Then we integrate over each of them separately.
	

	\section{Numerical results} \label{sec_num-results}	
	
	To perform our numerical benchmarks, we choose to work with the Euler equations of compressible gas dynamics which represent a well-known system of hyperbolic equations that can be cast in the form~\eqref{eq.generalform} by taking $\Q$ and $\F$ as follows
	\begin{equation*}
		\Q = \left( \begin{array}{c} \rho   \\ \rho u  \\ \rho v \\ \rho w   \\ \rho E \end{array} \right), \quad
		\mathbf{F} = \left( 
		\begin{array}{c} 
			\rho u \\ 
			\rho u^2 + p \\ 
			\rho v u \\ 
			\rho w u \\ 
			u(\rho E + p) 
		\end{array} , 
		\begin{array}{c} 
			\rho v \\ 
			\rho u v \\ 
			\rho v^2 + p \\ 
			\rho w v \\ 
			v(\rho E + p) 
		\end{array} , 
		\begin{array}{c} 
			\rho w \\ 
			\rho u w \\ 
			\rho v w \\ 
			\rho w^2 + p \\ 
			w(\rho E + p) 
		\end{array} 
		\right).
	\end{equation*}
	The vector $\Q = \Q(\x,t)$ of the conserved variables includes the fluid density $\rho$, the momentum  
	vector $(\rho u, \rho v, \rho w)$ and the total energy density $\rho E$. 
	The system is closed with the ideal gas equation of state which relates the 
	fluid pressure $p$ to $\Q$ by
	\begin{equation*}
		p = (\gamma-1)\left(\rho E - \frac{1}{2} \rho (u^2 + v^2 + w^2) \right), 
	\end{equation*}
	where $\gamma$ is the ratio of specific heats so that the speed of sound takes the form $c=\sqrt{\frac{\gamma p}{\rho}}$.

	The objective of the presented test cases is to verify that the fundamental properties of a high order conservative scheme are still satisfied in the context of our moving meshes and in the presence of single or multiple topology changes.	We consider two distinct scenarios.
	
	First, we focus on a series of \textit{sanity checks} where we provide the minimum geometry required to trigger an elementary flip, thus to construct a single \textit{hole-like} element, and then we perform one step of the method. 
	We numerically demonstrate that our scheme is i) volume and ii) mass conservative, 
	iii) exactly satisfies the Geometric Conservation Law (GCL), 
	iv) preserves constant states exactly for any order of accuracy, 
	and that v) our ALE ADER-DG scheme, using polynomials of degree $N$, exactly preserves stationary initial conditions where the density is represented by polynomials of the same degree $N$, a property that is numerically verified for each considered $N$ (e.g., $P_1$ preserves linear profiles exactly, $P_2$ preserves linear and quadratic profiles exactly, etc.). 
	Subsequently, these tests are repeated for a more extended rigid rotation that entails multiple topology changes.
	
	Secondly, we employ the classical \RIcolor{moving} isentropic Shu-vortex type benchmark to verify the robustness of the scheme under rotational flows. 
	This setup requires the generation of numerous \textit{hole-like} elements throughout the simulation. 
	We show that the method not only maintains all the relevant properties listed above but also vi) achieves the expected~$N+1$ order of accuracy over long times after the formation of many \textit{hole-like} elements, and vii) the mesh movement does not deteriorate the timestep size. 
	Furthermore, we compare our results with those obtained in a classical direct ALE setting without topology changes, which fails to complete even a quarter of a rotation. 
	Finally, we provide statistics regarding the number of \textit{hole-like} elements generated per timestep in relation to the total number of elements \RIIcolor{and the computational costs of each part of our algorithm}.
	
	\RIcolor{Finally, we close this section with a brief discussion acknowledging the pros and cons of our direct, remap-free, ALE algorithm compared with classical indirect ALE approaches.}

	\subsection{Sanity checks: single elementary flips}
	
\begin{figure}[!t]  
		\centering
		\includegraphics[width=1.0\linewidth]{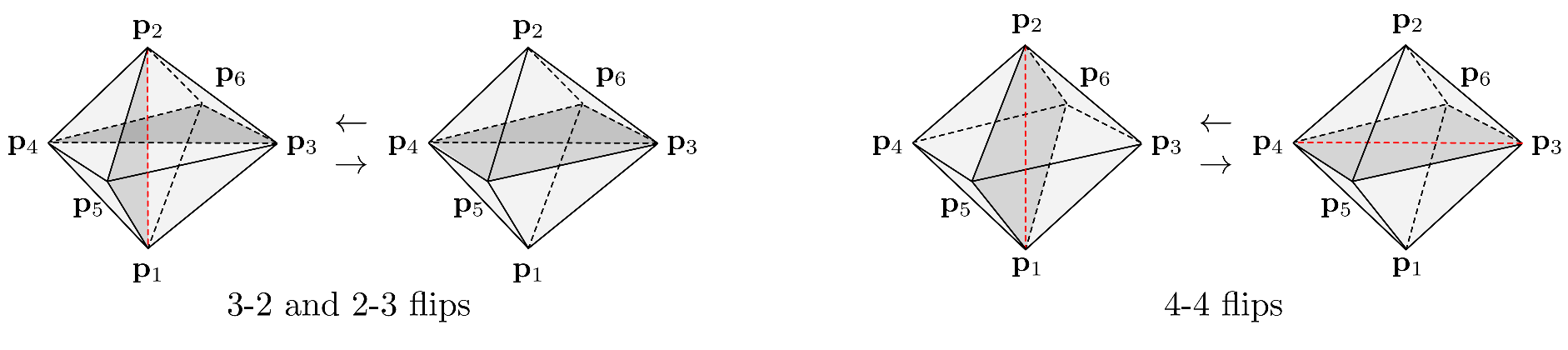}		
		\caption{Connectivity configurations of the diamond region of the tetrahedralization corresponding to elementary flips. 
			The left side of the figure illustrates the configuration for a 3-2 flip (reading from left to right) and a 2-3 flip (reading from right to left); for these cases, one may focus solely on the first five generator points to visualize the transformation. 
			We also refer the reader to Figure~\ref{fig.3-2flip-origin}, where this configuration was previously introduced. 
			Additionally, the two right panels show the configuration for a 4-4 flip, for which we refer the reader to Figure~\ref{fig.4-4flip-origin}.
		}
		\label{fig.san.gen}
	\end{figure}

	In this section, we focus on the study of a single flip at a time. 
	To create the simplest mesh capable of representing elementary flip transformations, 
	we consider a cubic domain $\Omega^0 = [-1,1]^3$ and we place $14$ generator points at these specific locations  
	\begin{equation*}
		\mathbf{p}_{1 \ldots 14} = \begin{pmatrix}
			0   & 0     & 1   & -1    &  0   & 0    & -1   & 1    & 1    & -1   & -1   & 1    & 1    & -1   \\
			y_s & y_s   & 0   &  0    & -1   & 1    & -1   & -1   & 1    &  1   & -1   & -1   & 1    &  1   \\
			-1   & 1     & 0   &  0    &  0   & 0    & -1   & -1   & -1   & -1   &  1   &  1   & 1    &  1   
		\end{pmatrix}.
	\end{equation*}		
	The $y$-axis offset is set to $y_s = -0.5$ for the 2-3 and 3-2 flips, and to $y_s = 0$ for the 4-4 flip. 
	In particular, in Figure~\ref{fig.san.gen}, we show the \textit{diamond} region of the tetrahedralization delimited by the points $\mathbf{p}_i, i =1, \ldots, 6$, where the elementary flips are performed. 
	Indeed, we construct the central tetrahedra as described in the figure, changing their connectivity from the configuration on the left at the beginning of the timestep to the one on the right at the end of the timestep (and \textit{vice versa}), whereas the rest of the domain is meshed regularly and remains fixed.
	
\begin{figure}[!t]
		\centering
		\includegraphics[width=1.0\linewidth]{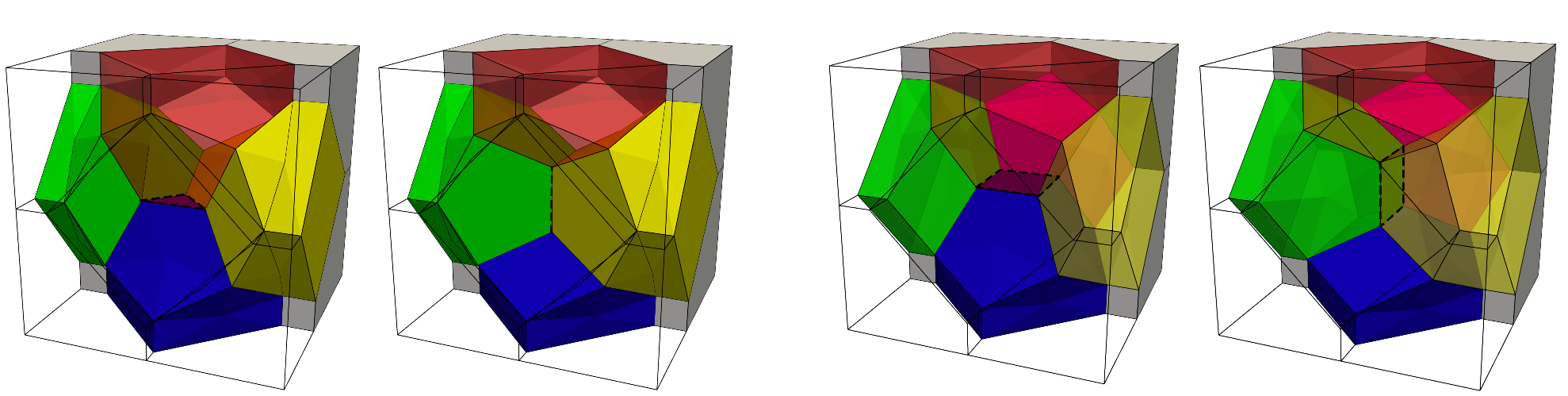}		
		\caption{The left part of this figure shows the minimal computational polyhedral mesh corresponding to a 3-2 flip, before (first image) and after (second image) the topology change. 
			Colored cells are linked to some of the generators of the diamond region; namely, the blue one corresponds to $\textbf{p}_1$, the red one to $\textbf{p}_2$, the yellow one to $\textbf{p}_3$, and the green one to $\textbf{p}_4$. 
			The front cell, corresponding to $\textbf{p}_5$, is not shown to make the transformation clearly visible. The eight corner cells corresponding to the remaining generators do not change their connectivity during the transformation.
			Similarly, the right part shows the minimal computational polyhedral mesh for a 4-4 flip, before (third image) and after (fourth image) the topology change.
		}
		\label{fig.san.flip23e44}
	\end{figure}

	Moreover, in Figure~\ref{fig.san.flip23e44} we report the corresponding effect on the polyhedral tessellation using the same color map of Figures~\ref{fig.3-2hole} and~\ref{fig.4-4hole}, specifically emphasizing the surfaces and edges removed and/or created by the topology changes. 
	These objects, shown with black dashed lines in the figures, represent the configurations with zero 3D volume into which our 4D \textit{hole-like} elements degenerate at the beginning and end of the timestep. Through these \textit{hole-like} elements, such objects are formally connected to each other, linking all surrounding polyhedra and enabling the correct exchange of fluxes between them.
	
	These simple setups represent the building blocks of our numerical scheme; 
	therefore, we use them to verify that all fundamental properties of an ALE ADER-DG scheme with $N=0,1,2,3$ (i.e., up to \textit{fourth} order) are satisfied over a single timestep with $\Delta t = \frac{0.1}{2N+1}$ while forcing a single flip.
	
	We note that 	
	the 3D volume of the cube is $ |\Omega^0| = |\Omega^1| = 8$ and 
	the expected 4D space-time volume is $|\Omega^0| \cdot \Delta t = \frac{0.8}{2N+1}$.
	We consider a stationary solution with various density profiles $\rho$ and with $(u,v,w,p) = (0,0,0,1)$ throughout, and we impose wall boundary conditions.

	
	\paragraph{3-2 flip} In the case of a 3-2 flip, we first verify that $|\Omega^1| = |\Omega^0|$ up to machine precision
	
	\begin{equation*}
		\sum_{i=1}^{14} \int_{P_i^1} 1 \, d\mathbf{x}  = \sum_{i=1}^{14} 
		\int_{P_i^0} 1 \, d\mathbf{x}, \quad \text{for }\, N=1,2,3,\Rall{4}.\\
	\end{equation*} 
	Next, we sum the 4D space-time volumes of all the classical control volumes $C_i^0, i = 1, \ldots, 14$, together with the 4D space-time volume of the unique \textit{hole-like} element $H_{15}^0$ required for this configuration, obtaining	
	\begin{equation*}
		\sum_{i=1}^{14} \int_{C_i^0} 1 \, d\mathbf{x} \, dt + \int_{H_{15}^0} 1 \, d\mathbf{x} \, dt = 
		\frac{\text{7.9974E-01}}{2N+1} +  \frac{\text{2.6042E-04}}{2N+1} = \frac{\text{0.8}}{2N+1}   \quad \text{for }\, N=1,2,3,\Rall{4},
	\end{equation*}	
	which proves that, thanks to the presence of the \textit{hole-like} element, we are indeed able, for each $N$, to cover the entire space-time domain and conserve the total \textit{volume}. 
	
	Next, we compute the difference between the integrals over the spatial control volumes at the beginning and the end of the timestep for various density profiles \RIIcolor{(including non polynomial ones)} and each degree $N$; the results are reported in Table~\ref{tab:conservation-errors-3-2}. This table demonstrates that \textit{mass} is conserved up to machine precision during the PDE evolution for all $N$. We emphasize that, although only a single timestep is performed here, all phases of the numerical method described in Section~\ref{sec.method-description} are fully executed and contribute to the final result.
	
	Furthermore, we verify that our scheme, using polynomials of degree $N$, \textit{exactly preserves} stationary initial conditions when the density is represented by \textit{polynomials} up to degree $N$, as shown in Table~\ref{tab:conservation-errors-3-2}. 
	It should be noted that achieving exactness on polynomial profiles of degree $N$ requires the quadrature rules employed on each domain and for each involved integrand to be exact, not just to have the correct order of accuracy. 
	In Section~\ref{s-sec.quadrature}, we have specified the number of quadrature points necessary to ensure exactness up to degree $2N$, which is sufficient to guarantee exactness on constants and the correct order of convergence; this is the standard number of points used throughout this work. 
	Specifically here, only for the verification of this exactness property, we have increased the number of quadrature points used for the flux computation on the 3D lateral surfaces to $\max(2, N+2) \cdot \max(4, (N+2)^2)$. This is because, the integrand here contains not only the product of two basis functions (which can reach degree up to $2N$) but also the Jacobian determinant of the transformation from the reference element through our set of bi-linear basis functions, see~(\ref{eq.surfaceintegral-term}-\ref{eq.surfaceintegral}), which has degree~2, reaching a total degree up to $2N+2$.	
	
	Finally, we remark that all the methods exactly preserve \textit{constant states}; this numerically verifies the satisfaction of the \textit{GCL condition} which, in fact, is ensured by construction by our methodology through integration over space-time control volumes.

	\begin{table}[tb]
		\centering
		\begin{minipage}{0.48\textwidth}		
			\centering
			\renewcommand{\arraystretch}{1.4} 
			\caption{3-2 flip. Conservation error for different density profiles $\rho(\mathbf{x}, t)$ and polynomial degrees $N$, \RIIcolor{including non polynomial profiles}.}
			\label{tab:conservation-errors-3-2}
			\resizebox{\textwidth}{!}{
				\begin{tabular}{lcccc}
					\toprule
					\parbox[c][1.2cm][c]{4cm}{$\displaystyle \left| \sum_{i=1}^{14} \left( \int_{P_i^1} \rho \, d\mathbf{x} - \int_{P_i^0} \rho \, d\mathbf{x} \right) \, \right|$} & $N=1$ & $N=2$ & $N=3$ & $\Rall{N=4}$ \\
					\midrule
					$\rho = 10$  & 8.52E-14   & 1.42E-14   & 5.68E-14 & \Rall{5.68E-14} \\ 
					$\rho = 10 + x + y + z$    & 4.47E-12   & 2.69E-11   & 3.64E-11 & \Rall{3.98E-13} \\ 
					$\rho = 10 + x^3 + y^2 + z^5 + xyz$   & 1.34E-12   & 2.36E-12   & 6.85E-12 & \Rall{3.27E-13} \\ 
					\RIIcolor{$\rho = 10 + e^{xy + y^3} + 1/(z+5)$}   & 7.11E-14   & 1.31E-12   & 5.97E-13 & \Rall{3.98E-13} \\  
					\bottomrule
				\end{tabular}
			}
		\end{minipage}
		\hfill
		\begin{minipage}{0.48\textwidth}			
			\centering
			\renewcommand{\arraystretch}{1.4} 
			\caption{3-2 flip. Numerical errors between numerical density $\rho$ at $t=t^1$ and exact density $\rho_E$ when expected to be of machine precision.}
			\label{tab:errors-3-2}
			\resizebox{\textwidth}{!}{
				\begin{tabular}{lcccc}
					\toprule
					\parbox[c][1.2cm][c]{4cm}{$\displaystyle \sqrt{\sum_{i=1}^{14} \left( \int_{P_i^1} \left( \rho -  \rho_E \right)^2 \, d\mathbf{x} \right)}$} & $N=1$ & $N=2$ & $N=3$ & \Rall{$N=4$} \\
					\midrule
					$\rho_E = 10 + x + y + z$ & 1.49E-12 & 1.61E-13 & 5.82E-13 & \Rall{3.99E-12} \\
					$\rho_E = 10 + x^2 + y^2 + z^2 + xz$ & - & 1.49E-13 & 8.25E-13 & \Rall{3.79E-12} \\
					$\rho_E = 10 + x^3 + y^3 + z^3 + xyz$ & - & - & 6.75E-13 & \Rall{4.49E-12} \\
					$\rho_E = 10 + x^4 + y^4 + z^4 + x y^2 z$ & - & - & - & \Rall{2.88E-12} \\
					\bottomrule
				\end{tabular}
			}
		\end{minipage}
	\end{table}
	\begin{table}[tb]
		\centering
		\begin{minipage}{0.48\textwidth}		
			\centering
			\renewcommand{\arraystretch}{1.4} 
			\caption{2-3 flip. Conservation error for different density profiles $\rho(\mathbf{x}, t)$ and polynomial degrees $N$, \RIIcolor{including non polynomial profiles}.}
			\label{tab:conservation-errors-2-3}
			\resizebox{\textwidth}{!}{
				\begin{tabular}{lcccc}
					\toprule
					\parbox[c][1.2cm][c]{4cm}{$\displaystyle \left| \sum_{i=1}^{14}  \left( \int_{P_i^1} \rho \, d\mathbf{x} - \int_{P_i^0} \rho \, d\mathbf{x} \right)\,\right|$} & $N=1$ & $N=2$ & $N=3$ & $N=4$ \\
					\midrule
					$\rho = 10$ & 9.95E-14   & 4.26E-14   & 1.42E-13 & \Rall{3.13E-13} \\ 
					$\rho = 10 + x + y + z$   & 5.68E-14   & 2.70E-11   & 1.84E-11 & \Rall{4.97E-13} \\ 
					$\rho = 10 + x^3 + y^2 + z^5 + xyz$ & 4.68E-12   & 0.00E-00   & 1.63E-12 & \Rall{2.84E-13}  \\ 
					\RIIcolor{$\rho = 10 + e^{xy + y^3} + 1/(z+5)$}   & 8.53E-14   & 6.96E-12   & 5.28E-11 & \Rall{5.68E-14} \\ 
					\bottomrule
				\end{tabular}
			}
		\end{minipage}
		\hfill
		\begin{minipage}{0.48\textwidth}			
			\centering
			\renewcommand{\arraystretch}{1.4} 
			\caption{2-3 flip. Numerical errors between numerical density $\rho$ at $t=t^1$ and exact density $\rho_E$ when expected to be of machine precision.}
			\label{tab:errors-2-3}
			\resizebox{\textwidth}{!}{
				\begin{tabular}{lcccc}
					\toprule
					\parbox[c][1.2cm][c]{4cm}{$\displaystyle \sqrt{ \sum_{i=1}^{14} \left( \int_{P_i^1} \left( \rho - \rho_E \right)^2\, d\mathbf{x} \right)}$} & $N=1$ & $N=2$ & $N=3$ & \Rall{$N=4$} \\
					\midrule
					$\rho_E = 10 + x + y + z$             & 2.95E-14 & 1.77E-13 & 4.87E-13 & \Rall{2.58E-12} \\
					$\rho_E = 10 + x^2 + y^2 + z^2 + xz$  & -        & 1.63E-13 & 5.32E-13 & \Rall{2.84E-12} \\
					$\rho_E = 10 + x^3 + y^3 + z^3 + xyz$ & -        & -        & 4.60E-13 & \Rall{3.29E-12} \\
					$\rho_E = 10 + x^4 + y^4 + z^4 + x y^2 z$  & -        & -        & -        & \Rall{4.30E-12} \\
					\bottomrule
				\end{tabular}
			}
		\end{minipage}
	\end{table}
	
	\paragraph{2-3 flip} The same elementary checks have been performed for the reverse operation, the 2-3 flip, to ensure the robustness of the general implementation. 
	In particular, we demonstrate that the 4D space-time \textit{volume}, exactly as in the previous case, is completely filled thanks to the presence of the \textit{hole-like} element, indeed
	\begin{equation*}
		\sum_{i=1}^{14} \int_{C_i^0} 1 \, d\mathbf{x} \, dt + \int_{H_{15}^0} 1 \, d\mathbf{x} \, dt = 
		\frac{\text{7.9974E-01}}{2N+1} +  \frac{\text{2.6042E-04}}{2N+1} = \frac{\text{0.8}}{2N+1}, \quad \text{for }\, N=1,2,3,\Rall{4}.		
	\end{equation*}
	Moreover, we report the numerical results concerning correct \textit{mass} conservation in Table~\ref{tab:conservation-errors-2-3}. We also show the difference between the numerical solution and the exact polynomial density profiles in Table~\ref{tab:errors-2-3}, which demonstrates the \textit{accuracy} of the method and, for constant states, the satisfaction of the \textit{GCL property}.

	\paragraph{4-4 flip} 
	To keep the presentation concise, regarding the 4-4 flip, we limit ourselves to show that even by using the \textit{hole-like} element of Figure~\ref{fig.4-4hole}, we are able to completely fill the 4D space-time \textit{volume} contained between $t^0$ and $t^1$, indeed	
	\begin{equation*}
		\sum_{i=1}^{14} \int_{C_i^0} 1 \, d\mathbf{x} \, dt + \int_{H_{15}^0} 1 \, d\mathbf{x} \, dt = 
		\frac{\text{7.9792E-01}}{2N+1} +  \frac{\text{2.0833E-03}}{2N+1} = \frac{\text{0.8}}{2N+1}   \quad \text{for }\, N=1,2,3,\Rall{4}.
	\end{equation*}

	\subsection{Sanity checks: rigid rotation of a spherical region within a stationary cubic enclosure}
	
	\begin{figure}[!t]
		\centering
		\begin{minipage}{1.0\textwidth} 
			\centering  
			\includegraphics[width=0.55\linewidth]{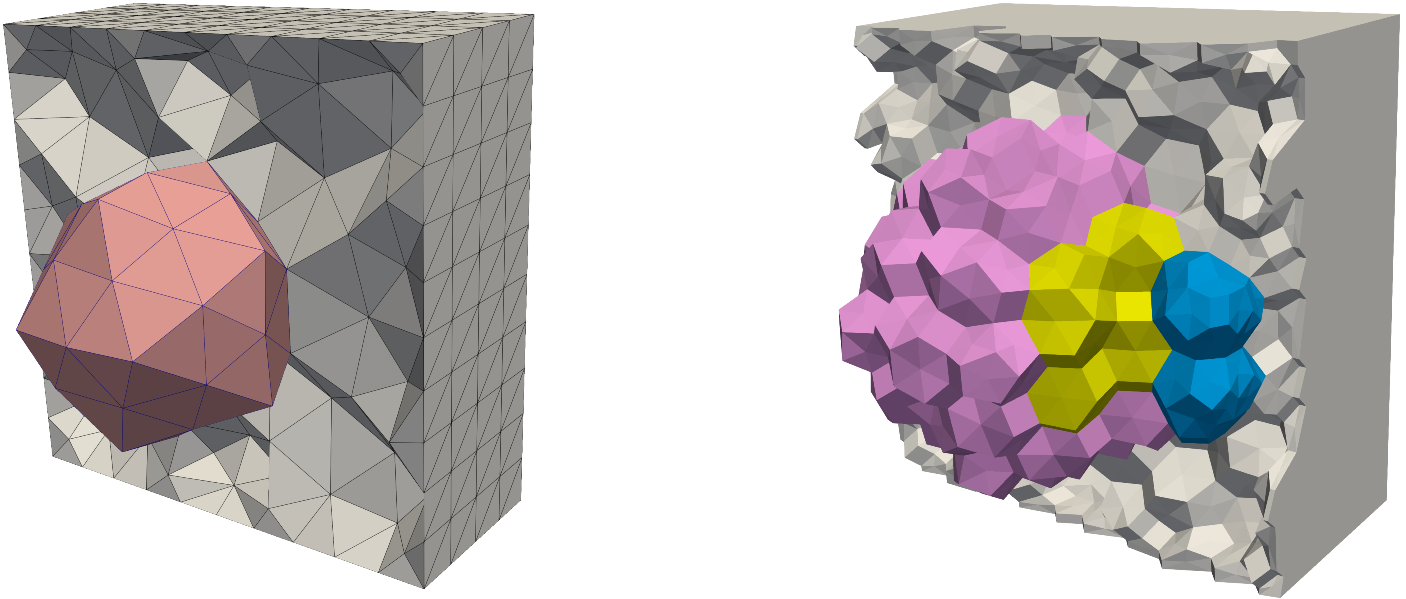}
			\caption{Rotating sphere test case. We show here the cutout of the initial tetrahedral mesh (left) with the elements having $r<0.3$ in red. On the right, we plot the corresponding polyhedral tessellation and we highlight a bunch of elements with colors depending on their numbering. We will then track throughout the simulation those colored elements: see Figure~\ref{fig.rot.genB}.}
			\label{fig.rot.gen}
		\end{minipage}		
		\vfill \vspace{1.5em} 
		\begin{minipage}{1.0\textwidth}
			\centering
			\includegraphics[width=1.0\linewidth]{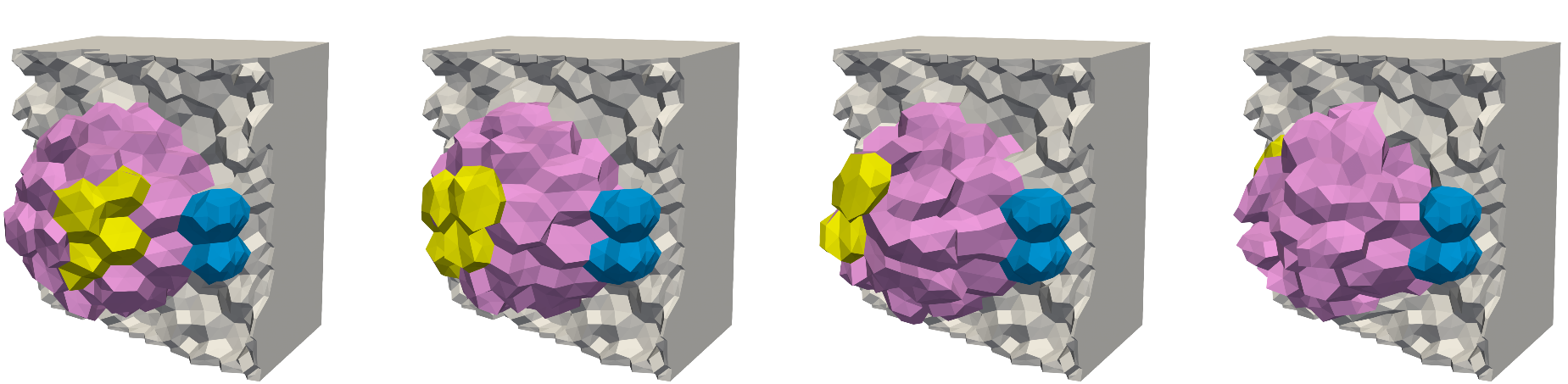}
			\caption{Rotating sphere test case. In this figure we highlight, maintaining the same color map of the previous Figure~\ref{fig.rot.gen}, the elements associated to chosen generator points. This shows the movement of the elements (yellow and violet) and the transformations that they are undergoing over time, specifically at times $t=0.25, 0.5, 0.75, 1.0$. The blue elements, outside the rotating region, are instead almost static with a few modifications due to the smoothing procedure.}
			\label{fig.rot.genB}
		\end{minipage}
	\end{figure}
	\begin{table}[!h] 
		\centering
		\caption{\RIIcolor{Rotating sphere test case solved up to $t=1$ on a mesh of $N_P=1021$ elements with our ALE ADER-DG scheme with $N=1,2,3$ generating a total of $\sum_n |\mathcal{H}^n|$ \textit{hole-like} elements during $n_f$ timesteps. 
				%
				%
				In the table below we report:
				i) the error between the numerical $\rho$ profile obtained at the end of the simulation and $\rho_0$, which is of machine accuracy for constant density and decreases when increasing the order of accuracy for the nonlinear profile;
				ii) the temporal mean and the temporal maximum of the mass conservation error over the timesteps, scaled by the total mass, which demonstrates that mass is conserved at every timestep, even for nonlinear, non-polynomial density profiles, independently of the mesh motion.
		}}
		\label{tab:ball-conservation-errors-1}
		\begin{tabular}{|c |c|c|c|c|c|}												
			\hline
			\rule{0pt}{6.5ex}\rule[-4.0ex]{0pt}{0pt} & $n_f$ & $\sum_n |\mathcal{H}^n|$ & $\displaystyle \sum_{i=1}^{N_P} \int_{P_i^{n_f}} \left| \rho - \rho_0 \right| $ & $\displaystyle \frac{ \frac{1}{n_f} \sum_{n=1}^{n_f} \left| \sum_{i=1}^{N_P} \left( \int_{P_i^n} \rho - \int_{P_i^{n-1}} \rho \right) \right|}{\sum_{i=1}^{N_P} \int_{P_i^0} \rho_0} \, $ & $\displaystyle \frac{ \max\limits_{n = 1,\dots,n_f} \left| \sum_{i=1}^{N_P} \left( \int_{P_i^n} \rho  - \int_{P_i^{n-1}} \rho \right) \right|}{\sum_{i=1}^{N_P} \int_{P_i^0} \rho_0 } \, $ \\[1.5ex] \hline \hline 
			$N$ & \multicolumn{5}{l|}{\rule{0pt}{2.1ex}\rule[-1.05ex]{0pt}{0pt}$\rho_0 = 1$} \\ \hline
			$1$ & 1551 & 3148 & 6.16\es{-15} & 7.80\es{-17} & 9.99\es{-16} \\ 
			$2$ & 2898 & 3156 & 4.18\es{-14} & 9.48\es{-17} & 9.99\es{-16} \\ 
			$3$ & 4922 & 3200 & 1.34\es{-13} & 1.51\es{-16} & 9.99\es{-16} \\ \hline \hline
			$N$ & \multicolumn{5}{l|}{\rule{0pt}{2.4ex}\rule[-1.05ex]{0pt}{0pt}\RIIcolor{$\rho_0 = 10 + e^{xy + y^3} + 1/(z+5)$}} \\ \hline
			$1$ & 1507 & 3161 & 1.13\es{-02} & 8.44\es{-13} & 1.27\es{-11} \\ 
			$2$ & 2591 & 3045 & 9.86\es{-04} & 5.91\es{-13} & 1.51\es{-11} \\ 
			$3$ & 4427 & 3082 & 7.35\es{-05} & 3.74\es{-13} & 9.14\es{-12} \\ \hline 
		\end{tabular}
	\end{table}	
	
	Now, we propose a more complex test case designed to verify that the above properties, in particular \RIIcolor{mass conservation and} GCL, are valid  even under mesh movements that trigger multiple successive topology changes. 
	
	We consider stationary initial condition\RIIcolor{s} \RIIcolor{with $\rho_0$ equal either to a constant $1.0$ or to a non polynomial profile  $10 + e^{xy + y^3} + 1/(z+5)$}, constant pressure and zero fluid velocity; however, in the spirit of the ALE methodology, we prescribe the motion of several generator points. 
	Specifically, we employ \Rall{$N_P = 1021$} internal generator points arranged in spherical layers. All generators within a radius $r < 0.3$ rotate around the origin in the x-y plane with a prescribed angular velocity $\omega = \pi$, while the remainder of the domain is associated with a zero mesh velocity field. In Figure~\ref{fig.rot.gen}, we illustrate the underlying tetrahedral mesh at time $t=0$ and the corresponding initial polyhedral tessellation. \Rall{The mesh smoothing parameter is set to $\kappa = 1/50$ and to gain greater control over the dihedral angles on the very coarse mesh the threshold of maximum dihedral angle in \eqref{eq_mesh_quality} is set to $\beta^{\min} = 90^{\circ}$ in this test.}
	
	We remark that most of the reconnections leading to the necessity of \textit{hole-like} elements take place in the spherical layer between the rotating and the stationary parts of the mesh. 
	In Figure~\ref{fig.rot.genB}, we highlight several polyhedral elements in that portion of the domain with different colors, and we show their position and shape at different simulation times. 
	We observe that the central elements are indeed rotating and undergoing numerous topology changes, while the external ones remain essentially static, except for minor variations due to the smoothing. 	
	All elements maintain excellent quality even with a mesh velocity field that is discontinuous.
	
	More importantly, this benchmark allows us to verify that even in the presence of many \textit{hole-like} elements the method is conservative for each $N=1,2,3$ and maintains exactly the constant density profile, as shown in Table~\ref{tab:ball-conservation-errors-1}, which is possible only when the GCL is effectively satisfied up to machine precision.

	\subsection{Shu-type \RIcolor{moving} isentropic vortex}
	\label{s-sec.Shu}
	
	\begin{table}[tb]
		\centering \small
		\caption{ 
			Convergence results for the Shu-type \RIcolor{moving} vortex solved with our ALE ADER-DG method
			of order $2$ ($N=1$). 
			Here, we compare the results obtained with a classical ALE scheme, 
			where the mesh is simply stretched and relocated without allowing for topology changes (left), 
			with the results obtained using our approach, which allows for topology changes and treats 
			them via the insertion of \textit{hole-like} elements (right). 
			In particular, we report the total number of timesteps $n_f$ and the cumulative 
			count of \textit{hole-like} elements $\sum_n |\mathcal{H}^n|$ generated during 
			the simulation; furthermore, we show the size of the central mesh cell at the 
			final time $h_c^{n_f}$, the $L_2$ error norms for $\rho$ and the corresponding 
			order of accuracy at times $t=0.2$, $t=1.0$, and $t=10.0$.
			We remark that these results were obtained after hundreds of \textit{hole-like} elements had been generated and a complete revolution had been performed. This demonstrates that the order of convergence is fully maintained even in the presence of large mesh deformations. 
			We also refer to Section~\ref{s-sec.Shu} for an additional and detailed discussion of the obtained results.
		}
		%
		%
		\begin{tabular}{|l |c |c | c c | c ||c  | c | c c |c|}
			\hline
			\rule{0pt}{3ex} \multirow{2}{2cm}{\centering ALE ADER-DG $N=1$} & & \multicolumn{4}{c||}{ALE \textit{without} topology changes} & \multicolumn{5}{c|}{ALE \textit{with} topology changes} \\[1.5ex] \cline{3-11}
			\rule{0pt}{3ex} & $N_P$ & $n_f$ & $h_c^{n_f}$ & $L_2(\rho)$ & $\mathcal{O}(\rho)$ & $\sum_{n} |\mathcal{H}^n|$ & $n_f$ & $h_c^{n_f}$ & $L_2(\rho)$ & $\mathcal{O}(\rho)$ \\[1ex] \hline \hline
			\rule{0pt}{3ex} \multirow{4}{*}{$t=0.2$} 
			& 2293  & 67 & 9.02\es{-01} & 2.34\es{-01} & -- &  2   & 55  & 9.26\es{-01} & 2.36\es{-01} & -- \\
			& 5475  & 91 & 5.80\es{-01} & 1.08\es{-01} & 1.74 & 134  & 77  & 6.04\es{-01} & 1.08\es{-01} & 1.83 \\
			& 9381 & 128 & 4.81\es{-01} & 7.67\es{-02} & 1.84 & 119 & 101 & 5.00\es{-01} & 7.57\es{-02} & 1.87 \\
			& 19199 & 141 & 3.60\es{-01} & 4.31\es{-02} & 1.99 & 311  & 123 & 3.70\es{-01} & 4.11\es{-02} & 2.02 \\ \hline \hline
			\rule{0pt}{3ex} \multirow{4}{*}{$t=1.0$} 
			& 2293  & 331 & 8.98\es{-01} & 2.95\es{-01} & -- & 114  & 225 & 9.45\es{-01} & 2.96\es{-01} & -- \\
			& 5475  & 454 & 5.82\es{-01} & 1.46\es{-01} & 1.64 & 481  & 330 & 6.19\es{-01} & 1.36\es{-01} & 1.84 \\
			& 9381 & 638 & 4.82\es{-01} & 1.02\es{-01} & 1.89 &  797 & 427 & 5.12\es{-01} & 9.38\es{-02} & 1.96 \\
			& 19199 & 705 & 3.61\es{-01} & 5.50\es{-02} & 2.13 & 2076 & 563 & 3.79\es{-01} & 4.87\es{-02} & 2.15 \\ \hline \hline
			\rule{0pt}{3ex} \multirow{4}{*}{$t=10.0$} 
			& 2293  & \multicolumn{4}{c||}{{Simulation crashes at $t\sim1.60$}}  & 1116  & 2011  & 8.92\es{-01} & 1.02\es{-00} & -- \\ 
			& 5475  &  \multicolumn{4}{c||}{{Simulation crashes at $t\sim1.59$}}  &  9704 & 3220  & 5.64\es{-01} & 4.83\es{-01} & 1.63 \\
			& 9381 &  \multicolumn{4}{c||}{{Simulation crashes at $t\sim1.53$}}  &  18704 & 3990 & 4.66\es{-01} & 3.13\es{-01} & 2.27 \\
			& 19199 &  \multicolumn{4}{c||}{{Simulation crashes at $t\sim1.72$}}  &  57288 &  5447 &  3.34\es{-01} & 1.65\es{-01} & 2.07 \\ \hline
		\end{tabular}
		\label{tab.shu.orderN1}
	\end{table}
	
	\begin{table}[htb]	
		\centering \small
		\caption{ 
			Convergence results for the Shu-type \RIcolor{moving} vortex solved with our ALE ADER-DG method
			of order $3$ ($N=2$). 
			We refer to the caption of Table~\ref{tab.shu.orderN1} and to Section~\ref{s-sec.Shu} for a more detailed discussion of the results presented here.
		}
		%
		%
		\begin{tabular}{|l |c |c | c c | c ||c  | c | c c |c|}
			\hline
			\rule{0pt}{3ex} \multirow{2}{2cm}{\centering ALE ADER-DG $N=2$} & & \multicolumn{4}{c||}{ALE \textit{without} topology changes} & \multicolumn{5}{c|}{ALE \textit{with} topology changes} \\[1.5ex] \cline{3-11}
			\rule{0pt}{3ex} & $N_P$ & $n_f$ & $h_c^{n_f}$ & $L_2(\rho)$ & $\mathcal{O}(\rho)$ & $\sum_{n} |\mathcal{H}^n|$ & $n_f$ & $h_c^{n_f}$ & $L_2(\rho)$ & $\mathcal{O}(\rho)$ \\[1ex] \hline \hline
			\rule{0pt}{3ex} \multirow{4}{*}{$t=0.2$} 
			& 2293  & 129 & 9.04\es{-01} & 5.77\es{-02} & -- &  6   & 103 & 9.32\es{-01} & 5.76\es{-02} & -- \\
			& 5475  & 177 & 5.81\es{-01} & 1.86\es{-02} & 2.56 & 147  & 144 & 6.08\es{-01} & 1.84\es{-02} & 2.67 \\
			& 9381  & 248 & 4.81\es{-01} & 1.14\es{-02} & 2.59 & 128 & 188 & 5.02\es{-01} & 1.12\es{-02} & 2.62 \\
			& 19199 & 274 & 3.60\es{-01} & 4.93\es{-03} & 2.90 & 345  & 233 & 3.72\es{-01} & 4.61\es{-03} & 2.94 \\ \hline \hline
			\rule{0pt}{3ex} \multirow{4}{*}{$t=1.0$} 
			& 2293  & 642 & 9.12\es{-01} & 7.68\es{-02} & -- & 114  & 422 & 9.49\es{-01} & 7.05\es{-02} & -- \\
			& 5475  & 881 & 5.82\es{-01} & 2.56\es{-02} & 2.44 & 425  & 626 & 6.25\es{-01} & 2.15\es{-02} & 2.84 \\
			& 9381  & 1238 & 4.82\es{-01} & 1.61\es{-02} & 2.48 &  721 & 811 & 5.14\es{-01} & 1.24\es{-02} & 2.82 \\
			& 19199 & 1368 & 3.60\es{-01} & 7.19\es{-03} & 2.77 & 1795 & 1082 & 3.80\es{-01} & 5.23\es{-03} & 2.84 \\ \hline \hline
			\rule{0pt}{3ex} \multirow{4}{*}{$t=10.0$} 
			& 2293  &  \multicolumn{4}{c||}{{Simulation crashes at $t\sim1.37$}}   &   711 & 3900  & 8.57\es{-01} & 1.34\es{-01} & -- \\ 
			& 5475  &  \multicolumn{4}{c||}{{Simulation crashes at $t\sim1.58$}}   &  7521 & 6069  & 5.86\es{-01} & 4.59\es{-02} & 2.81 \\
			& 9381  &  \multicolumn{4}{c||}{{Simulation crashes at $t\sim1.51$}}   &  16648 & 7728  & 4.42\es{-01} & 2.30\es{-02} & 2.46 \\
			& 19199 &  \multicolumn{4}{c||}{{Simulation crashes at $t\sim1.74$}}   &  50150 & 10565 &  3.23\es{-01} & 7.69\es{-03} & 3.52 \\ \hline
		\end{tabular}
		\label{tab.shu.orderN2}
	\end{table}
	\begin{table}[htb]	
		\centering \small
		\caption{ 
			Convergence results for the Shu-type \RIcolor{moving} vortex solved with our ALE ADER-DG mehod
			of order $4$ ($N=3$).  
			We refer to the caption of Table~\ref{tab.shu.orderN1} and to Section~\ref{s-sec.Shu} for a more detailed discussion of the results presented here.
		}
		%
		%
		\begin{tabular}{|l |c |c | c c | c ||c  | c | c c |c|}
			\hline
			\rule{0pt}{3ex} \multirow{2}{1.5cm}{\centering ALE ADER-DG $N=3$} & & \multicolumn{4}{c||}{ALE \textit{without} topology changes} & \multicolumn{5}{c|}{ALE \textit{with} topology changes} \\[1.5ex] \cline{3-11}
			\rule{0pt}{3ex} & $N_P$ & $n_f$ & $h_c^{n_f}$ & $L_2(\rho)$ & $\mathcal{O}(\rho)$ & $\sum_{n} |\mathcal{H}^n|$ & $n_f$ & $h_c^{n_f}$ & $L_2(\rho)$ & $\mathcal{O}(\rho)$ \\[1ex] \hline \hline
			\rule{0pt}{3ex} \multirow{4}{*}{$t=0.2$} 
			& 2293  & 219 & 9.03\es{-01} & 1.30\es{-02} & -- &  18  & 168 & 9.36\es{-01} & 1.33\es{-02} & -- \\
			& 5475  & 300 & 5.80\es{-01} & 2.78\es{-03} & 3.49 & 153  & 238 & 6.11\es{-01} & 2.76\es{-03} & 3.68 \\
			& 9381  & 421 & 4.81\es{-01} & 1.41\es{-03} & 3.59 & 131  & 310 & 5.04\es{-01} & 1.38\es{-03} & 3.60 \\
			& 19199 & 465 & 3.60\es{-01} & 4.63\es{-04} & 3.85 & 732  & 382 & 3.80\es{-01} & 4.70\es{-04} & 3.80 \\ \hline \hline
			\rule{0pt}{3ex} \multirow{4}{*}{$t=1.0$} 
			& 2293  & 1092 & 9.05\es{-01} & 1.81\es{-02} & -- &  92  & 699 & 9.45\es{-01} & 1.63\es{-02} & -- \\
			& 5475  & 1498 & 5.80\es{-01} & 4.22\es{-03} & 3.27 & 356  & 1050 & 6.28\es{-01} & 3.22\es{-03} & 3.97 \\
			& 9381  & 2104 & 4.81\es{-01} & 2.28\es{-03} & 3.27 &  594 & 1362 & 5.16\es{-01} & 1.62\es{-03} & 3.50 \\
			& 19199 & 2325 & 3.60\es{-01} & 7.84\es{-04} & 3.68 & 1900 & 1823 & 3.86\es{-01} & 5.42\es{-04} & 3.77 \\ \hline \hline
			\rule{0pt}{3ex} \multirow{4}{*}{$t=10.0$} 
			& 2293 & \multicolumn{4}{c||}{{Simulation crashes at $t\sim1.29$}} &   308 & 6463  & 8.97\es{-01} & 4.18\es{-02} & -- \\ 
			& 5475 & \multicolumn{4}{c||}{{Simulation crashes at $t\sim1.54$}} &  3174 & 10083 & 5.45\es{-01} & 5.26\es{-03} & 4.16 \\
			& 9381 & \multicolumn{4}{c||}{{Simulation crashes at $t\sim1.50$}} &  9217 & 13144 & 4.79\es{-01} & 2.51\es{-03} &  5.72\\
			& 19199 & \multicolumn{4}{c||}{{Simulation crashes at $t\sim1.72$}} &  41564 & 17957 &  3.56\es{-01} & 7.74\es{-04} & 3.97 \\ \hline
		\end{tabular}
		\label{tab.shu.orderN3}
	\end{table}

	\begin{figure}[!htb]
		\centering
		\includegraphics[width=1.0\linewidth]{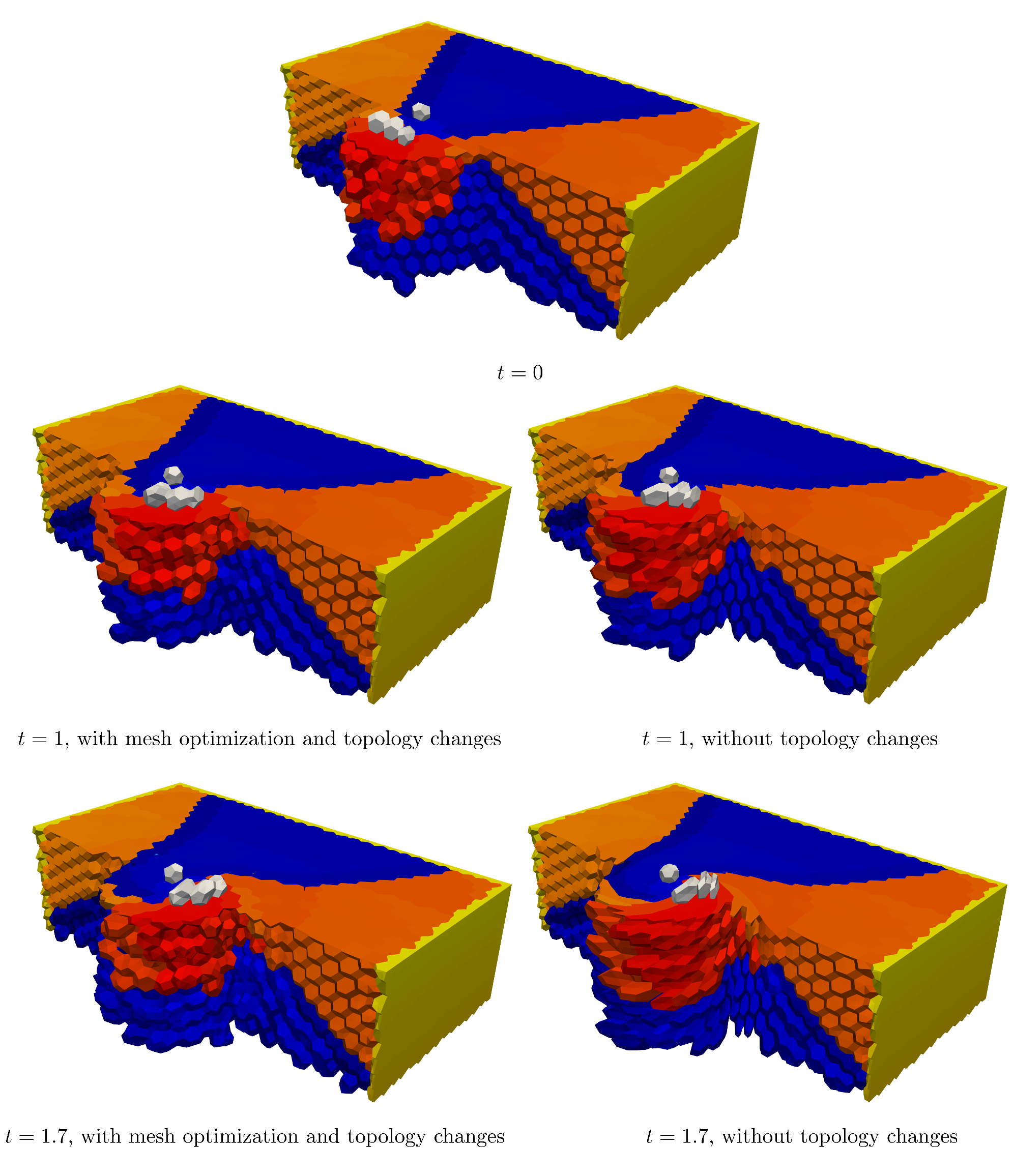}
		\caption{\Rall{ 
				In this figure, we report a three-dimensional view of the Shu vortex simulation ($\mathbf{u}_T = [0,0,0]$) obtained with our ALE ADER-DG method with $N=2$ with (left) and without (right) topology changes at $t=1$, i.e., well before the simulation would crash due to excessive mesh distortion, and at $t=1.7$, i.e., just before the crash. 
				The coloring depicts the element numbering on a vertical mesh cutout shown to illustrate the internal movement of the elements during the rotation. Several standalone elements (central and three rotating) are highlighted above the mesh cutout. The element shape is visibly elongated and distorted when preventing topology changes.}
		}
		\label{fig:shu_p2_mesh}
	\end{figure}
	
	\begin{figure}[!htb]
		\centering
		\includegraphics[width=1.0\linewidth]{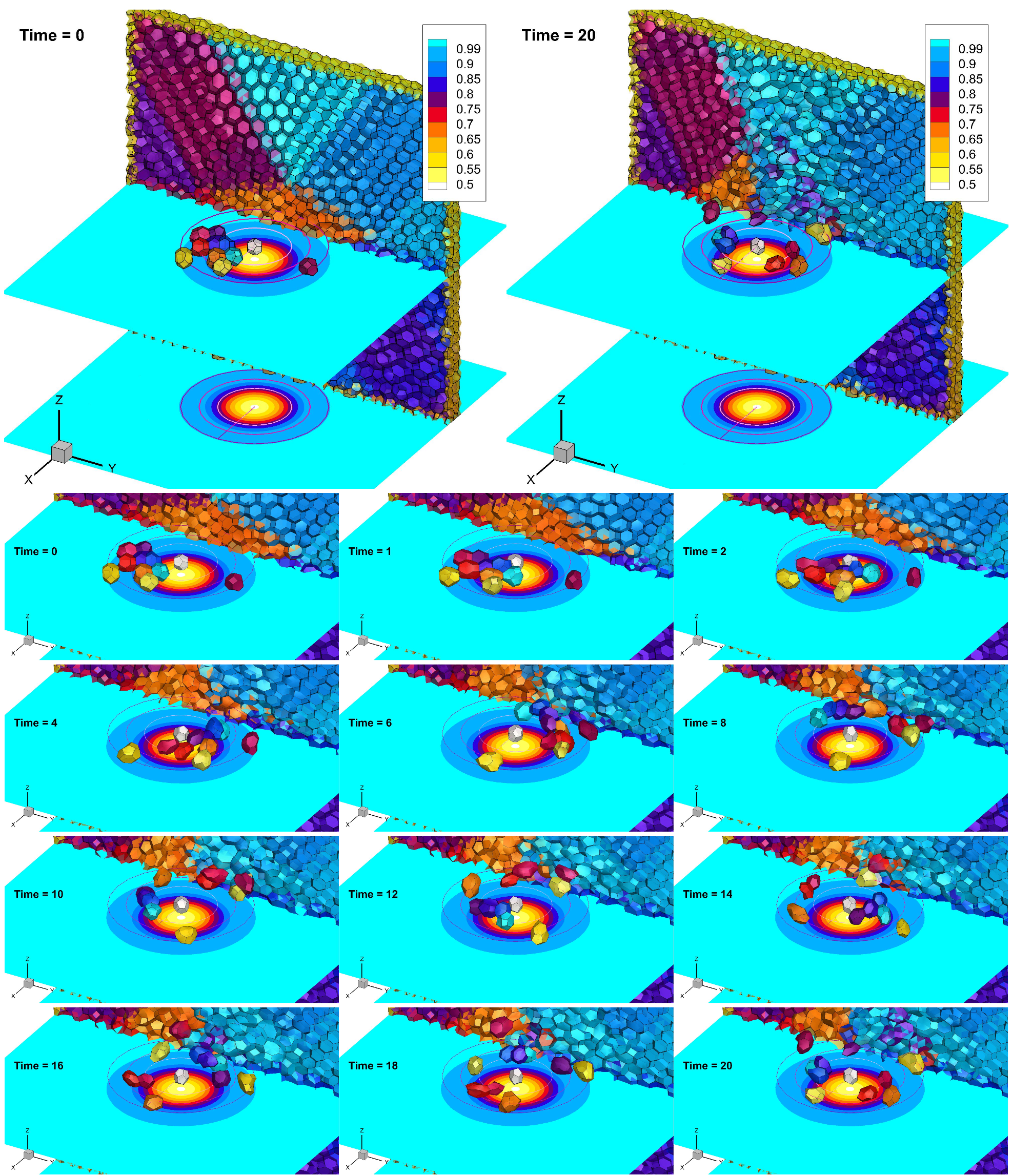}
		\caption{In this figure, we report a three-dimensional view of the Shu vortex simulation ($\mathbf{u}_T = [0,0,0]$)  obtained with our ALE ADER-DG method with $N=3$ at several time instances. 
			The horizontal cuts ($z=0, 4.6$) highlight the perfectly preserved vortical profile of the density $\rho$, which remains visually identical to its initial state throughout the entire simulation. Simultaneously, the filled vertical section  ($2.5 < x < 3.5$) and the bunch of tracked colored elements at the center demonstrate the effective rotational movement of the mesh, showing the continuous changes in position, shape, and topology of the polyhedral cells, together with their optimal quality.}
		\label{fig:goodshusolp3}
	\end{figure}

	Our main benchmark to validate the capabilities of our scheme in the presence of large mesh deformations due to vortical dynamics is represented by a classical smooth isentropic vortex, such as the one proposed in~\cite{HuShuVortex1999}, which constitutes a stationary solution of the Euler equations \RIcolor{and undergoes a translating movement with constant velocity $\mathbf{u}_T = [1,1,0]$}.
	
	As computational domain, we take the cube $\Omega=[0,10]^3$ with wall boundary conditions 
	\RIcolor{on top and bottom boundaries (normal to $z$-coordinate) and transmissive boundary conditions 
		on the remaining part of the boundary.} 
	We cover the domain with an increasingly refined set of polyhedral tessellations. At any given time, the characteristic mesh size is chosen as the size~$h_c^n$, computed as in~\eqref{eq.char.length} for the central polyhedron of the domain at that time.	
	The initial conditions are defined as:
	\begin{equation*}
		\left\{
		\begin{aligned}
			& \rho^0 (x,y,z) = (1 + \delta T)^{\frac{1}{\gamma-1}} \\
			& u^0 (x,y,z)    = \RIcolor{\mathbf{u}_{T,x}}-(y-5\, \RIcolor{- \,\mathbf{u}_{T,y}} \, t) \frac{\epsilon}{2\pi}e^{\frac{1-r^2}{2}} \\
			& v^0 (x,y,z)    = \RIcolor{\mathbf{u}_{T,y}}+(x-5\, \RIcolor{- \,\mathbf{u}_{T,x}} \, t) \frac{\epsilon}{2\pi}e^{\frac{1-r^2}{2}} \\
			& w^0 (x,y,z)    = 0 \\
			& p^0 (x,y,z)    = (1 + \delta T)^{\frac{\gamma}{\gamma-1}}
		\end{aligned}
		\right.
		\qquad \text{with } \qquad
		\begin{aligned}
			&\delta T = -\frac{(\gamma-1)\epsilon^2}{8\gamma\pi^2}e^{1-r^2}, \\
			&r        = \sqrt{(x-5\ \RIcolor{-\, \mathbf{u}_{T,x}} \, t)^2 + (y-5\ \RIcolor{-\, \mathbf{u}_{l,y}} \, t)^2}, \\
			&\epsilon = 5, 
		\end{aligned}
	\end{equation*}
	which consists of a translationally symmetric extension of the 2D vortex, where the flow field remains invariant along the $z$-coordinate. Moreover, we note that the vortex perturbations and thus its rotational velocity vanish as $x, y \rightarrow 0$~or~$10$.
	Next, we move the computational domain according to the strategy described in Section~\ref{s-sec.mesh_movement} for $z \in  [2, 8]$, thus using here a velocity as close as possible to the fluid velocity \RIcolor{with the smoothing parameter $\kappa = 1/100$}. Outside this region, we progressively and smoothly damp the fluid velocity so as to obtain zero mesh velocity as $z \rightarrow 0$ or $10$. In this manner, we can study a vortical problem involving a challenging rotational velocity inside the domain without the need to move the external boundaries \RIcolor{in a different manner than pure translation}.
	
	We report in Tables~\ref{tab.shu.orderN1}, \ref{tab.shu.orderN2} and~\ref{tab.shu.orderN3}, respectively for order of accuracy $N=1,2$ and $3$, several results which allow us to prove the capabilities of our scheme as follows. 
	
	\begin{itemize} 
		\item First, we remark that without a full optimization including topology changes, the mesh would be completely destroyed even before a quarter of a rotation, forcing the algorithm to break down by $t \approx 2$ at the latest; conversely, with the techniques introduced in this paper, we can complete full rotations without any mesh distortion issues.
		
		\item Moreover, if we compare the numerical errors obtained on a mesh that cannot be fully optimized (even well before it becomes completely tangled) with the errors produced by our new method, we observe that ours are always remarkably lower, especially as $N$ increases and the mesh is refined. 
		
		\item Additionally, focusing on the right-hand side of the tables, one can see the number of \textit{hole-like} elements effectively required during the simulation of the rotational phenomenon. Their number remains consistently low per timestep, which underlines that only a few locations need optimization at each step; this demonstrates that, although we allow only one flip per generator per timestep, our approach remains feasible as the necessary optimization operations do not accumulate. On the other hand, their presence is frequent enough to justify a dedicated strategy for their management.
		
		\item We can also observe that the number of timesteps required to simulate a given time interval remains fairly constant throughout the entire simulation: this indicates that the mesh is well-optimized \RIIcolor{even though our strategy allows at most one \textit{hole-like} element to be created per generator per timestep}.
		
		\item As a final observation, the order of accuracy matches the theoretical $N+1$ for each $N$, even in the presence of multiple \textit{hole-like} elements, mesh deformation, and very long simulation times. Regarding this last point, we wish to specify that the values of our tables are particularly noteworthy because, on unstructured meshes, and especially with a mesh size that even varies in time, the $N+1$ trend does not typically appear so precisely. Indeed, in the literature, it is common to report only an average order even for stationary cases; therefore, we consider the results presented here to be quite remarkable.
	\end{itemize}
	
	Moreover, we report in Figure~\ref{fig:goodshusolp3} a three-dimensional view of a Shu vortex simulation, 
	\Rall{executed with with $\mathbf{u}_T = [0,0,0]$ for the sake of a clearer visualization}, 
	at several time instances within \Rall{$t \in [0, 20]$}, obtained with our ALE ADER-DG method for $N=3$.
	In particular, the horizontal cuts along the planes at $z=0$ and $z=4.6$ display the vortical density solution, demonstrating that its quality remains perfectly preserved throughout the simulation; indeed, the solution remains visually identical to its initial state, showing no signs of degradation. 
	This notable result is enabled by the combination of our Lagrangian method, which allows for mesh optimization, and the high order of accuracy that we achieve even in the presence of topology changes.
	We also present two vertical cuts along the planes $x=3.6$ and $y=3.6$, where the colors represent the mesh numbering, to show that all cells effectively move during the simulation. To further illustrate the movement of these polyhedral cells, specifically their change in position, connectivity, and shape, we highlight a bunch of elements at the center of the domain and track them throughout the simulation, keeping the color associated with their numbering fixed.
	
	\Rall{Furthermore, to show the difference in mesh quality when allowing or not allowing topology changes while employing smoothing in both cases, we refer to Figure~\ref{fig:shu_p2_mesh}, where the mesh configuration is presented right before the crash as well as a certain time before it occurs.}

\begin{figure}[!t]
		\centering   \vspace{-30pt}
		
		\begin{minipage}[c]{0.48\textwidth}
			\centering
			\small
			\begin{tabular}{|l |c |c | c |}
				\hline
				\rule{0pt}{3ex} $t_{\textrm{WALL}} [s]$ & $N=1$ & $N=2$ & $N=3$ \\[1ex] \hline \hline
				\rule{0pt}{3ex} Mesh smoothing, optimization & 10.9 & 9.03 & 7.34 \\ \hline
				\rule{0pt}{3ex} Space-time mesh construction & 0.839 & 0.886 & 0.93 \\ \hline
				\rule{0pt}{3ex} Predictor, std. elements & 3.12 & 33.5 & 652 \\ \hline
				\rule{0pt}{3ex} Predictor, \textit{hole-like} elements & 0.0708 & 0.0827 & 0.444 \\ \hline
				\rule{0pt}{3ex} Corrector & 4.22 & 28.8 & 145 \\ \hline
			\end{tabular}
			\vspace{0.2em}
			
			{\small Table \the\numexpr\value{table}+1\relax} 
		\end{minipage}
		\hfill
		\begin{minipage}[c]{0.48\textwidth}
			\centering
			\includegraphics[width=\linewidth, trim=5cm 10cm 5cm 10cm, clip]{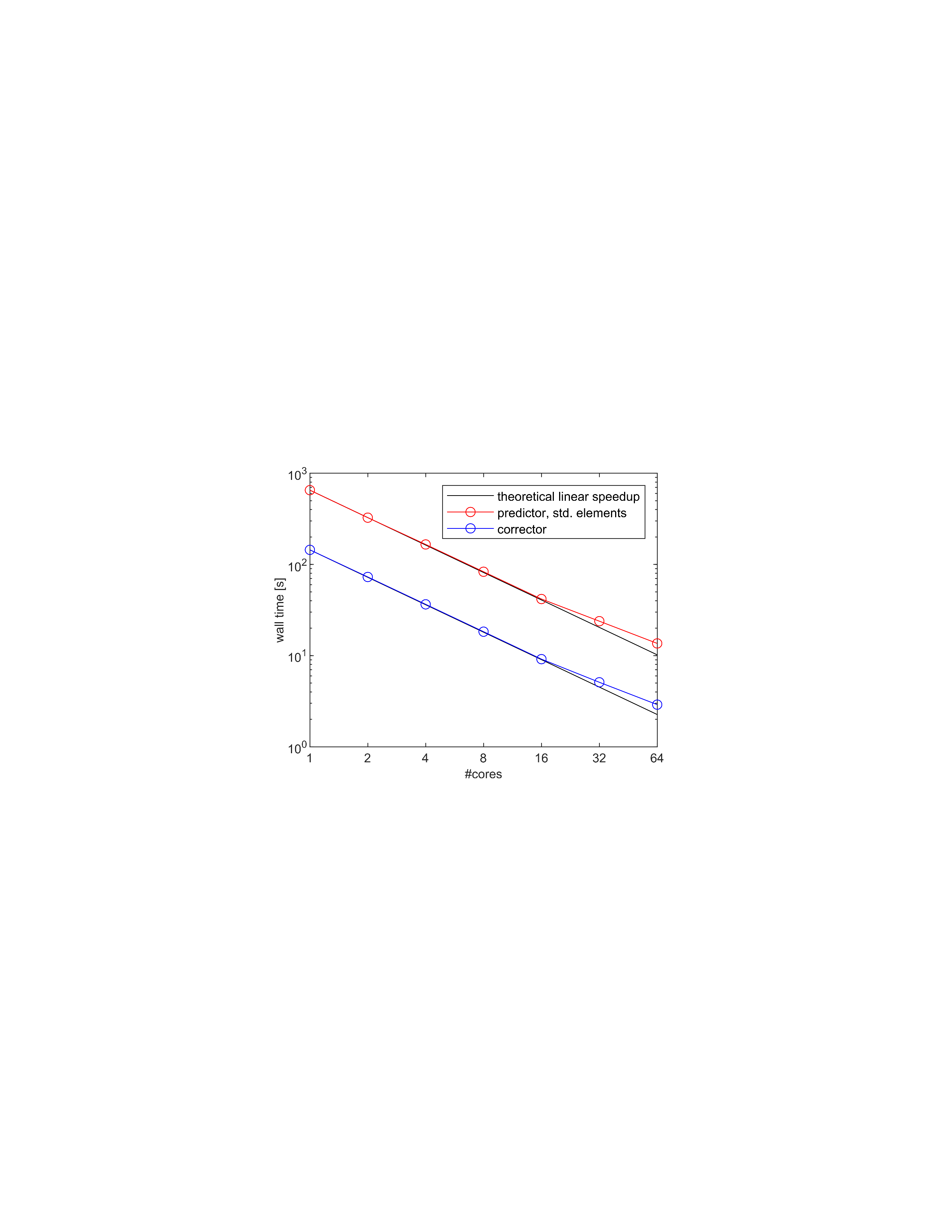}
			
			{\small Figure \the\numexpr\value{figure}+1\relax} 
		\end{minipage}
		
		\captionof{table}{%
			\RIIcolor{Approximate CPU computational time of each part of our ALE ADER-DG method, respectively with $N=1, 2$, and $3$ when applied to the Shu-type vortex simulation on a fine mesh of $N_P=19199$ elements.
				The time displayed here has been measured as the averaged cost of one timestep over 100 consecutive timesteps starting from $t=9$ when the movement of the mesh already plays its role. 
				The timing measurement was realized on a single core although in real computations OpenMP is employed for the parallel execution of the predictor and corrector steps over standard elements, as shown in the next Figure~\ref{fig:shu_scaling}
				We note that the cost related to \textit{hole-like} elements is negligible with respect to the entire algorithm.}
		}
		\label{tab.shu.cost}
		
		\captionof{figure}
		{		
			\RIIcolor{CPU time scaling for the parallel execution of our ALE ADER-DG method with $N=3$ on standard space-time control volumes for both the predictor and corrector steps, measured in the same conditions of the previous Table~\ref{tab.shu.cost}. 
				The code is parallelized via the OpenMP paradigm, and the test has been run on a shared memory system with two AMD EPYC 9654 96-core processors.
			}
		}
		\label{fig:shu_scaling}
		
	\end{figure}

	\RIIcolor{Finally, in Table~\ref{tab.shu.cost} we show the approximate computational cost related to the different parts of our whole numerical algorithm for $N=1, 2$, and $3$. 
		The shown results refer to the time necessary to perform one timestep on a single core and have been computed as the average time over 100 timesteps executed on a fine mesh with $N_P=19199$ starting from $t=9$ in such a way to consider a situation where the movement already plays an important role. 
		
		The computationally most intensive parts of the algorithm are, as expected, those connected with the solver (predictor and corrector steps) over standard elements and not those referred to the mesh preparation. 
		Thus, to guarantee the feasibility of our simulation, the solver part has also been parallelized with the OpenMP paradigm and, as depicted in Figure~\ref{fig:shu_scaling}, it shows very good scaling with respect to the number of cores typically available for shared memory multi-thread execution. The parallelization of the rest of the algorithm is possible and will be implemented soon.
		
		We just want to remark that the construction of the space-time geometry of \textit{hole-like} elements and the locally implicit strategy performed over them have a decidedly negligible cost compared to the entire algorithm.}

	\RIcolor{
		\subsection{Brief discussion on differences, advantages, and disadvantages of our direct ALE algorithm versus classical versions of the indirect ALE approach}

		While our current validation has focused on checking the effects of our direct ALE approach with and without topology changes, leaving a comprehensive comparison between direct and indirect ALE approaches to future work, we provide here a preliminary qualitative comparison between the two methodologies, acknowledging that this discussion is intended to be concise and does not aim for an exhaustive analysis.
		
		In classical indirect ALE algorithms, the process begins with a purely Lagrangian phase during which the mesh is moved and the solution evolved accordingly. This approach offers the \textit{distinct advantage} of \textit{fully exploiting} the strengths of Lagrangian methods, providing \textit{zero} errors at contact discontinuities and material interfaces, alongside the absence of mass transfer between cells, which is critical for multi-material simulations. However, it may happen that additional limitations affect the selection of the timestep size in order to guarantee that the mesh does not become tangled during the Lagrangian phase.
		Moreover, even if several purely Lagrangian steps may be performed, the resulting mesh \textit{must} eventually be optimized to enable subsequent motion either once it begins to show the first signs of distortion or when it becomes tangled (depending on the specific method employed). During this optimization stage, even in indirect ALE algorithms, the formation of multi-material cells is often unavoidable.
		
		Conversely, our direct ALE algorithm does not benefit from the advantages of a \textit{strictly pure} Lagrangian evolution from the beginning, because we consistently apply a controlled amount of smoothing to the mesh motion to ensure its optimal quality.
		Thus, we cannot decouple the process into a distinct purely Lagrangian phase;
		nevertheless, this strategy prevents the accumulation of geometric errors associated with mesh distortion and ensures that the timestep selection remains unconstrained by tangling issues while guaranteeing a large part of the advantages of the Lagrangian motion.
		
		So far, we believe the choice between the two methodologies depends on the specific requirements of the application: indirect ALE is generally preferred when the exact preservation of interfaces in multi-material simulations is paramount, whereas direct ALE is more suitable for applications prioritizing the reduction of convective errors and the accurate capture of contact discontinuities.
		Furthermore, we acknowledge that indirect methods, having been developed for many more years than ours, are often equipped with more advanced mesh smoothing algorithms. 
		
		However, we believe the primary strength of our algorithm lies in the fact that it does not require a \textit{remap phase}, which instead characterizes indirect methods. In our approach, the space-time evolution through integration over space-time control volumes, completes the timestep directly.
		
		In contrast, in indirect algorithms, every time the mesh is optimized, it is necessary to remap the previous Lagrangian solution onto the new mesh. This phase is particularly delicate because the \textit{remap} must ensure
		conservativity, oscillation-free solutions, compliance with the GCL, and high order of accuracy, which are instead automatically derived by construction in our approach. A remap with such properties is in general highly challenging because it requires exact mesh intersection and high order accurate integration within every sub-portion, especially when the solution is represented by high order discontinuous polynomials rather than simple cell averages. 
		In particular, mesh intersection in 3D is a computationally expensive operation that can also be prone to inaccuracies, as two segments connecting floating-point coordinates might intersect at a point not representable within floating-point precision. Additionally, quadrature in variable geometric regions is costly, requiring a large number of non-precomputable quadrature points.
		Therefore, we believe our approach is particularly interesting and worthy of future investigation precisely because it is \textit{remap-free}.
	}

	
	\section{Conclusions} \label{sec_concl}	
	
	Although this work does not yet showcase complex hyperbolic simulations, 
	we believe it introduces a groundbreaking perspective on handling topology changes within the context of \textit{direct} ALE schemes. 
	The approach presented here follows a path already explored in 2D by authors such as 
	Springel, Alauzet, Gaburro \textit{et al.} in~\cite{springel2010pur,olivier2011new,gaburro2020high,gaburro2025high} and recently provided with theoretical foundations in~\cite{bonafini2026stability}. 
	However, to the best of our knowledge, this is the first time in the literature that a complete and well-justified proof of concept is demonstrated in a full 3D setting.
	Thus, we feel that a dedicated description of the main building blocks, specifically 
	the construction and integration over \textit{hole-like} space-time elements and their neighborhood, 
	is essential to clarify the strategy we have developed. 
	
	A key aspect that sets this work apart from the current state of the art is the simultaneous management 
	of topology changes and the integration of non-linear PDEs around them maintaining conservativity up to machine precision and the desired order of convergence. 
	While many existing methods achieve impressive results in terms of mesh quality and optimization, 
	they often do not couple these techniques with the actual evolution of the physical solution. 
	Furthermore, when they do so, they typically follow an indirect approach that requires projecting and reconstructing 
	the solution between different meshes. 
	This is a well-explored path, that however must face the significant computational difficulties related to mesh intersection \cite{levy2025exact}. 
	Our novel approach, instead, embeds the mesh connectivity changes directly into the space-time integration framework, 
	avoiding the need for projection-reconstruction procedures between the meshes. 
	
	Moreover, via introducing the \textit{implicit} treatment of the \textit{hole-like} elements alongside the \textit{explicit} update of classical elements, we have proposed an entirely new technology that mixes implicit and explicit schemes in a geometry-based rather than PDE-based way. 
	This approach may find significant application also in other fields where it is necessary to work with degenerate or near-degenerate geometries, 
	such as in cut-cell methods~\cite{may2022dod,may2024accuracy,birke2026error}.	
	
	Currently, we are committed to improving our mesh optimization techniques, which remains a time-consuming task because the specific information required by our space-time algorithm prevents us from using standard optimizers as \textit{black boxes}. 
	Instead, we must have full control over the optimization routines to ensure they are compatible with our ALE formulation. Nevertheless, this is a feasible path, and simulations obtained through more 
	advanced optimization techniques will be the subject of future publications.
	\Rall{Furthermore, future work will include a comprehensive comparison between direct and indirect ALE methods. As a longer-term objective, we intend to decouple our solver to enable a purely Lagrangian phase for the solution evolution, while retaining the benefits of integration over \textit{hole-like} elements to achieve an intersection-free, automatically conservative, and high order accurate remap.}

	This work also serves as the formal introduction of the main building blocks of 
	our novel code, \textit{TheALcHEMiST} (\textit{The} reads as \textit{D} and stays for Direct), which stands for \textit{Direct Arbitrary Lagrangian Eulerian High order Schemes on Moving Space Time polyhedral Meshes}.


	\section*{Acknowledgments}
	
	E.~Gaburro and M.~Tavelli are members of the INdAM GNCS group in Italy;
	M.~Bonafini is member of the INdAM GNAMPA group in Italy.	
	All the authors gratefully acknowledge the support received from the European Union 
	with the ERC Starting Grant \textit{ALcHyMiA} (grant agreement No. 101114995).
	Views and opinions expressed are however those of the author only and do not necessarily 
	reflect those of the European Union or the European Research Council Executive Agency. 
	Neither the European Union nor the granting authority can be held responsible for them.

	\section*{Declaration of generative AI and AI-assisted technologies in the manuscript preparation process}
	
	The authors declare that AI-assisted technologies, specifically ChatGPT and Gemini, were utilized exclusively as interactive dictionaries to improve the linguistic quality of the manuscript. These tools were solely employed to check English grammar, correct typographical errors, enhance text fluidity, and identify appropriate synonyms. 
	Crucially, the AI tools were only prompted to refine phrases and sentences entirely conceived, structured, and drafted by the authors; no automated text generation was requested or used. Every correction proposed by the AI was meticulously reviewed, critically evaluated, and manually accepted only upon strict verification by the authors. 
	
	No other part of this work, including, but not limited to, the core ideas, scientific concepts, methodologies, bibliography, code, or figures, was generated or influenced by Artificial Intelligence. The responsibility for the content and integrity of the final manuscript remains entirely with the authors.

	\bibliographystyle{abbrv}
	\bibliography{referencesTheALcHEMiST_Inception_black}

\begin{thebibliography}{100}

\bibitem{alauzet2014changing}
F.~Alauzet.
\newblock A changing-topology moving mesh technique for large displacements.
\newblock {\em Engineering with Computers}, 30(2):175--200, 2014.

\bibitem{anderson2018high}
R.~Anderson, V.~Dobrev, T.~Kolev, R.~Rieben, and V.~Tomov.
\newblock {High-order multi-material ALE hydrodynamics}.
\newblock {\em SIAM Journal on Scientific Computing}, 40(1):B32--B58, 2018.

\bibitem{balsara2012self}
D.~Balsara.
\newblock Self-adjusting, positivity preserving high order schemes for
  hydrodynamics and magnetohydrodynamics.
\newblock {\em Journal of Computational Physics}, 231(22):7504--7517, 2012.

\bibitem{BARLOW2018ALE}
A.~Barlow, M.~Klima, and M.~Shashkov.
\newblock {Constrained optimization framework for interface-aware sub-scale
  dynamics models for voids closure in Lagrangian hydrodynamics}.
\newblock {\em Journal of Computational Physics}, 371:914--944, 2018.

\bibitem{indALE-AWE2016}
A.~Barlow, P.~H. Maire, , W.~Rider, R.~Rieben, and M.~Shashkov.
\newblock {Arbitrary Lagrangian–Eulerian methods for modeling high-speed
  compressible multimaterial flows}.
\newblock {\em Journal of Computational Physics}, 322:603--665, 2016.

\bibitem{belytschko1978computer}
T.~Belytschko and J.~Kennedy.
\newblock {Computer models for subassembly simulation}.
\newblock {\em Nuclear Engineering and Design}, 49(1-2):17--38, 1978.

\bibitem{benson1992computational}
D.~Benson.
\newblock {Computational methods in Lagrangian and Eulerian hydrocodes}.
\newblock {\em Computer methods in Applied mechanics and Engineering},
  99(2-3):235--394, 1992.

\bibitem{MaireMM2}
M.~Berndt, J.~Breil, S.~Galera, M.~Kucharik, P.~H. Maire, and M.~Shashkov.
\newblock {Two--step hybrid conservative remapping for multimaterial arbitrary
  Lagrangian--Eulerian methods}.
\newblock {\em Journal of Computational Physics}, 230:6664--6687, 2011.

\bibitem{birke2026error}
G.~Birke, C.~Engwer, J.~Giesselmann, and S.~May.
\newblock {Error Analysis of a First-Order DoD Cut Cell Method for 2D Unsteady
  Advection: G. Birke et al.}
\newblock {\em Journal of Scientific Computing}, 106(1):1, 2026.

\bibitem{ReALE2015}
W.~Bo and M.~Shashkov.
\newblock {Adaptive reconnection-based arbitrary Lagrangian Eulerian method}.
\newblock {\em Journal of Computational Physics}, 299:902--939, 2015.

\bibitem{ShashkovRemap1}
P.~Bochev, D.~Ridzal, and M.~Shashkov.
\newblock Fast optimization-based conservative remap of scalar fields through
  aggregate mass transfer.
\newblock {\em Journal of Computational Physics}, 246:37--57, 2013.

\bibitem{bonafini2026stability}
M.~Bonafini, D.~Torlo, and E.~Gaburro.
\newblock {Stability analysis of Arbitrary-Lagrangian-Eulerian ADER-DG methods
  on classical and degenerate spacetime geometries}.
\newblock {\em arXiv preprint arXiv:2602.09198}, 2026.

\bibitem{Lagrange2D}
W.~Boscheri and M.~Dumbser.
\newblock {Arbitrary--Lagrangian--Eulerian One--Step WENO Finite Volume Schemes
  on Unstructured Triangular Meshes}.
\newblock {\em Communications in Computational Physics}, 14:1174--1206, 2013.

\bibitem{Lagrange3D}
W.~Boscheri and M.~Dumbser.
\newblock A direct {Arbitrary-Lagrangian-Eulerian ADER-WENO} finite volume
  scheme on unstructured tetrahedral meshes for conservative and
  non-conservative hyperbolic systems in 3d.
\newblock {\em Journal of Computational Physics}, 275:484 -- 523, 2014.

\bibitem{boscheriAFE2022}
W.~Boscheri, M.~Dumbser, and E.~Gaburro.
\newblock {Continuous Finite Element Subgrid Basis Functions for Discontinuous
  Galerkin Schemes on Unstructured Polygonal Voronoi Meshes}.
\newblock {\em Communications in Computational Physics}, 32(1):259--298, 2022.

\bibitem{boscheri2013semi}
W.~Boscheri, M.~Dumbser, and M.~Righetti.
\newblock {A semi-implicit scheme for 3D free surface flows with high-order
  velocity reconstruction on unstructured Voronoi meshes}.
\newblock {\em International journal for numerical methods in fluids},
  72(6):607--631, 2013.

\bibitem{ALELTS2D}
W.~Boscheri, M.~Dumbser, and O.~Zanotti.
\newblock {High order cell-centered Lagrangian-type finite volume schemes with
  time-accurate local time stepping on unstructured triangular meshes}.
\newblock {\em Journal of Computational Physics}, 291:120--150, 2015.

\bibitem{breil20253d}
J.~Breil, G.~Damour, S.~Guisset, and A.~Cola{\"\i}tis.
\newblock 3d mesh regularization within an ale code using a weighted line
  sweeping method.
\newblock {\em Computers \& Fluids}, 292:106591, 2025.

\bibitem{busto2020high}
S.~Busto, S.~Chiocchetti, M.~Dumbser, E.~Gaburro, and I.~Peshkov.
\newblock {High order ADER schemes for continuum mechanics}.
\newblock {\em Frontiers in Physics}, 8:32, 2020.

\bibitem{busto2022new}
S.~Busto and M.~Dumbser.
\newblock A new family of thermodynamically compatible discontinuous {Galerkin}
  methods for continuum mechanics and turbulent shallow water flows.
\newblock {\em Journal of Scientific Computing}, 93(2):56, 2022.

\bibitem{CaramanaShashkov1998}
E.~Caramana and M.~Shashkov.
\newblock Elimination of artificial grid distorsion and hourglass–type
  motions by means of {Lagrangian} subzonal masses and pressures.
\newblock {\em Journal of Computational Physics}, 142:521--561, 1998.

\bibitem{Castro2008}
M.~Castro, J.~Gallardo, J.~L\'opez, and C.~Par\'es.
\newblock Well-balanced high order extensions of godunov's method for
  semilinear balance laws.
\newblock {\em SIAM Journal of Numerical Analysis}, 46:1012--1039, 2008.

\bibitem{Castro2006}
M.~Castro, J.~Gallardo, and C.~Par\'es.
\newblock High-order finite volume schemes based on reconstruction of states
  for solving hyperbolic systems with nonconservative products. {A}pplications
  to shallow-water systems.
\newblock {\em Mathematics of Computation}, 75:1103--1134, 2006.

\bibitem{chiocchetti2021high}
S.~Chiocchetti, I.~Peshkov, S.~Gavrilyuk, and M.~Dumbser.
\newblock {High order ADER schemes and GLM curl cleaning for a first order
  hyperbolic formulation of compressible flow with surface tension}.
\newblock {\em Journal of Computational Physics}, 426:109898, 2021.

\bibitem{ciallella2024very}
M.~Ciallella, S.~Clain, E.~Gaburro, and M.~Ricchiuto.
\newblock {Very high order treatment of embedded curved boundaries in
  compressible flows: ADER discontinuous Galerkin with a space-time
  Reconstruction for Off-site data}.
\newblock {\em Computers \& Mathematics with Applications}, 175:1--18, 2024.

\bibitem{colaitis2026cell-1}
A.~Cola{\"\i}tis, S.~Guisset, and J.~Breil.
\newblock {A cell-centered Amr-Ale framework for 3d multi-material
  hydrodynamics. Part I: Lagrangian and indirect Euler Amr algorithms}.
\newblock {\em Journal of Computational Physics}, page 114701, 2026.

\bibitem{colaitis2026cell-2}
A.~Cola{\"\i}tis, S.~Guisset, and J.~Breil.
\newblock {A cell-centered AMR-ALE framework for 3D multi-material
  hydrodynamics. Part II: linesweep ALE rezoning for nonconformal
  block-structured AMR meshes}.
\newblock {\em Journal of Computational Physics}, page 114702, 2026.

\bibitem{cremonesi2010lagrangian}
M.~Cremonesi, A.~Frangi, and U.~Perego.
\newblock A lagrangian finite element approach for the analysis of
  fluid--structure interaction problems.
\newblock {\em International Journal for Numerical Methods in Engineering},
  84(5):610--630, 2010.

\bibitem{cremonesi2017explicit}
M.~Cremonesi, S.~Meduri, U.~Perego, and A.~Frangi.
\newblock An explicit lagrangian finite element method for free-surface weakly
  compressible flows.
\newblock {\em Computational Particle Mechanics}, 4(3):357--369, 2017.

\bibitem{DASSI20182}
F.~Dassi, L.~Kamenski, P.~Farrell, and S.~Hang.
\newblock {Tetrahedral mesh improvement using moving mesh smoothing, lazy
  searching flips, and RBF surface reconstruction}.
\newblock {\em Computer-Aided Design}, 103:2--13, 2018.
\newblock 25th International Meshing Roundtable Special Issue: Advances in Mesh
  Generation.

\bibitem{delisleGeorge}
E.~B. de~L'isle and P.~L. George.
\newblock Optimization of tetrahedral meshes.
\newblock In {\em Modeling, Mesh Generation, and Adaptive Numerical Methods for
  Partial Differential Equations. The IMA Volumes in Mathematics and its
  Applications.}, volume~75, pages 97--127. Springer New York, 1995.

\bibitem{del2023triangular}
S.~Del~Pino and I.~Marmajou.
\newblock {Triangular metric-based mesh adaptation for compressible
  multi-material flows in semi-Lagrangian coordinates}.
\newblock {\em Journal of Computational Physics}, 478:111975, 2023.

\bibitem{despres2017numerical}
B.~Despr{\'e}s.
\newblock {\em {Numerical methods for Eulerian and Lagrangian conservation
  laws}}.
\newblock Birkh{\"a}user, 2017.

\bibitem{despres2024lagrangian}
B.~Despr{\'e}s.
\newblock Lagrangian vorono{\"\i} meshes and particle dynamics with shocks.
\newblock {\em Computer Methods in Applied Mechanics and Engineering},
  418:116427, 2024.

\bibitem{di20243d}
D.~Di~Cristofaro, A.~Frangi, and M.~Cremonesi.
\newblock {3d fluid--structure interaction simulation with an
  Arbitrary--Lagrangian--Eulerian approach with applications to flying
  objects}.
\newblock {\em Engineering with Computers}, pages 1--23, 2024.

\bibitem{Dobrev3}
V.~Dobrev, T.~Ellis, T.~Kolev, and R.~Rieben.
\newblock {High order curvilinear finite elements for axisymmetric Lagrangian
  hydrodynamics}.
\newblock {\em Computers and Fluids}, 83:58--69, 2013.

\bibitem{donea1977lagrangian}
J.~Don{\'e}a, P.~Fasoli-Stella, and S.~Giuliani.
\newblock {Lagrangian and Eulerian finite element techniques for transient
  fluid-structure interaction problems}.
\newblock In {\em Structural mechanics in reactor technology}. 1977.

\bibitem{donea1982arbitrary}
J.~Donea, S.~Giuliani, and J.-P. Halleux.
\newblock {An arbitrary Lagrangian-Eulerian finite element method for transient
  dynamic fluid-structure interactions}.
\newblock {\em Computer methods in applied mechanics and engineering},
  33(1-3):689--723, 1982.

\bibitem{ALELTS1D}
M.~Dumbser.
\newblock {Arbitrary--Lagrangian--Eulerian ADER--WENO finite volume schemes
  with time--accurate local time stepping for hyperbolic conservation laws}.
\newblock {\em Computer Methods in Applied Mechanics and Engineering},
  280:57--83, 2014.

\bibitem{dumbser2008unified}
M.~Dumbser, D.~Balsara, E.~Toro, and C.-D. Munz.
\newblock {A unified framework for the construction of one-step finite volume
  and discontinuous Galerkin schemes on unstructured meshes}.
\newblock {\em Journal of Computational Physics}, 227(18):8209--8253, 2008.

\bibitem{dumbser2019glm}
M.~Dumbser, F.~Fambri, E.~Gaburro, and A.~Reinarz.
\newblock {On GLM curl cleaning for a first order reduction of the CCZ4
  formulation of the Einstein field equations}.
\newblock {\em Journal of Computational Physics}, page 109088, 2019.

\bibitem{dumbser2016high}
M.~Dumbser, I.~Peshkov, E.~Romenski, and O.~Zanotti.
\newblock {High order {ADER} schemes for a unified first order hyperbolic
  formulation of continuum mechanics: viscous heat-conducting fluids and
  elastic solids}.
\newblock {\em Journal of Computational Physics}, 314:824--862, 2016.

\bibitem{dumbser2017high}
M.~Dumbser, I.~Peshkov, E.~Romenski, and O.~Zanotti.
\newblock {High order {ADER} schemes for a unified first order hyperbolic
  formulation of {Newtonian }continuum mechanics coupled with
  electro-dynamics}.
\newblock {\em Journal of Computational Physics}, 348:298--342, 2017.

\bibitem{dumbser2023WBGR}
M.~Dumbser, O.~Zanotti, E.~Gaburro, and I.~Peshkov.
\newblock {A well-balanced discontinuous Galerkin method for the first--order
  Z4 formulation of the Einstein--Euler system}.
\newblock {\em Journal of Computational Physics}, page 112875, 2024.

\bibitem{eder20263d}
D.~Eder, C.~Parisuana, A.~Fisher, P.~Yip, T.~Schwartzentruber, and A.~Koniges.
\newblock {A 3D Arbitrary Lagrangian Eulerian Adaptive Mesh Refinement
  Framework for Simulating Hypersonic Raindrop Impacts}.
\newblock In {\em AIAA SCITECH 2026 Forum}, page 0730, 2026.

\bibitem{fambri2020discontinuous}
F.~Fambri.
\newblock {Discontinuous Galerkin methods for compressible and incompressible
  flows on space--time adaptive meshes: toward a novel family of efficient
  numerical methods for fluid dynamics}.
\newblock {\em Archives of Computational Methods in Engineering},
  27(1):199--283, 2020.

\bibitem{fernandez2022arbitrary}
E.~G. Fern{\'a}ndez, M.~C. D{\'\i}az, M.~Dumbser, and T.~M. De~Luna.
\newblock An arbitrary high order well-balanced {ADER-DG} numerical scheme for
  the multilayer shallow-water model with variable density.
\newblock {\em Journal of Scientific Computing}, 90(1):52, 2022.

\bibitem{freitagGooch}
L.~A. Freitag and C.~Ollivier-Gooch.
\newblock Tetrahedral mesh improvement using swapping and smoothing.
\newblock {\em International Journal for Numerical Methods in Engineering},
  40(21):3979--4002, 1997.

\bibitem{gaburro2021unified}
E.~Gaburro.
\newblock {A unified framework for the solution of hyperbolic PDE systems using
  high order direct Arbitrary-Lagrangian--Eulerian schemes on moving
  unstructured meshes with topology change}.
\newblock {\em Archives of Computational Methods in Engineering},
  28(3):1249--1321, 2021.

\bibitem{gaburro2025high}
E.~Gaburro.
\newblock {High order well-balanced Arbitrary-Lagrangian-Eulerian ADER
  discontinuous Galerkin schemes on general polygonal moving meshes}.
\newblock {\em Computers \& Fluids}, page 106764, 2025.

\bibitem{gaburro2020high}
E.~Gaburro, W.~Boscheri, S.~Chiocchetti, C.~Klingenberg, V.~Springel, and
  M.~Dumbser.
\newblock High order direct arbitrary-lagrangian-eulerian schemes on moving
  voronoi meshes with topology changes.
\newblock {\em Journal of Computational Physics}, 407:109167, 2020.

\bibitem{gaburro2024discontinuous}
E.~Gaburro, W.~Boscheri, S.~Chiocchetti, and M.~Ricchiuto.
\newblock {Discontinuous Galerkin schemes for hyperbolic systems in
  non-conservative variables: quasi-conservative formulation with subcell
  finite volume corrections}.
\newblock {\em Computer Methods in Applied Mechanics and Engineering},
  431:117311, 2024.

\bibitem{gaburro2021bookchapter}
E.~Gaburro and S.~Chiocchetti.
\newblock {High-order Arbitrary-Lagrangian-Eulerian schemes on crazy moving
  Voronoi meshes}.
\newblock In {\em Young Researchers Conference}, pages 99--119. Springer, 2021.

\bibitem{gaburro2017direct}
E.~Gaburro, M.~Dumbser, and M.~Castro.
\newblock {Direct Arbitrary-Lagrangian-Eulerian finite volume schemes on moving
  nonconforming unstructured meshes}.
\newblock {\em Computers \& Fluids}, 159:254--275, 2017.

\bibitem{gaburro2026innovations}
E.~Gaburro, M.~Kazolea, S.~Chiocchetti, M.~Dumbser, and R.~Abgrall.
\newblock Innovations in modeling and numerical methods for evolutionary pdes:
  Theory and applications, 2026.

\bibitem{Gaburro2026}
E.~Gaburro, M.~Ricchiuto, and M.~Dumbser.
\newblock On general and complete multidimensional riemann solvers for
  nonlinear systems of hyperbolic conservation laws.
\newblock {\em Computers and Fluids}, 311, 2026.

\bibitem{gallice2022entropy}
G.~Gallice, A.~Chan, R.~Loub{\`e}re, and P.-H. Maire.
\newblock {Entropy stable and positivity preserving Godunov-type schemes for
  multidimensional hyperbolic systems on unstructured grid}.
\newblock {\em Journal of Computational Physics}, 468:111493, 2022.

\bibitem{guisset2024cell}
S.~Guisset, G.~Damour, and J.~Breil.
\newblock {Cell-centered indirect Arbitrary Lagrangian-Eulerian numerical
  strategy for solving 3D gas dynamics equations}.
\newblock {\em Journal of Computational Physics}, 505:112903, 2024.

\bibitem{han2021dec}
M.~Han~Veiga, P.~{\"O}ffner, and D.~Torlo.
\newblock Dec and ader: similarities, differences and a unified framework.
\newblock {\em Journal of Scientific Computing}, 87(1):1--35, 2021.

\bibitem{hang2015tetgen}
S.~Hang.
\newblock {TetGen, a Delaunay-based quality tetrahedral mesh generator}.
\newblock {\em ACM Trans. Math. Softw}, 41(2):11, 2015.

\bibitem{hidalgo2011ader}
A.~Hidalgo and M.~Dumbser.
\newblock {ADER schemes for nonlinear systems of stiff
  advection--diffusion--reaction equations}.
\newblock {\em Journal of Scientific Computing}, 48(1-3):173--189, 2011.

\bibitem{hirt1974arbitrary}
C.~Hirt, A.~Amsden, and J.~Cook.
\newblock {An arbitrary Lagrangian-Eulerian computing method for all flow
  speeds}.
\newblock {\em Journal of computational physics}, 14(3):227--253, 1974.

\bibitem{HuShuVortex1999}
C.~Hu and C.~Shu.
\newblock {A high-order WENO finite difference scheme for the equations of
  ideal magnetohydrodynamics.}
\newblock {\em Journal of Computational Physics}, 150:561 -- 594, 1999.

\bibitem{hughes1981lagrangian}
T.~J. Hughes, W.~K. Liu, and T.~Zimmermann.
\newblock {Lagrangian-Eulerian finite element formulation for incompressible
  viscous flows}.
\newblock {\em Computer methods in applied mechanics and engineering},
  29(3):329--349, 1981.

\bibitem{jackson2017eigenvalues}
H.~Jackson.
\newblock {On the eigenvalues of the ADER-WENO Galerkin predictor}.
\newblock {\em Journal of Computational Physics}, 333:409--413, 2017.

\bibitem{kemm2020simple}
F.~Kemm, E.~Gaburro, F.~Thein, and M.~Dumbser.
\newblock {A simple diffuse interface approach for compressible flows around
  moving solids of arbitrary shape based on a reduced Baer--Nunziato model}.
\newblock {\em Computers \& fluids}, 204:104536, 2020.

\bibitem{kenamond2021positivity}
M.~Kenamond, D.~Kuzmin, and M.~Shashkov.
\newblock A positivity-preserving and conservative
  intersection-distribution-based remapping algorithm for staggered ale
  hydrodynamics on arbitrary meshes.
\newblock {\em Journal of Computational Physics}, 435:110254, 2021.

\bibitem{kincl2025numerical}
O.~Kincl, I.~Peshkov, and W.~Boscheri.
\newblock {A numerical method based on quasi-Lagrangian Voronoi cells for
  two-phase flows with large density contrast}.
\newblock {\em Computers \& Fluids}, page 106813, 2025.

\bibitem{KLIMA2020ALE}
M.~Klima, A.~Barlow, M.~Kucharik, and M.~Shashkov.
\newblock {An interface-aware sub-scale dynamics multi-material cell model for
  solids with void closure and opening at all speeds}.
\newblock {\em Computers \& Fluids}, 208:104578, 2020.

\bibitem{KLIMA2024LAG}
M.~Klima, M.~Kucharik, and R.~Liska.
\newblock {Cell-centered Lagrangian Lax-Wendroff HLL hybrid scheme in 3D}.
\newblock {\em Journal of Computational Physics}, 498:112675, 2024.

\bibitem{klingerSchewchuk}
B.~M. Klingner and J.~R. Shewchuk.
\newblock Aggressive tetrahedral mesh improvement.
\newblock In {\em Proceedings of the 16th International Meshing Roundtable
  (Seattle, Washington)}, pages 3--23, 2007.

\bibitem{KUCHARIK2008CONNECTIVITY}
M.~Kucharik and M.~Shashkov.
\newblock {Extension of efficient, swept-integration-based conservative
  remapping method for meshes with changing connectivity}.
\newblock {\em International Journal for Numerical Methods in Fluids},
  56(8):1359--1365, 2008.

\bibitem{lakiss2024ader}
A.~Lakiss, T.~Heuz{\'e}, M.~Tannous, and L.~Stainier.
\newblock {ADER discontinuous Galerkin material point method}.
\newblock {\em International Journal for Numerical Methods in Engineering},
  125(1):e7365, 2024.

\bibitem{lei2023high}
N.~Lei, J.~Cheng, and C.-W. Shu.
\newblock A high order positivity-preserving polynomial projection remapping
  method.
\newblock {\em Journal of Computational Physics}, 474:111826, 2023.

\bibitem{levy2025exact}
B.~L{\'e}vy.
\newblock {Exact predicates, exact constructions and combinatorics for mesh
  CSG.}
\newblock {\em ACM Transactions on Graphics}, 44(5):1--27, 2025.

\bibitem{chengshu2}
W.~Liu, J.~Cheng, and C.~Shu.
\newblock {High order conservative Lagrangian schemes with Lax--Wendroff type
  time discretization for the compressible Euler equations}.
\newblock {\em Journal of Computational Physics}, 228:8872--8891, 2009.

\bibitem{ReALE2011}
R.~Loub{\`e}re, P.~H. Maire, and M.~Shashkov.
\newblock {ReALE: A Reconnection Arbitrary-Lagrangian–Eulerian method in
  cylindrical geometry}.
\newblock {\em Computers and Fluids}, 46:59--69, 2011.

\bibitem{ReALE2010}
R.~Loub{\`e}re, P.~H. Maire, M.~Shashkov, J.~Breil, and S.~Galera.
\newblock {ReALE: A reconnection-based arbitrary-Lagrangian–Eulerian method}.
\newblock {\em Journal of Computational Physics}, 229:4724--4761, 2010.

\bibitem{LoubereShashkov2004}
R.~Loub\`ere and M.~Shashkov.
\newblock {A subcell remapping method on staggered polygonal grids for
  arbitrary--lagrangian--eulerian methods}.
\newblock {\em Journal of Computational Physics}, 23:155--160, 2004.

\bibitem{Maire2010}
P.~H. Maire.
\newblock A unified sub-cell force-based discretization for cell-centered
  {Lagrangian} hydrodynamics on polygonal grids.
\newblock {\em International Journal for Numerical Methods in Fluids},
  65:1281--1294, 2011.

\bibitem{may2024accuracy}
S.~May and F.~Laakmann.
\newblock Accuracy analysis for explicit-implicit finite volume schemes on cut
  cell meshes.
\newblock {\em Communications on Applied Mathematics and Computation},
  6(4):2239--2264, 2024.

\bibitem{may2022dod}
S.~May and F.~Streitb{\"u}rger.
\newblock {DoD stabilization for non-linear hyperbolic conservation laws on cut
  cell meshes in one dimension}.
\newblock {\em Applied Mathematics and Computation}, 419:126854, 2022.

\bibitem{micalizzi2023efficient}
L.~Micalizzi, D.~Torlo, and W.~Boscheri.
\newblock Efficient iterative arbitrary high-order methods: an adaptive bridge
  between low and high order.
\newblock {\em Communications on Applied Mathematics and Computation}, pages
  1--38, 2023.

\bibitem{mill}
R.~Millington, E.~Toro, and L.~Nejad.
\newblock {\em Arbitrary High Order Methods for Conservation Laws I: The One
  Dimensional Scalar Case}.
\newblock PhD thesis, Manchester Metropolitan University, Department of
  Computing and Mathematics, June 1999.

\bibitem{morgan2021origins}
N.~Morgan and B.~Archer.
\newblock {On the origins of Lagrangian hydrodynamic methods}.
\newblock {\em Nuclear Technology}, 207(sup1):S147--S175, 2021.

\bibitem{munz94}
C.~Munz.
\newblock {On Godunov--type schemes for Lagrangian gas dynamics}.
\newblock {\em SIAM Journal on Numerical Analysis}, 31:17--42, 1994.

\bibitem{muzzolon2026high}
S.~Muzzolon, M.~Dumbser, O.~Zanotti, and E.~Gaburro.
\newblock High order numerical discretizations of the einstein-euler equations
  in the generalized harmonic formulation.
\newblock {\em Journal of Computational Physics}, 563:115084, 2026.

\bibitem{olivier2011new}
G.~Olivier and F.~Alauzet.
\newblock A new changing-topology ale scheme for moving mesh unsteady
  simulations.
\newblock In {\em 49th AIAA Aerospace Sciences Meeting including the New
  Horizons Forum and Aerospace Exposition}, page 474, 2011.

\bibitem{Pares2006}
C.~Par\'es.
\newblock Numerical methods for nonconservative hyperbolic systems: a
  theoretical framework.
\newblock {\em SIAM Journal on Numerical Analysis}, 44:300--321, 2006.

\bibitem{plessier2023implicit}
A.~Plessier, S.~Del~Pino, and B.~Despr{\'e}s.
\newblock Implicit discretization of lagrangian gas dynamics.
\newblock {\em ESAIM: Mathematical Modelling and Numerical Analysis},
  57(2):717--743, 2023.

\bibitem{popov2023space}
I.~Popov.
\newblock {Space-Time Adaptive ADER-DG Finite Element Method with LST-DG
  Predictor and a posteriori Sub-cell WENO Finite-Volume Limiting for
  Simulation of Non-stationary Compressible Multicomponent Reactive Flows}.
\newblock {\em Journal of Scientific Computing}, 95(2):44, 2023.

\bibitem{rannabauer2018ader}
L.~Rannabauer, M.~Dumbser, and M.~Bader.
\newblock {ADER-DG} with a-posteriori finite-volume limiting to simulate
  tsunamis in a parallel adaptive mesh refinement framework.
\newblock {\em Computers \& Fluids}, 173:299--306, 2018.

\bibitem{rio2024high}
L.~R{\'\i}o-Mart{\'\i}n and M.~Dumbser.
\newblock {High-order ADER Discontinuous Galerkin schemes for a symmetric
  hyperbolic model of compressible barotropic two-fluid flows}.
\newblock {\em Communications on Applied Mathematics and Computation},
  6(4):2119--2154, 2024.

\bibitem{Rusanov:1961a}
V.~V. Rusanov.
\newblock {Calculation of Interaction of Non--Steady Shock Waves with
  Obstacles}.
\newblock {\em J. Comput. Math. Phys. USSR}, 1:267--279, 1961.

\bibitem{ShashkovCellCentered}
S.~Sambasivan, M.~Shashkov, and D.~Burton.
\newblock {A finite volume cell-centered Lagrangian hydrodynamics approach for
  solids in general unstructured grids}.
\newblock {\em International Journal for Numerical Methods in Fluids},
  72:770--810, 2013.

\bibitem{scovazzi2}
G.~Scovazzi.
\newblock {Lagrangian shock hydrodynamics on tetrahedral meshes: A stable and
  accurate variational multiscale approach}.
\newblock {\em Journal of Computational Physics}, 231:8029--8069, 2012.

\bibitem{springel2010pur}
V.~Springel.
\newblock {E pur si muove: Galilean-invariant cosmological hydrodynamical
  simulations on a moving mesh}.
\newblock {\em Monthly Notices of the Royal Astronomical Society},
  401(2):791--851, 2010.

\bibitem{stroud}
A.~Stroud.
\newblock {\em {Approximate Calculation of Multiple Integrals}}.
\newblock Prentice-Hall Inc., Englewood Cliffs, New Jersey, 1971.

\bibitem{tavellihigh}
M.~Tavelli and W.~Boscheri.
\newblock {A high order parallel Eulerian-Lagrangian algorithm for
  advection-diffusion problems on unstructured meshes}.
\newblock {\em International Journal for Numerical Methods in Fluids}, 2019.
\newblock in press.

\bibitem{tavelli2020space}
M.~Tavelli, S.~Chiocchetti, E.~Romenski, A.-A. Gabriel, and M.~Dumbser.
\newblock {Space-time adaptive ADER discontinuous Galerkin schemes for
  nonlinear hyperelasticity with material failure}.
\newblock {\em Journal of computational physics}, 422:109758, 2020.

\bibitem{toro3}
V.~Titarev and E.~Toro.
\newblock {ADER}: Arbitrary high order {Godunov} approach.
\newblock {\em Journal of Scientific Computing}, 17(1-4):609--618, December
  2002.

\bibitem{toro1}
E.~Toro, R.~Millington, and L.~Nejad.
\newblock Towards very high order {Godunov} schemes.
\newblock In E.~Toro, editor, {\em Godunov Methods. Theory and Applications},
  pages 905--938. Kluwer/Plenum Academic Publishers, 2001.

\bibitem{toro2005ader}
E.~Toro and V.~Titarev.
\newblock {ADER schemes for scalar non-linear hyperbolic conservation laws with
  source terms in three-space dimensions}.
\newblock {\em Journal of Computational Physics}, 202(1):196--215, 2005.

\bibitem{trulio1966air}
J.~Trulio.
\newblock Air force weapons laboratory.
\newblock {\em AFWL-Tr-66-19, June}, 1966.

\bibitem{vargas2025multi}
A.~Vargas, V.~Z. Tomov, M.~A. Skinner, V.~Dobrev, J.~Nikl, T.~Kolev, and R.~N.
  Rieben.
\newblock {Multi-material ALE remap with interface sharpening using high-order
  matrix-free finite element methods}.
\newblock {\em Journal of Computational Physics}, page 114367, 2025.

\bibitem{Neumann1950}
J.~von Neumann and R.~Richtmyer.
\newblock A method for the calculation of hydrodynamics shocks.
\newblock {\em Journal of Applied Physics}, 21:232--237, 1950.

\bibitem{wilkins1964methods}
M.~Wilkins, B.~Alder, S.~Fernbach, and M.~Rotenberg.
\newblock Methods in computational physics.
\newblock {\em Calculation of elastic--plastic flow}, pages 211--263, 1964.

\bibitem{Wilkins1964}
M.~L. Wilkins.
\newblock {Calculation of elastic-plastic flow}.
\newblock {\em Methods in Computational Physics}, 3, 1964.

\bibitem{wu2021cell}
W.~Wu, A.-M. Zhang, and M.~Liu.
\newblock {A cell-centered indirect Arbitrary-Lagrangian-Eulerian discontinuous
  Galerkin scheme on moving unstructured triangular meshes with topological
  adaptability}.
\newblock {\em Journal of Computational Physics}, 438:110368, 2021.

\bibitem{yuan2022hybrid}
D.~Yuan, P.~Tsoutsanis, and K.~Jenkins.
\newblock Hybrid high-order finite volume discontinuous {Galerkin} methods for
  turbulent flows.
\newblock In {\em World Congress in Computational Mechanics and ECCOMAS
  Congress}, 2022.

\bibitem{zanotti2025new}
O.~Zanotti, M.~Dumbser, D.~Balsara, and D.~Bhoriya.
\newblock {A new first-order formulation of the {Einstein} equations:
  comparison among different high order numerical schemes}.
\newblock In {\em Journal of Physics: Conference Series}, volume 2997, page
  012015. IOP Publishing, 2025.

\end{thebibliography}

\end{document}